\documentclass[regno]{elsarticle}

\usepackage{lineno,hyperref}
\usepackage{graphicx}
\usepackage{wrapfig,lipsum}
\usepackage{{ntheorem}}
\usepackage[T1]{fontenc}
\usepackage{appendix}
\usepackage{amssymb}
\usepackage{amsfonts}
\usepackage{multicol}
\usepackage{anyfontsize}
\usepackage[utf8]{inputenc}
\usepackage{calligra}
\usepackage[x11names,table]{xcolor}
\usepackage{multirow}
\usepackage{pdfpages}
\usepackage{nicefrac}
\usepackage{subcaption}
\usepackage{amsmath}
\usepackage{natbib}
\usepackage{esint}
\usepackage{algorithm}
\usepackage{algorithmic}
\usepackage{comment}
\modulolinenumbers[5]

\journal{Journal of \LaTeX\ Templates}

\newcommand{\ha}{\frac{1}{2}}
\newcommand{\R}{\mathbb{R}}

\newcommand{\F}{\mathcal{F}}
\newcommand{\tF}{\Tilde{F}}
\renewcommand{\epsilon}{\varepsilon}
\definecolor{Ema}{RGB}{0, 140, 250}

\newcounter{dctr}[section]
\numberwithin{equation}{section}

\newtheorem{remark}{Remark}[section]

\begin{document}
\begin{frontmatter}
\title{Semi-Implicit-type {order-adaptive} CAT2 schemes for systems of balance laws with relaxed source term}


\author[catania]{Emanuele Macca$^{*}$}
\ead{emanuele.macca@unict.it}
\cortext[cor1]{Corresponding author}

\author[catania]{Sebastiano Boscarino}
\ead{boscarino@dmi.unict.it}


\address[catania]{Dipartimento di Matematica ed Informatica Universit{\'a} di Catania, Viale Andrea Doria 6, 95125, Catania, Italy}


\begin{abstract}
{In this paper we present two semi-implicit-type second order Compact Approximate Taylor (CAT2) numerical schemes and blend them with a local \textit{a posteriori} Multi-dimensional Optimal Order Detection (MOOD) paradigm to solve hyperbolic systems of balance laws with relaxed source term. The resulting scheme presents high accuracy when applied to smooth solutions, essentially non-oscillatory behavior for irregular ones, and offers a nearly fail-safe property in terms of ensuring positivity. The numerical results obtained from a variety of test cases, including smooth and non-smooth well-prepared and unprepared initial condition, assessing the appropriate behavior of the semi-implicit-type second order CATMOOD schemes. These results have been compared in accuracy and efficiency with a second order semi-implicit Runge-Kutta (RK) method.}

\end{abstract}



\end{frontmatter}

\section{Introduction}
In the late 1960s, Lax and Wendroff introduced a numerical scheme based on divided differences for the hyperbolic resolution of systems of conservation laws \cite{LaxWendroff,LeVeque2007,Toro2009}. Over the years, this scheme has been extensively used in the field of numerical methods, providing valuable solutions for various conservation law systems. However, extending this technique to nonlinear conservation laws has been a longstanding challenge.

In the pursuit of extending the Lax-Wendroff scheme to nonlinear systems, significant efforts have been made by researchers in the past. Previously, with the 2-step schemes introduced by Richtmeyer et al. ~\cite{Richtmeyer} and MacCormack ~\cite{MacCormack}, successively by Qui and Shu, in 2003 \cite{Qiu-Shu}, attempted to achieve this extension by exploring the Cauchy-Kowalesky identity. During the same period, Toro and collaborators developed ADER methods \cite{Titarev,Schwartzkopff}, which offered higher-order solutions for the scheme.

Subsequent progress led to the work of Zorio et al. ~in 2017 \cite{Zorio}, who devised a numerical approach to circumvent the use of the Cauchy-Kowalesky identity. However, it's important to note that this method was not an exact extension of the original Lax-Wendroff scheme.

A noteworthy advancement came in 2019 when Carrillo et al. ~proposed a compact version of AT schemes, which served as a proper extension of the Lax-Wendroff scheme \cite{Carrillo-Pares}. These Compact Approximate Taylor (CAT) schemes exhibited remarkable performance when dealing with systems of conservation laws.

Despite these efforts, the extension of the Lax-Wendroff scheme to nonlinear conservation laws remained a challenge. In light of this, the present article aims to introduce two {\em semi-implicit-type} variants of the second order CAT2 \cite{CPZMR2020} scheme with the MOOD approach \cite{CDL1,CDL2,CATMOOD_2}. These variants not only overcome the limitations of the original Lax-Wendroff scheme but also offer an efficient and accurate solution to systems of balance laws, incorporating relaxation terms in the source. 

We are aware that when applying an explicit schemes to these systems, the time step size is limited by the CFL condition. This condition requires that the time step is proportional to the grid spacing and the maximum characteristic speed of the system. However, this requirement can be excessively restrictive, particularly in scenarios involving relaxation terms where may lead to excessively small time steps \cite{MaccaExner}. Consequently, explicit schemes can become computationally expensive and impractical for such cases.

By introducing these semi-implicit-type variants, this article not only addresses the extension of the Lax-Wendroff scheme to nonlinear systems but also emphasizes the practical advantage of circumventing the restrictive CFL condition in the presence of relaxation terms. The ability to use larger time steps in numerical simulations provides significant computational benefits, allowing for faster and more efficient solutions for systems of balance laws with relaxation terms.

In this paper, we compare these two semi-implicit-type CAT2 schemes with a semi-implicit method designed within an IMEX Runge-Kutta (RK) framework, \cite{Boscarino-Filbet,MaccaRussoBumi}. 

Finally, in the context of balance laws with a stiff source term, characterized by a parameter $\epsilon$, the asymptotic preserving (AP) property represents a crucial point for a numerical scheme \cite{BoscarinoRusso, ImexBosca, PareschiRusso}. The formal definition of AP schemes was introduced first by S. Jin, \cite{jin1995runge, jin2010asymptotic} and roughly speaking, the AP property ensures that the numerical scheme accurately captures the correct limit as $\epsilon \to 0$, while keeping $\Delta t$ and $\Delta x$ (time and space discretizations, respectively) fixed.  Therefore an additional focus of this paper is to present numerical schemes within a semi-implicit approach \cite{Boscarino-Filbet} that extend the Lax-Wendroff method and respect the AP property when applied to systems of  balance laws with stiff source terms.

In the subsequent sections, we present the detailed formulation and implementation of the proposed {\em semi-implicit-type} variants, highlighting their capability to overcome the CFL limitation and their enhanced accuracy and versatility compared to a traditional semi-implicit RK method \cite{Boscarino-Filbet}. Through numerical experiments and comparisons, we demonstrate the practical benefits of these variants in accurately resolving nonlinear systems with relaxation terms.

The rest of this paper is organized as follows. 
The second section introduces the governing equations.
Next, we derive second order semi-implicit-type variants of CAT and CATMOOD ans we present a second order semi-implicit Runge-Kutta scheme for comparison with these two variants.
The fourth section presents the models of system ob balance laws with stiff source terms considered in this work. Section 5 is related to the stability analysis of the two semi-implicit-type scheme. Finally, numerical results are gathered in the sixth section to assess the good behavior of the semi-implicit-type CAT, CATMOOD compared to second order semi-implicit Runge-Kutta methods.
Conclusions and perspectives are finally drawn.

\section{Governing equations} \label{sec:gov_equ}

\subsection{Non-linear scalar conservation laws}
To simplify the description of the numerical methods we consider the non-linear scalar conservation law on the $Oxt$-Cartesian frame
\begin{equation}
    \label{sec:CAT_gov_equ}
    \partial_t u + \partial_x f(u) = 0,
\end{equation} 
where $u=u(x,t):\mathbb{R}\times\mathbb{R}^+\rightarrow\mathbb{R}$ denotes the scalar variable, and $f(u)=f(u(x,t))$ the non-linear flux depending on $u$.
$u(x,0)=u_0(x)$ denotes the initial condition, while the boundary conditions are prescribed depending on the test case, for instance periodic ones, Dirichlet or Neumann ones. 
Eq ~(\ref{sec:CAT_gov_equ}) represents the generic model of non-linear scalar equation, the simplest one being probably Burgers' equation for which the flux is given by: $f(u)=u^2/2$.
~(\ref{sec:CAT_gov_equ}) also represents the generic model of linear scalar advection equation if $f(u)=a u$ with $a\in \mathbb{R}$ being the advection velocity. \\ 
In the following we denote the partial derivative in time and space with under-script letters as 
$u_t \equiv \partial_t u$ and $u_x \equiv \partial_x u$.

\subsection{1D system of balance laws with relaxation term}
This paper deals with the design of high-order  methods for hyperbolic systems of balance laws with relaxation term
\begin{equation}
    \label{hyp_law}
    U_t + F(U)_x = S(U,\varepsilon)
\end{equation}
with initial condition $U(x,0) = U_0(x),$ where $U:\R\times[0,+\infty] \rightarrow\R^d$ is the unknown vector field, $F:\R^d\rightarrow\R^d$ is the flux function, and $S:\R^d\times]0,1]\rightarrow\R^d$ is the source term which contains stiff parameter $\varepsilon.$ PDE systems of this form appear in many fluid models in different contexts: multi-phase flow models, traffic flow models, gas dynamic, shallow water equations, etc.

Before introducing the numerical schemes in this paper, we define the time domain as $\mathcal{T}=[0,T]$, where $T>0$, and divide it into time intervals $[t^n,t^{n+1}]$, $n\in \mathbb{N}$, with time-steps $\Delta t=t^{n+1}-t^n$, subject to a CFL (Courant-Friedrichs-Lewy) condition\footnote{{Although it is better to assign $\Delta t$ dynamically at each time step by imposing some CFL condition, here we shall adopt a constant time step $\Delta t$, in order to simplify the notation in the description of the method.}} \cite{CFL}. The computational spatial domain, denoted as $\Omega$, represents a 1D segment divided into $N_x$ cells. The generic cell is denoted $\omega_i$ and indexed by a unique label $1\leq i \leq N_x$. Classically we identify the cell end-points by half indexes so that $\omega_i=[x_{i-1/2},x_{i+1/2}]$ and the cell center is given by $x_i=\frac12(x_{i+1/2}+x_{i-1/2})$.
The size of a cell is given by $\Delta x_i = x_{i+1/2}-x_{i-1/2}$, that we simply denote as $\Delta x$, since we assume uniform grid.

\section{Numerical schemes}\label{sec:schemes}
In the forthcoming section, we will introduce CAT2 schemes for systems of conservation laws \eqref{sec:CAT_gov_equ}, \cite{TesiPhD} and two semi-implicit-type variants of the CAT2 schemes for system of balance laws with relaxation term \eqref{hyp_law}. Through a detailed presentation of the mathematical formulations and numerical techniques, we will explore the advantages and effectiveness of these semi-implicit-type variants in improving accuracy and expanding the applicability of the CAT2 scheme.

Furthermore, we will show that that the advantage of semi-implicit-type numerical methods applied to the system of the form (\ref{hyp_law}) offers advantages for accuracy, stability, robustness, asymptotic preserving proprieties, and computational efficiency. These advantages make this approach a valuable tool for investigating the behavior of rarefied gases with relaxation effects and analyzing the impact of relaxation processes on models of balance law with relaxation terms \eqref{hyp_law} for a wide range of parameter $\epsilon$, ensuring accurate solutions for different degrees of relaxation.

\subsection{CAT2 scheme for conservation laws}
First of all, let us recap the explicit expression of CAT2 method for a 1D system of conservation laws with initial condition $U(x,0) = U_0(x)$
\begin{equation}
    \label{sys_con}
    U_t + F(U)_x = 0.
\end{equation}
The formal expression, written in conservative form, of the second order CAT2 scheme is 
\begin{equation}
    \label{Cat2_formal}
    U_i^{n+1} = U_i^n - \frac{\Delta t}{\Delta x}\Bigl(F_{i+\ha} - F_{i-\ha}\Bigr),
\end{equation}
where 
\begin{equation}
    \label{F_CAT_2}
    F_{i+\ha} = \ha\Bigl(F_i^{(0)} + F_{i+1}^{(0)}\Bigr) + \frac{\Delta t}{4}\Bigl(F_{i,0}^{(1)} + F_{i,1}^{(1)}\Bigr).  
\end{equation}
The time derivatives of flux, for each position $s=0,1,$ are computed with a Taylor expansion in time as follows
\begin{align}
    F_{i+s}^{(0)} &= F(U_{i+s}^n) \\ 
    F_{i,s}^{(1)} &= \frac{1}{\Delta t}\Bigl(F(U_{i,s}^{1,n+1})-F(U_{i+s}^n)\Bigr),
\end{align}
in which 
\begin{equation}
    \label{U_tay_1}
    U_{i,s}^{1,n+1} = U_{i+s}^n + \Delta t U_{i,s}^{(1)}
\end{equation}
and {the first time derivatives}  $U_{i,s}^{(1)} $ are computed throughout the numerical version of the Cauchy-Kowalesky identity
\begin{equation}
    U_{i,s}^{(1)} = -\frac{1}{\Delta x}\Bigl(F(U_{i+1}^n) - F(U_i^n)\Bigr).
\end{equation}
Finally, the second order flux reconstruction is obtained as:
\begin{equation}
    \label{F2_fast}
    F_{i+\ha} = \frac{1}{4}\Bigl(F(U_{i}^n)+F(U_{i+1}^n) + F(U_{i,0}^{1,n+1})+F(U_{i,1}^{1,n+1})\Bigr).
\end{equation}

\subsection{{Semi-implicit-type CAT2 schemes}} 
In this section, we will introduce two semi-implicit-type versions of the CAT2 scheme for systems of balance laws \eqref{B_W_relax} each possessing distinct stability properties (see Section \ref{sec:ode_like_stability}).

As it is known, the CAT2 scheme is a spatial and temporal method that operates in a single step. This characteristic means that it does not differentiate between space and time reconstructions. However, we can discuss the reconstruction of the flux $F_{i+1/2}$ in the context of second-order applications for systems of balance laws (\ref{hyp_law}).

Let us consider the 1D systems of balance laws with relaxation term
\begin{equation}
    \label{Bal_laws}
    U_t + F(U)_x = S(U,\epsilon).
\end{equation}
A preliminary representation of the {semi-implicit-type} CAT2 scheme can be obtained by:
\begin{equation}
    \label{CAT2_preliminary}
    U^{n+1}_i = U^n_i + \frac{\Delta t}{\Delta x}\Bigl(\F_{i+\ha} - \F_{i-\ha}\Bigr) + \Delta t\mathcal{S}(U^{n+1}_{i}, U^{n}_i, \epsilon)
\end{equation}
in which $\mathcal{S}(U^{n+1}_{i},, U^{n}_i, \epsilon)$ is a suitable approximation of second-order of the relaxed source term  $S(U,\epsilon).$ 

The flux function $\F_{i+\ha}$ is obtained applying an extension of the CAT2 method presented before. For this reasoning, let us consider the flux \eqref{F2_fast} we need to change the expression of $U_{i+s}^{1,n+1}$ in a $\tilde{U}_{i+s}^{1,n+1}$ for all $s = 0,1.$ Indeed, we have to add the source term to \eqref{U_tay_1} in the following implicit way
\begin{equation}
    \label{U_tay_1_imp}
    \tilde{U}_{i,s}^{1,n+1} = U_{i+s}^n + \Delta t U_{i,s}^{(1)} + \Delta t S(\tilde{U}_{i,s}^{1,n+1}).
\end{equation}
Then, the second order flux reconstruction is obtained as:
\begin{equation}
    \label{F2_fast_bal}
    \F_{i+\ha} = \frac{1}{4}\Bigl(F(U_{i}^n)+F(U_{i+1}^n) + F(\tilde{U}_{i,0}^{1,n+1})+F(\tilde{U}_{i,1}^{1,n+1})\Bigr).
\end{equation}
\begin{remark}
{Note that in order to maintain simplicity, we kept the notation $U_{i,s}^{(1)}$ in (\ref{U_tay_1_imp})  as an approximation of the first time derivative $U$. However, it is important to note that while this notation holds true for systems of conservation laws, it is not applicable to systems of balance laws (\ref{Bal_laws}). This is due to the fact that in balance laws, the equation takes the form $u_t = -f(u)_x + S(u)$.}
\end{remark} 

\paragraph{{Cat2-Trap scheme}}
Let consider the {semi-implicit-type method} applied to system of balance laws \eqref{Bal_laws} in which the relaxed source term is treated with the implicit-explicit trapezoidal rule
\begin{equation}
    \label{Met_trap_rule}
    U^{n+1}_i = U^n_i + \frac{\Delta t}{\Delta x}\Bigl(\F_{i+\ha} - \F_{i-\ha}\Bigr) + \frac{\Delta t}{2}\Bigl(S^{n+1}_{i,\epsilon} + S^n_{i,\epsilon}\Bigr)
\end{equation}
where $S^{n+1}_{i,\epsilon}$ and $S^{n}_{i,\epsilon}$ represent, respectively, $S(U^{n+1}_{i},\epsilon)$ and $S(U^n_{i},\epsilon)$ and $\F_{i+\ha}$ is obtained by \eqref{F2_fast_bal}. We refer to this scheme (\ref{Met_trap_rule}) as the CAT2-Trap.

\paragraph{{CAT2-Tay scheme}}
{Now, in (\ref{Met_trap_rule}), we adopt a different approximation for the source term $S(U, \epsilon)$, resulting in}
\begin{equation}
    \label{Met_L-stalbe}
    U^{n+1}_i = U^n_i - \frac{\Delta t}{\Delta x}\Bigl(\F_{i+\ha} - \F_{i-\ha}\Bigr) + \frac{\Delta t}{2}\Bigl(S^{n+1}_{i,\epsilon} + S^{n+1}_{i,\epsilon} - \Delta t\nabla S^n_{i,\epsilon}(S^{n+1}_{i,\epsilon}-D_x(F^n))\Bigr)
\end{equation}
where $S^n_{i,\epsilon}$ of the trapezoidal rule was replaced by
$$ S^n_{i,\epsilon} = S^{n+1}_{i,\epsilon} - \Delta t\nabla S^n_{i,\epsilon}(S^{n+1}_{i,\epsilon}-D_x(F^n))$$ while $\F_{i+\ha}$ is exactly the same flux reconstruction of the previous scheme \eqref{Met_trap_rule} and $D_x(F^n)$ has been approximated with central derivatives as:
$$D_x(F^n) \approx \frac{1}{2\Delta x}\Bigl(F(U_{i+1}^n) - F(U_{i-1}^n)\Bigr).$$ 
Furthermore, $\nabla S^n_{i,\epsilon}$ represents $\partial S /\partial U$ valuated at time $t_n$ and position $x_i.$
{We refer to this scheme (\ref{Met_L-stalbe}) as the CAT2-Tay scheme, and in Section \ref{sec:ode_like_stability}, we will show that it exhibits better stability properties compared to CAT2-trap.}

\subsection{CAT-MOOD}
\label{sec:CAT-MOOD}
The main aim of this section is to combine the \textit{a posteriori} shock capturing technique (MOOD) \cite{CDL1,CATMOOD_1D} with an high-order semi-implicit-type CAT schemes, which enables simultaneous spatial and temporal reconstructions.

In addition to preserving the desired properties of the exact solutions and mitigating numerical oscillations, the MOOD algorithm also plays a crucial role in enhancing the robustness of the high-order CAT methods. By dynamically adapting the order of accuracy based on the admissibility criteria, the MOOD algorithm ensures that the numerical solution remains reliable and stable, even in the presence of complex flow phenomena. This adaptability allows for a seamless transition between high-order accuracy in smooth regions and lower-order accuracy in regions with discontinuities or strong gradients. As a result, the CATMOOD schemes provide an effective balance between accuracy and stability, making them suitable for a wide range of applications in computational fluid dynamics, including simulations involving shocks, boundary layers, and other challenging flow features. Furthermore, this procedure helps reduce and eliminate numerical oscillations that are typically introduced by high-order methods when dealing with shocks or large gradients, as observed in CAT methods \cite{Carrillo-Pares,MCPR2022}.

The core concept underlying the MOOD procedure involves the application of a high-order method across the entire computational domain for a given time step. Subsequently, the behavior of the solution is assessed locally in each cell $i$ using admissibility criteria. If the solution computed in cell $i$ at time $t^{n+1}$ satisfies the prescribed criteria, it is retained. However, if it fails to meet the requirements, it undergoes recalculation using a numerical method of lower order. This iterative process persists until an acceptable solution is obtained or, in cases where the admissibility criteria are not fulfilled in any of the previous reconstructions, the lowest order scheme (first order) is employed.

Therefore, the primary aim of this research is to construct a hierarchy of CAT methods in which the local order of accuracy is dynamically adapted using \textit{a posteriori} admissibility criteria. This approach gives rise to a novel class of adaptive CAT methods referred to as CATMOOD schemes \cite{CATMOOD_1D,CATMOOD_2}. By incorporating the MOOD algorithm and its admissibility criteria, the CATMOOD schemes offer the flexibility to adjust the order of accuracy based on the specific characteristics and requirements of the solution.

\subsubsection{MOOD Admissibility Criteria} \label{ssec:MOODdetection}
In this investigation, we consider three specific admissibility criteria \cite{CDL2,CDL3} to evaluate the candidate numerical solution $\left\{ u_i^{n+1} \right\}_{1\leq i \leq I}$:
\begin{enumerate}
    \item \textit{Computer Admissible Detector (CAD)}: This criterion identifies undefined or unrepresentative quantities, such as \texttt{NaN} (not-a-number) or infinite values, which typically arise from division by zero or similar situations.
    \item \textit{Physical Admissible Detector (PAD)}: The second detector ensures the physical validity of the candidate solution. This criterion is crucial to avoid non-physical sound speeds, imaginary time steps, and other undesired effects. It is important to note that the assessment of physicality is dependent on the specific model of the solved PDEs and should be tailored accordingly.
    \item \textit{Numerical Admissible Detector (NAD)}: This criterion represents a relaxed version of the Discrete Maximum Principle, ensuring the essential non-oscillatory (ENO) nature of the solution \cite{Shu1997,Ciarlet}. It verifies that the numerical solution $w_i^*$ obtained using a scheme of order $2P$ satisfies the following inequality:
    $$\min_{c\in \mathcal{C}_i^P}(w^n_c) - \delta_i \le w^*_i \le \max_{c\in \mathcal{C}_i^P}(w^n_c) + \delta_i,$$
    where $\mathcal{C}_i^P = \{-P,\ldots,P\}$ represents the local centered stencil of order $2P$, and $\delta_i$ denotes a relaxation term introduced to mitigate issues in regions with minimal variations. The value of $\delta_i$ is determined as follows:
    \begin{equation}
        \label{eq:delta}
        \delta_i = \max\Big(\varepsilon_1,\varepsilon_2\left[\max_{c\in \mathcal{C}_i^P}w^n_c - \min_{c\in \mathcal{C}_i^P}w^n_c\right]\Big),
    \end{equation}
    where $\varepsilon_1$ and $\varepsilon_2$ represent small dimensional constants. The $(NAD)$ criterion guarantees the absence of significant spurious minima or maxima in the local solution. 
\end{enumerate}
The iterative nature of the MOOD algorithm ensures that the numerical solution maintains desirable properties while mitigating potential issues caused by numerical oscillations. By promptly detecting and addressing problems such as \texttt{NaN} values and violation of positivity requirements, the algorithm promotes the reliability and stability of the computations. Additionally, the dynamic adjustment of the local order of accuracy through the use of an auxiliary scheme with reduced precision allows for a more robust and adaptable approach to solving complex problems. This integration of admissibility criteria within the MOOD framework enhances the overall performance of the CATMOOD schemes, making them effective tools for a wide range of applications in scientific and engineering domains.

\paragraph{CAT scheme with MOOD limiting} 
The primary goal of this research is to attain a second order of accuracy in regions of the computational domain where a smooth solution is present. To achieve this, our focus is on utilizing the CAT2 scheme as extensively as possible. In contrast, for cells that contain a discontinuous solution, we intend to employ a robust but less accurate first order scheme, such as Rusanov or HLL fluxes, only when and where it is required. By employing a combination of these schemes, we aim to strike a balance between accuracy and robustness, ensuring efficient and reliable numerical computations across a wide range of scenarios.

\subsection{Second order semi-implicit 2-stages IMEX Runge-Kutta (RK) scheme.}

In this section, we introduce a second-order semi-implicit IMEX-RK scheme with 2 stages applied to system \eqref{Bal_laws}, which will be used for comparison with the two previously introduced CAT2 schemes.

Following the idea proposed in \cite{Boscarino-Filbet}, the system \eqref{Bal_laws} can be reformulated as a larger system of ordinary differential equations (ODEs) by utilizing appropriate discrete operators for approximating the spatial derivatives. The key concept in \cite{Boscarino-Filbet} is to determine which specific terms should be treated implicitly and which can be handled explicitly.

In accordance with \cite{Boscarino-Filbet}, we rewrite the system \eqref{Bal_laws} by identifying the terms that will be treated explicitly and those that require an implicit treatment. This approach involves doubling system \eqref{Bal_laws}, i.e., 
\begin{equation}
        \label{ode_bosca_semi}
        U'_E = {K}(U_E,U_I) \quad \quad U'_I = {K}(U_E,U_I),
\end{equation}
where ${K}(U_E,U_I)$ is given by 
\begin{equation}
        \label{ode_form_semi}
        {K}(U_E,U_I) = 
        \begin{bmatrix}
            -D_x(F(U_E)) & S(U_I,\epsilon)
        \end{bmatrix},
\end{equation} 
in which for the spacial discretization $$D_x(F_E)|_{x=x_i} = \frac{\tF_{i+1/2} - \tF_{i-1/2}}{\Delta x},$$ 
with  $\tF_{i+1/2} = \ha\Bigl( F(U_{i+1/2}^{-}) + F(U_{i+1/2}^{+}) - \alpha_{i+1/2}\big(U_{i+1/2}^{+} - U_{i+1/2}^{-}\big)\Bigr)$ is the Rusanov flux and $\alpha_{i+1/2}$ is related to the eigenvalues of the explicit hyperbolic part. $U_{i\pm1/2}^{\pm}$ are obtained by piecewise linear reconstruction with  MinMod slope limiter.

Now to obtain a semi-implicit time discretization, we apply an IMEX Runge-Kutta method for partitioned system to the resulting double system (\ref{ode_bosca_semi}).

A semi-implicit RK scheme is defined by a double Butcher tableau of the form 
\begin{equation*}
\label{IMEX_tableau}
    \begin{array}{c|cc}
            \tilde{c} & \tilde{A} \\ \hline \\
                      & \tilde{b}^\top
        \end{array} 
        \hspace{3 cm}
    \begin{array}{c|cc}
            {c} & {A} \\ \hline \\
                      & {b}^\top
        \end{array} .
\end{equation*}
The explicit part of the scheme is characterized by the lower triangular matrix $\tilde{A}$, the vector $\tilde{c}$, and the vector $\tilde{b}$, while the implicit part is described by the triangular matrix $A$, the vector $c$, and the vector $b$.

Applying the IMEX-RK scheme to the system \eqref{ode_bosca_semi}, and assuming that $\tilde{b}_i = b_i$ $i = 1,...,2$, for we obtain the following steps:
\begin{enumerate}
    \item  Compute stage values for $i=1,\ldots,s$:
   \begin{align*}
   U^{(i)}_E & = U^n + \Delta t\sum_{j=1}^{i-1}\tilde{a}_{i,j} K\left(U_E^{(j)},U_I^{(j)}\right), \\
   U^{(i)}_I & = U^n + \Delta t\left(\sum_{j=1}^{i-1}a_{i,j} K\left(U_E^{(j)},U_I^{(j)}\right) + a_{i,i} K\left(U_E^{(i)},U_I^{(i)}\right)\right).
   \end{align*}
   \item Compute the numerical solution:
    \begin{align*}
   U^{n+1} & = U^n + \Delta t\sum_{i=1}^{s}b_{i} K\left(U_E^{(j)},U_I^{(j)}\right),
   \end{align*}
   \item If the implicit part of the scheme is stiffly accurate \cite{wanner1996solving}, i.e., $a_{si} = b_i$, for $i=1,..,s$ we get for the numerical solution:
   $$U^{n+1} = U_I^{(s)}$$.
\end{enumerate}

For a second-order semi-implicit time discretization, we consider the IMEX scheme defined by the following double Butcher tableau \cite{Boscarino-Filbet}:
\begin{equation*}
\begin{array}{c|cc}
        & 0 &  \\
        c & c & 0\\ \hline
        & 1-\gamma & \gamma
    \end{array} 
    \hspace{3 cm}
    \begin{array}{c|cc}
        \gamma & \gamma &  \\
        1 & 1-\gamma & \gamma\\ \hline
        & 1-\gamma & \gamma
    \end{array}
\end{equation*}
where $\displaystyle \gamma = 1 - \frac{1}{\sqrt{2}}$ and $\displaystyle c = \frac{1}{2\gamma}$.

Applying this scheme to the system, we can proceed as follows:
\begin{enumerate}
    \item $U_E^{(1)} = U^n;$
    \item $U_I^{(1)} = U^n + \Delta t\gamma H(U_E^{(1)},U_I^{(1)});$
    \item $U_E^{(2)} =  U^n + \Delta tc H(U_E^{(1)},U_I^{(1)});$
    \item $U_I^{(2)} = U^n + \Delta t(1-\gamma) H(U_E^{(1)},U_I^{(1)}) + \Delta t\gamma H(U_E^{(2)},U_I^{(2)});$
    \item $U^{n+1} = U_I^{(2)}.$
\end{enumerate}

\begin{remark}    
Note that steps 3 and 4, can be rewritten as follows:
\begin{align*}
U_E^{(2)} &= (1-\frac{c}{\gamma})U^n + \frac{c}{\gamma}U_I^{(1)},  \\
U_I^{(2)} &= (1-\frac{1-\gamma}{\gamma})U^n + \frac{1-\gamma}{\gamma}U_I^{(1)} + \Delta t\gamma K(U_E^{(2)},U_I^{(2)}).
\end{align*}
For further details on this method, please refer to \cite{Boscarino-Filbet}.
\end{remark}
We will refer this scheme in the numerical tests section as Runge-Kutta2.

\section{Relaxed model} \label{sec:Models}
This section is dedicated to the presentation of models of balance laws with relaxation terms. Specifically, we will study the Xin Jin model \cite{jin1995runge}, which focuses on the mathematical modeling of fluid flows with relaxation effects. The Broadwell model \cite{PareschiRusso, broadwell1964shock}, on the other hand, is a mathematical framework that incorporates relaxation processes into the study of rarefied gas dynamics. Lastly, the Euler system with heat transfer \cite{jin1995runge} examines the behavior of fluids in the presence of heat exchange, considering the conservation of mass, momentum, and energy. These models provide valuable insights and mathematical tools for understanding and analyzing complex physical phenomena involving fluid dynamics and energy transfer.

\subsection{Xin Jin}
The relaxation model developed by Xin Jin is a mathematical framework used to describe fluid flows with relaxation effects. This model incorporates the concept of relaxation time, which accounts for the time required for a fluid to reach a state of equilibrium after being subjected to an external disturbance. By considering the relaxation time in the governing equations of fluid dynamics, the model allows for a more accurate representation of the fluid behavior, especially in situations where the flow transitions from one state to another.

The Xin Jin model introduces relaxation terms that capture the gradual adjustment of the fluid properties, such as density, velocity, and temperature, towards their equilibrium values. These relaxation terms are formulated based on empirical observations and physical principles, and they play a crucial role in capturing the transient behavior and dynamics of the fluid system.

The model of relaxation proposed by Xin Jin has found applications in various fields, including computational fluid dynamics, chemical engineering, and environmental sciences. It provides a valuable tool for studying complex flow phenomena, such as shock waves, boundary layer transitions, and flow instabilities, where relaxation effects are significant and need to be taken into account for accurate predictions and analysis.

The 1D Xin-Jin relaxed model can be written as:
\begin{equation}
    \label{Xin_Jin_model}
    \begin{cases}
    u_t + v_x = 0; \\ v_t + u_x = \displaystyle -\frac{1}{\varepsilon}\Bigl(v-g(u)\Bigr)
    \end{cases}
\end{equation}
in which $\epsilon$ represents the relaxed parameter. 
System \eqref{Xin_Jin_model} can be written in matrix form as \eqref{Bal_laws} where the constitutive variable is  $U = [u,v]^T,$ the flux is $F(U) = [v,u]^T$ and the source term is $S(U,\epsilon) = [0,v-g(u)]^T/\epsilon.$
In system \eqref{Xin_Jin_model} when $\epsilon\rightarrow 0,$ $v = g(u)$ then system \eqref{Xin_Jin_model} relaxes to a 1D conservation law of the form:
\begin{equation}
    \label{X_J_relax} 
    u_t + g(u)_x = 0.
\end{equation}

\subsection{Broadwell model} 
Usually, the Broadwell model \cite{broadwell1964shock, PareschiRusso}  is used to study the behavior of conserved variables in the presence of relaxation processes. The conserved variables in the Broadwell model are represented by the vector \(U = (\rho, \rho v, z)\), where \(\rho\) denotes density, \(v\) represents velocity, and \(z\) is an additional state variable. 
The flux functions are $F(U) = [\rho v,z,\rho v]^T$ and the relaxed source terms are \(S(U,\epsilon) = [0,0,\rho^2 + (\rho v)^2 - 2\rho z]^T/2\epsilon\).

Explicitly, the Broadwell model becomes:
\begin{equation} \label{Broadwell_model}
\begin{cases}
\displaystyle \frac{\partial \rho}{\partial t} + \frac{\partial}{\partial x} (\rho v) = 0 \\
\displaystyle\frac{\partial}{\partial t} (\rho v) + \frac{\partial}{\partial x} (z) = 0 \\ 
\displaystyle \frac{\partial z}{\partial t} + \frac{\partial}{\partial x} \left( \rho v\right) = \frac{1}{2\epsilon}\left(\rho^2 + (\rho v)^2 - 2 \rho z\right)
\end{cases}
\end{equation}

 In  \eqref{Broadwell_model}  when $\epsilon \rightarrow 0,$ $z = (\rho + \rho v^2)/2$, system \eqref{Broadwell_model} relaxes to a 1D system of conservation laws of the form:
    \begin{equation}
        \label{B_W_relax}
        \begin{cases}
            \rho_t + (\rho v)_x = 0 \\
            (\rho v)_t + \ha\Bigl( \rho + \rho v^2\Bigr) = 0.
        \end{cases}
    \end{equation}

\subsection{Euler equations with heat transfer model}
The Euler equations with heat transfer are a fundamental set of equations used to describe the behavior of compressible fluids while accounting for the transfer of thermal energy. They are derived from the Euler equations, which are a simplified form of the Navier-Stokes equations, neglecting viscous and thermal diffusion effects.

In the Euler equations with heat transfer, an additional term is included to account for the exchange of thermal energy between the fluid and its surroundings. This heat transfer term incorporates Fourier's law of heat conduction, which states that the heat flux is proportional to the negative temperature gradient. These equations find wide-ranging applications in various fields of engineering and science. They are extensively used in aerodynamics, combustion analysis, astrophysics, and many other areas where compressible fluid flows and heat transfer play a crucial role. The Euler equations with heat transfer provide insights into the propagation of pressure and temperature waves, the behavior of shock waves, and the interaction between fluid flow and thermal effects.

The Euler equations with heat transfer in one dimension without gravity acceleration can be written as follows:
\begin{equation}
    \label{Euler_heat}
    \begin{cases}    
     \displaystyle   \frac{{\partial \rho}}{{\partial t}} + \frac{{\partial (\rho u)}}{{\partial x}} = 0 \\[2mm]
     \displaystyle   \frac{{\partial (\rho u)}}{{\partial t}} + \frac{{\partial (\rho u^2 + p)}}{{\partial x}} = 0 \\[2mm]
      \displaystyle  \frac{{\partial (\rho E)}}{{\partial t}} + \frac{{\partial (\rho u E + p u)}}{{\partial x}} = -K\rho(T-T_0)
    \end{cases}
\end{equation}
where \(\rho\) represents the density of the fluid,  \(u\) the velocity of the fluid, \(E = e + u^2/2\) the total energy per unit volume and $e$ the internal energy,  \(p\) the pressure, \(T\) the temperature, $K$ the thermal conductivity and $T_0$ is the temperature of the thermal bath. Away from the equilibrium the gas could be considered as $\gamma-$law gas, then $p=(\gamma-1)\rho e.$ The internal energy is proportional to the absolute temperature\footnote{From perfect gas law, $p = \mathcal{R}\rho T,$ where $\mathcal{R} = R/m$ is the universal gas constant $R$ divided by the gas molecular mass $m.$ For a polytropic gas it is $e = p/((\gamma-1)\rho) = \mathcal{R}T/(\gamma-1) = c_v T.$}, $e = \mathcal{R}T/(\gamma-1).$

By solving these equations, one can analyze and understand the complex dynamics of fluids, including the propagation of shock waves, the formation of boundary layers, and the behavior of thermal effects in various engineering and scientific applications.

\section{Numerical properties}
In this section, we will delve into some theoretical aspects of the methods described in the previous section. Specifically, we will explore some stability properties for the semi-implicit-type CAT2 schemes applied to a linear scalar equation. Additionally, we will conduct an analysis of linear stability using Fourier modes. 

\subsection{{Linear stability analysis}}
\label{sec:ode_like_stability}
{In this section, our goal is to develop a theoretical stability analysis, following a similar approach used for linear ordinary differential equation (ODE) systems, i.e., $u' = \lambda u$ with $\lambda$ positive parameter. To simplify the discussion, we will focus on studying a scalar linear balance law that includes a linear source term $S(u) = \lambda u$. With this in mind, let us consider}

\begin{equation}
    \label{equ:lin_bal_stiff}
    u_t + u_x = \lambda u,
\end{equation}
where $\lambda\in \R$ and for hypothesis we assume $\lambda \gg 0$. In this latter, the exact solution converges to zero. 

Now we analyze the behaviour of the numerical solution for each of the CAT2 scheme, \eqref{Met_trap_rule} and \eqref{Met_L-stalbe}.
Since for the exact solution in the limit case $\lambda \to \infty$, we have $u = 0$,  we will show which scheme achieves correctly the behaviour of the exact solution. 

The semi-implicit scheme CAT2-Trap \eqref{Met_trap_rule} applied to equation \eqref{equ:lin_bal_stiff} becomes:
\begin{equation}
    u_i^{n+1} = u_i^n - \frac{\Delta t}{\Delta x}\Bigl(\mathcal{F}_{i+\ha} - \mathcal{F}_{i-\ha}\Bigr) + \frac{\Delta t}{2}\Bigl(\lambda u_i^{n+1} + \lambda u_i^n\Bigr),
\end{equation}
and defining $z = \lambda\Delta t$ we obtain:
\begin{equation}
    \label{equ:a_stable_0}
    \Bigl(1-\frac{z}{2}\Bigr)u_i^{n+1} = \Bigl(1+\frac{z}{2}\Bigr)u_i^{n} - \frac{\Delta t}{\Delta x}\Bigl(\mathcal{F}_{i+\ha} - \mathcal{F}_{i-\ha}\Bigr). 
\end{equation}

    By construction, $\mathcal{F}_{i+1/2} - \mathcal{F}_{i-1/2}$ \eqref{F2_fast_bal} depends on $z.$ Indeed, 
    $$\mathcal{Q}(z) = \Bigl(\mathcal{F}_{i+\ha} - \mathcal{F}_{i-\ha}\Bigr) = \frac{1}{4}\Big(u_{i+1}^n - u_{i-1}^n + \frac{1}{1+z}\Bigl(u_{i+1}^n-u_{i-1}^n + \frac{2\Delta t}{\Delta x}\Bigl(2u_i^n - u_{i-1}^n - u_{i+1}^n\Bigr)\Bigr)\Big).$$ 
    Then, $\mathcal{Q}(z) = c_1 + c_2/(1+z)$ where $c_1$ and $c_2$ are constants. 

{Hence, we get the following numerical solution}
{\begin{equation}
    \label{equ:a_stable_func}
    u_i^{n+1} =\mathcal{R}(z) u_i^{n} - \frac{\frac{\Delta t}{\Delta x}\Bigl(c_1 + \frac{c_2}{1+z}\Bigr)}{1-\frac{z}{2}}
\end{equation}
with 
$$
    \mathcal{R}(z) = \frac{1+\frac{z}{2}}{1-\frac{z}{2}}.
$$
}
{Now setting $\Delta t$ and $\Delta x$ such that the CFL stability condition is satisfied \cite{CFL}, when $z\rightarrow\infty$, 
$$  u_i^{n+1} \rightarrow - u^n_i \neq 0.$$
Therefore, we can observe that in the limit, 
the scheme \eqref{Met_trap_rule} does not converge correctly to the limit solution, i.e., $u = 0$.}

On the contrary the semi-implicit scheme CAT2-Tay \eqref{Met_L-stalbe} applied to equation \eqref{equ:lin_bal_stiff} turns into 
\begin{equation}
    u_i^{n+1} = u_i^n - \frac{\Delta t}{\Delta x}\Bigl(\mathcal{F}_{i+\ha} - \mathcal{F}_{i-\ha}\Bigr) + \Delta t\Bigl(\lambda u_i^{n+1} -\frac{\Delta t}{2}\lambda(\lambda u_i^{n+1} - \frac{1}{2\Delta x}(u_{i+1}^n - u_{i-1}^n)   \Bigr).
\end{equation}
and similarly to the previous scheme we get:
{\begin{equation}
    \label{equ:l_stable_func}
    u_i^{n+1} =\mathcal{R}(z) u_i^{n},
\end{equation}
where
\begin{equation}
    \label{equ:l_stable_tay}
    \mathcal{R}(z) = \frac{1+\frac{z}{4}c_3 - \frac{\Delta t}{\Delta x}\Bigl(c_1 + \frac{c_2}{1+z}\Bigr)}{1-z+\frac{z^2}{2}}.
\end{equation}}
{Now setting $\Delta t$ and $\Delta x$ such that the CFL stability condition is satisfied \cite{CFL}, when $z\rightarrow\infty$, we get 
$ u^{n+1}_i \rightarrow 0,$, i.e., the correct limit solution.}

This simple analysis, in the context of balance laws with a source term, provides us with a  insight to understand the {\em asymptotic preserving} property behavior for the two semi-implicit-type CAT2 schemes. It's straightforward to show that in the limit case $\varepsilon \to 0$ the CAT2-Trap scheme \eqref{Met_trap_rule} performs well only in the case of well-prepared initial conditions\footnote{A well-prepared initial condition in the context of a numerical simulation or mathematical model refers to an initial state that is carefully selected to reflect the problem being studied and to ensure the stability and accuracy of the simulation. It means that the initial condition is chosen in a way that is physically meaningful and conducive to obtaining reliable results from the simulation. For instance, in the Xin Jin model, given $u(x,0)$ the initial condition are well-prepared if and only if $v(x,0) = g(u(x,0)).$}, i.e., given $u_0(x)$ and $g(u) = bu$ in \eqref{Xin_Jin_model} (linear case), $v_0(x) = bu_0(x) $. On the contrary, the CAT2-Tay scheme \eqref{Met_L-stalbe} performs effectively with both well-prepared and unprepared initial conditions.
Hence, in the numerical tests section, we will adopt well-prepared initial conditions for the Jin-Xin model test, whereas for the Broadell model, we considered both well-prepared and unprepared initial conditions (see the Sec.~\ref{sec:numerics} for numerical experiments and \cite{BoscarinoRusso} for further details).

\subsection{Fourier stability analysis}
In this section, we focus on investigate the linear stability analysis throughout the Fourier mode.
For this reason we will consider the Xin Jin relaxed model \eqref{Xin_Jin_model} in which $g(u) = au.$ 

Let us assume that $$u_j^{n} = \rho^ne^{{\bf i} jk\Delta x} \; {\rm and} \; v_j^{n} = \rho^ne^{{\bf i} jk\Delta x}.  $$ After some algebraic manipulation, the semi-implicit CAT2 scheme with trapezoidale rule \eqref{Met_trap_rule} becomes:
\begin{equation}
\label{Fourier_Trap}
    \begin{cases}
        \rho = 1 + \frac{\Delta t}{2\Delta x}\Big[-{\bf i}\sin(k\Delta x)\Bigl(1 + \frac{1}{1 + \frac{\Delta t}{\epsilon}}\Bigl(1 + a\frac{\Delta t}{\epsilon}\Bigr)\Bigr) +  2(\cos(k\Delta x) - 1)\frac{\Delta t}{\Delta x}\Bigl(1 + a\frac{\Delta t}{\epsilon}\Bigr)  \frac{1}{1 + \frac{\Delta t}{\epsilon}} \Big] \\ \vspace{-0.2cm}\\
        \rho = \frac{1 + \frac{\Delta t}{\Delta x}\Bigl({\bf i}\sin(k\Delta x) + \frac{\Delta t}{\Delta x}\Bigl(\cos(k\Delta x) - 1\Bigr)\Bigr) - \frac{\Delta t}{2\epsilon}(1-a)}{1 + \frac{\Delta t}{2\epsilon}(1-a)}.    
    \end{cases}
\end{equation}
The k-th Fourier mode is stable if the amplification coefficient $\rho = \rho(k)$ satisfies $|\rho|\le 1$. Therefore, given that $\rho = \alpha + \textbf{i}\beta$ and defining $\mu = \Delta t/\Delta x$, we impose the condition $|\rho|^2\le 1$ to obtain:
\begin{equation}
    \label{Four_Trap_k_mode}
    \begin{cases}
        \mu \le \frac{\sqrt{2}}{\sqrt{y_{\epsilon}}} \\ 
        \mu \le \frac{\sqrt{2}}{\sqrt{\sqrt{1+\frac{4}{z_{\epsilon}}} + 1}}
    \end{cases}    
\end{equation}
where $$y_{\epsilon} = \frac{1 + a\frac{\Delta t}{\epsilon} }{1 + \frac{\Delta t}{\epsilon}} \quad {\rm and} \quad z_{\epsilon} = \frac{\Delta t}{2\epsilon}(1-a).$$
\begin{figure}[!ht]
    \centering
         \hspace{-0.76cm}
         \includegraphics[width=0.5\textwidth]{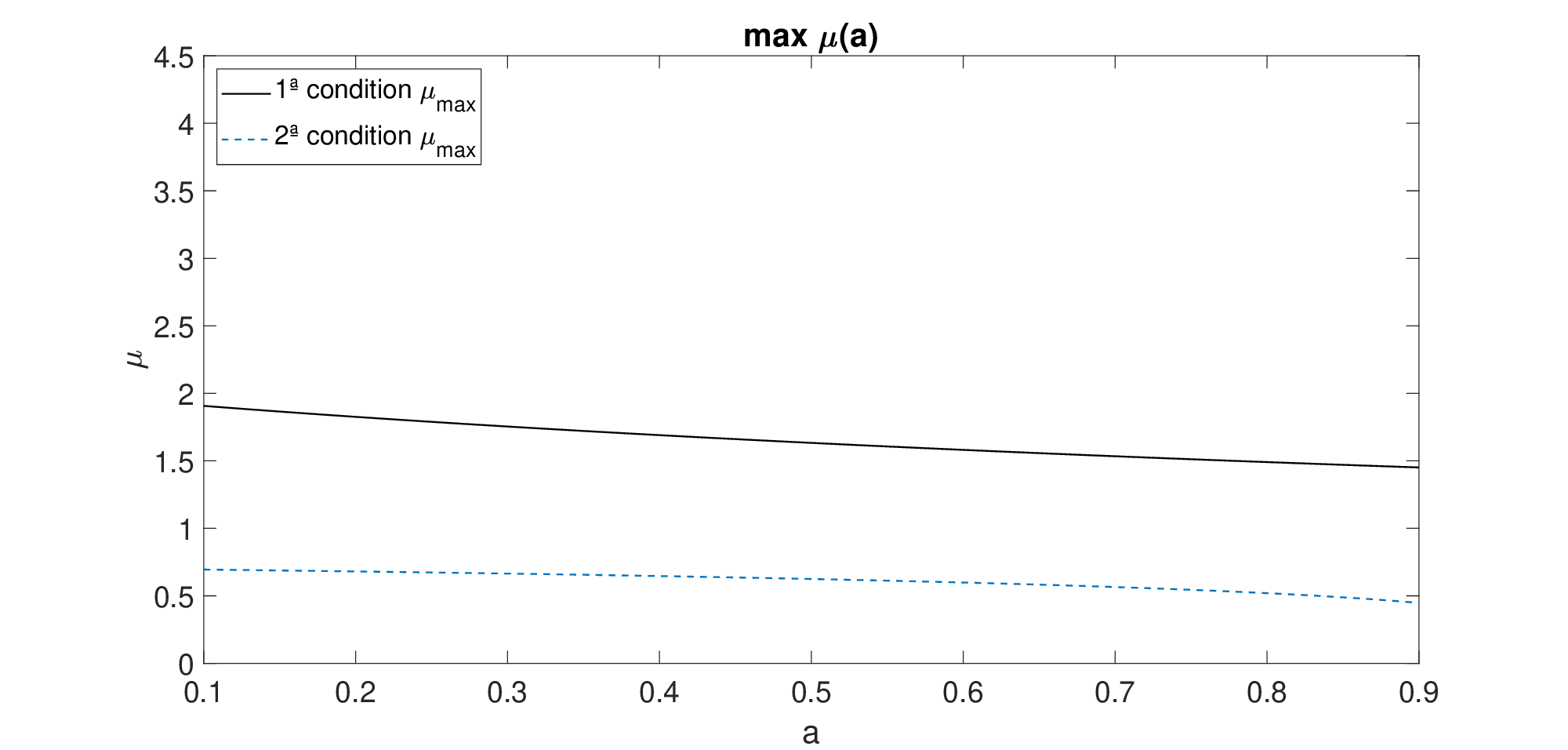}
         \includegraphics[width=0.5\textwidth]{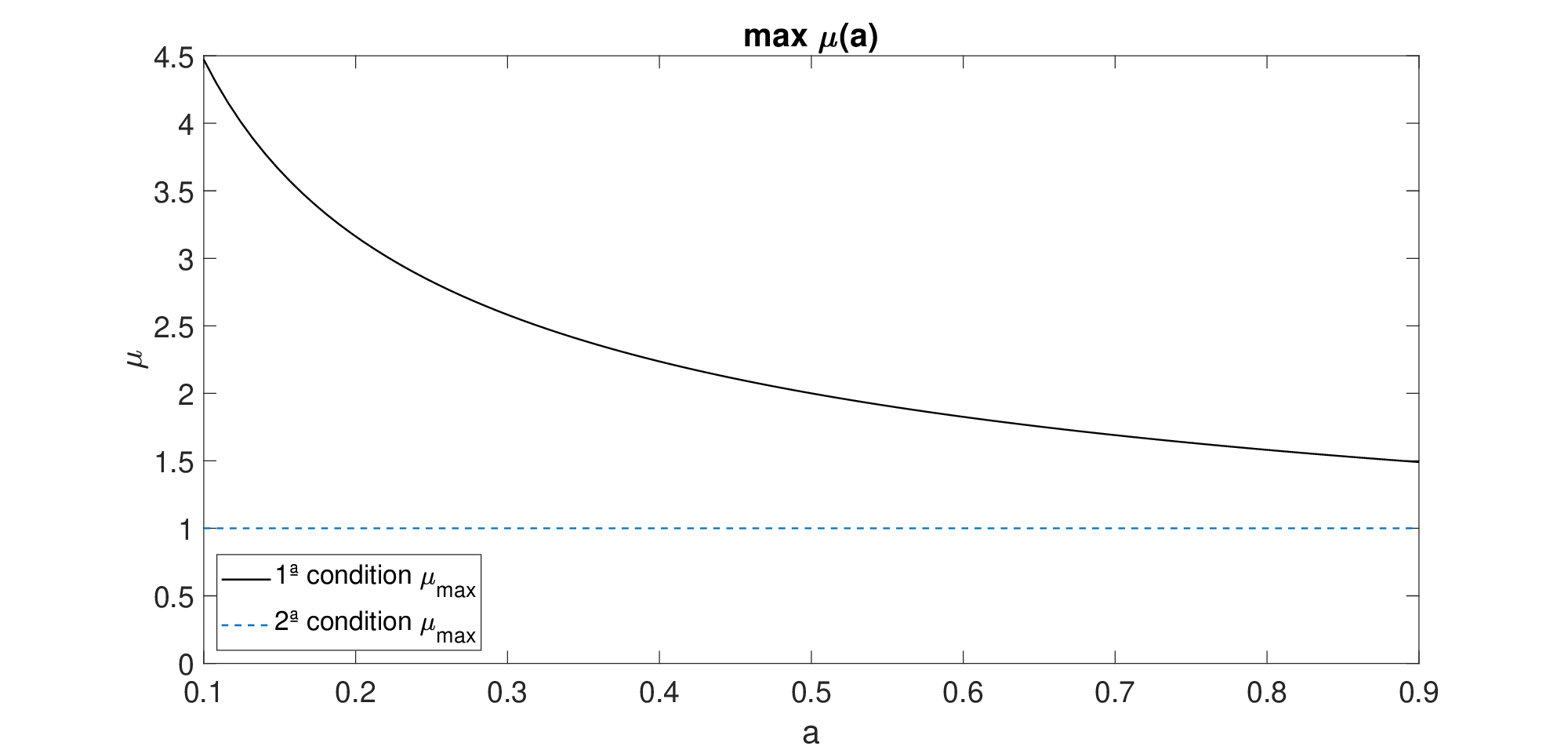}        
    \caption{Von Neumann stability analysis. Linear stability conditions for the semi-implicit CAT2 method with trapezoidale rule \eqref{Met_trap_rule} applied linear Xin Jin model \eqref{Xin_Jin_model} for $\epsilon= 1$ (left) and for $\epsilon\rightarrow0,\ (10^{-14})$ (right).}
   \label{fig:XinJin_Analysis_Trap}
\end{figure}
\begin{remark}
    When $\epsilon \rightarrow 0,$ assuming $a\neq0$\footnote{We can assume that $a\neq0$ since if $a=0$ (when $\epsilon\rightarrow0$), the only solution to Xin Jin's model is $v=0$ and $u=c.$}, the stability conditions are given by 
    \begin{equation}
    \label{Stab_mu_Trap}
        \begin{cases}
            \mu \le \frac{\sqrt{2}}{\sqrt{a}} \\ 
            \mu \le 1
        \end{cases}
    \end{equation}    
\end{remark}
Figure \ref{fig:XinJin_Analysis_Trap} displays the contour of the stability region for the semi-implicit CAT2 method utilizing the trapezoidale rule \eqref{Met_trap_rule} on the linear Xin Jin system \eqref{Xin_Jin_model}, while varying the parameters $a$ and $\epsilon$. Specifically, it showcases two conditions derived from equation \eqref{Four_Trap_k_mode} for large values of $\epsilon$, while exploring the parameter $a$ within the range of $[0.1, 0.9]$, as well as the condition derived from equation \eqref{Stab_mu_Trap} in the $\epsilon$ limit, while varying $a$ within the interval $[0.1, 0.9]$. It is apparent that in the extreme case, there is no dependence on $a$, satisfying the classical stability requirement $\mu\le1$. However, when $\epsilon$ is significantly large, an additional restriction arises when $a$ is greater than 0.5\footnote{In order to obtain explicit formulations, the stability conditions have been made more stringent. In fact, it has been empirically observed that less restrictive conditions can be applied to the CFL. In general, it has been observed empirically that $\mu\le0.95$.}.

In a similar manner to the approach used for CAT with trapezoidal rule, the formulation of CAT with IMEX Taylor can be derived as follows:
\begin{equation}
\label{Fourier_Tay}
    \begin{cases}
        \rho = 1 + \frac{\Delta t}{2\Delta x}\Big[-{\bf i}\sin(k\Delta x)\Bigl(1 + \frac{1}{1 + \frac{\Delta t}{\epsilon}}\Bigl(1 + a\frac{\Delta t}{\epsilon}\Bigr)\Bigr) +  2(\cos(k\Delta x) - 1)\frac{\Delta t}{\Delta x}\Bigl(1 + a\frac{\Delta t}{\epsilon}\Bigr)  \frac{1}{1 + \frac{\Delta t}{\epsilon}} \Big] \\ \vspace{-0.2cm}\\
        \rho = \frac{1 + \frac{\Delta t}{\Delta x}\Bigl({\bf i}\sin(k\Delta x) + \frac{\Delta t}{\Delta x}\Bigl(\cos(k\Delta x) - 1\Bigr)\Bigr) - \textbf{i}\frac{\Delta t}{2\epsilon\Delta x}(1-a)\sin(k\Delta x)}{1 + (\frac{\Delta t}{\epsilon}+\frac{\Delta t^2}{2\epsilon^2})(1-a)}.    
    \end{cases}
\end{equation}
\begin{figure}[!ht]
    \centering
         \hspace{-0.76cm}
         \includegraphics[width=0.5\textwidth]{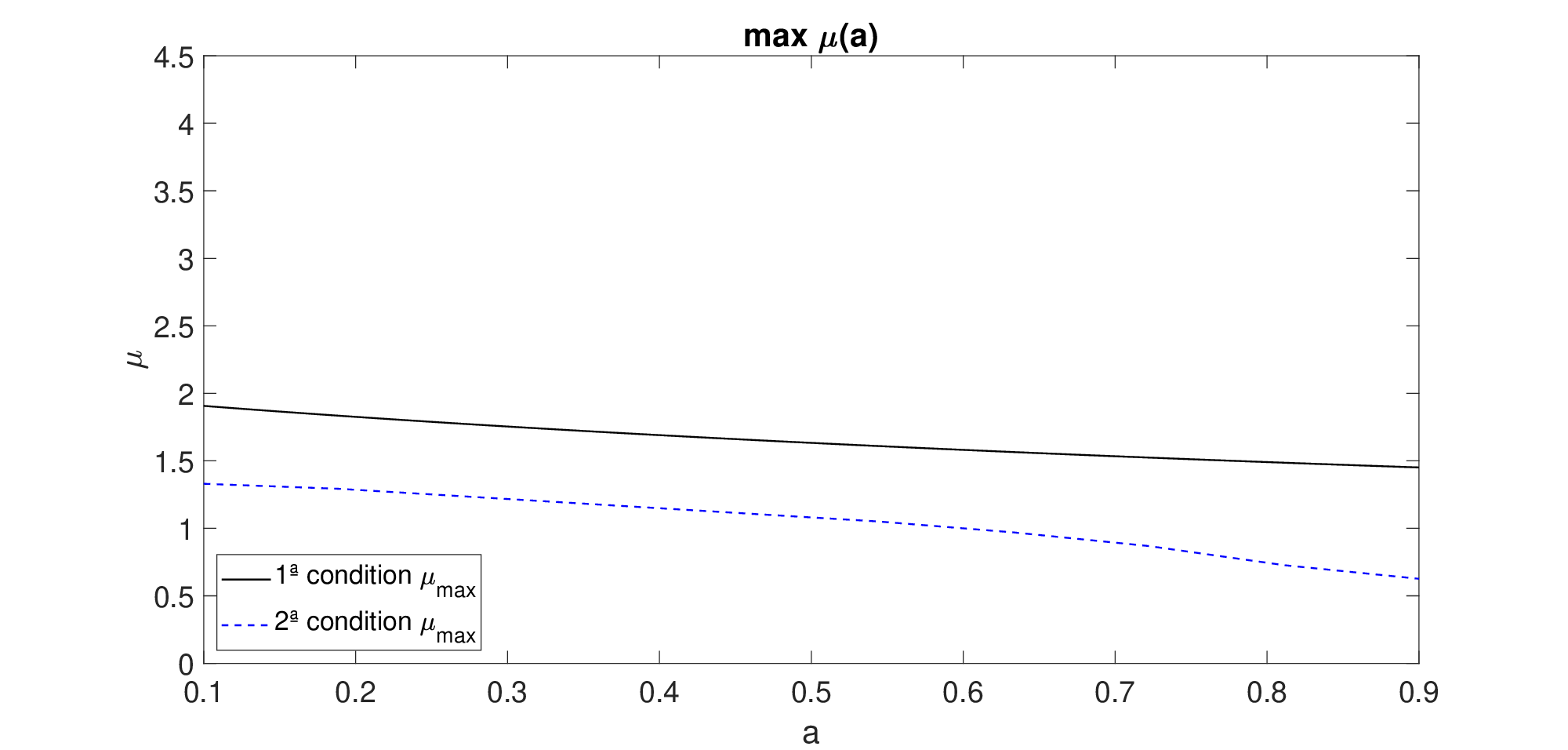}
         \includegraphics[width=0.5\textwidth]{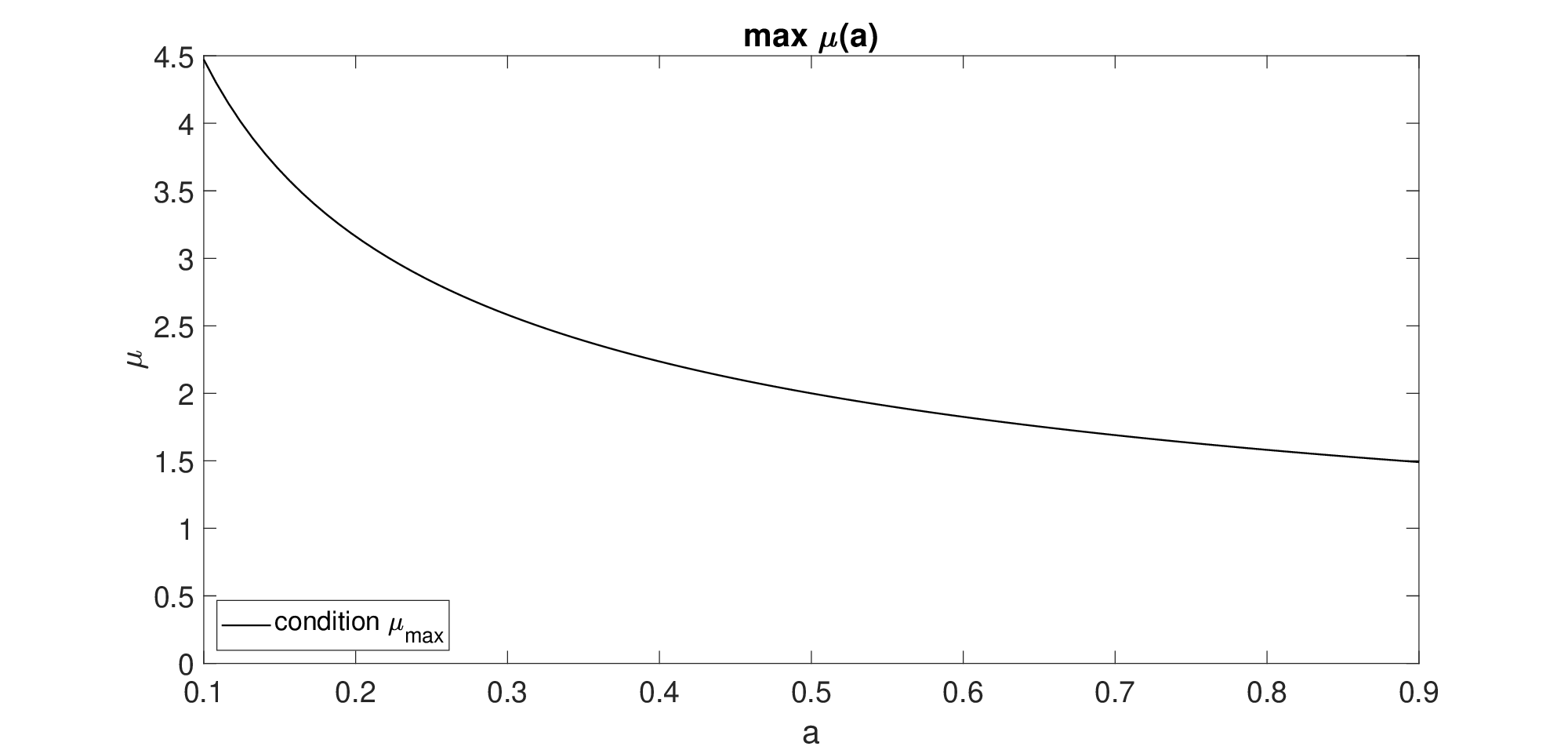}        
    \caption{Von Neumann stability analysis. Linear stability conditions for the semi-implicit CAT2 method with Taylor \eqref{Met_L-stalbe} applied linear Xin Jin model \eqref{Xin_Jin_model} for $\epsilon \approx 1$ (left) and for $\epsilon\rightarrow0 \, 10^{-14}$ (right).}
   \label{fig:XinJin_Analysis_Taylor}
\end{figure}
Analogously, we obtain:
\begin{equation}
    \label{Four_Tay_k_mode}
    \begin{cases}
        \mu \le \frac{\sqrt{2}}{\sqrt{y_{\epsilon}}} \\ 
        \mu \le  \frac{1}{\sqrt{2}}\sqrt{2 - w_{\epsilon,a} + \sqrt{(2 - w_{\epsilon,a})^2 + 2\gamma_{\epsilon,a} }}
    \end{cases}    
\end{equation}
where $$ w_{\epsilon,a} = (1+\cos(k\Delta x))\Bigl(1 - \frac{1-a}{2\epsilon}\Bigr)^2 $$ and $$\gamma_{\epsilon,a} = \zeta^2 + 2\zeta \quad {\rm with} \quad \zeta = \Bigl(\frac{\Delta t}{\epsilon} + \frac{\Delta t^2}{2\epsilon^2} \Bigr)(1-a).$$
\begin{remark}
    When $\epsilon \rightarrow 0,$ assuming $a\neq0$, the stability conditions are given by 
    \begin{equation}
    \label{Stab_mu_Tay}
        \begin{cases}
            \mu \le \frac{\sqrt{2}}{\sqrt{a}} \\ 
            \forall \mu\in\R.
        \end{cases}
    \end{equation}
\end{remark}
Figure \ref{fig:XinJin_Analysis_Taylor} shows the stability region contour for the semi-implicit CAT2 method with Taylor \eqref{Met_L-stalbe} applied to the linear Xin Jin system \eqref{Xin_Jin_model}, varying the parameters $a$ and $\epsilon$. Specifically, it displays two conditions derived from equation \eqref{Four_Tay_k_mode} for $\epsilon\gg0$ while varying the parameter $a$ in the range of $[0.1, 0.9]$, as well as the condition derived from equation \eqref{Stab_mu_Tay} in the limit of $\epsilon$ while varying $a$ in the range of $[0.1, 0.9]$. It is evident that in the limiting case, there is no dependence on $a$, falling into the classical stability condition $\mu\le1$; whereas, when $\epsilon\gg0$, there is an additional restriction when $a>0.5.$

\section{Numerical experiments}
\label{sec:numerics}
In this paper, our numerical tests specifically focus on 1D systems of balance laws with relaxed source terms. We apply the semi-implicit-type CATMOOD schemes, and show to be advantageous for systems with relaxed source terms, offering improved efficiency.

Our testing methodology involves evaluating the performance of CATMOOD schemes using a range of classical and challenging test cases. These test cases are carefully selected to assess the ability of the schemes to handle various types of flow phenomena and capture essential features accurately. Our specific focus lies on three main test cases introduced in Sec.~\ref{sec:Models}.

In the Xin Jin model \cite{jin1995runge}, we examine the behavior of the CATMOOD schemes for both smooth and well-prepared initial conditions. This allows us to obtain numerical solutions with high-order accuracy for both large and small values of $\epsilon$, thereby assessing the performance of the schemes in capturing the solution accurately under different parameter regimes.

For the Broadwell model, we investigate the capabilities of the schemes to capture the relaxed solution, particularly in non-smooth cases. By utilizing smooth and well-prepared initial conditions, we aim to obtain accurate numerical solutions that faithfully represent the relaxation phenomenon.

Finally, in the Euler system with heat transfer, we examine the ability of the schemes to handle non-smooth solutions while overcoming the introduction of spurious oscillations. This test case allows us to evaluate the robustness and accuracy of the schemes in scenarios involving complex flow phenomena and heat transfer effects.

Through these carefully selected test cases, we can  assess the performance and capabilities of the semi-implicit-type second order CATMOOD schemes in capturing a wide range of flow phenomena accurately.

The MOOD technique offers various degrees of freedom that, in general, do not excessively impact the solution quality despite their slight differences. Following the recommendations of the authors  for using this technique \cite{CDL0_FVCA,CDL1,CDL2,CDL3,CATMOOD_1D}, we decide to keep the values of $\epsilon_1$ and $\epsilon_2$  in \eqref{eq:delta} adaptable in relation to the specific test ($NAD$, see Sec.~\ref{ssec:MOODdetection}). {Additionally, numerical checks for positive or non-oscillatory solutions are conducted on the first equation of the system ($PAD$, see Sec.~\ref{ssec:MOODdetection}). For instance, this involves $u$ in the Xin Jin model, $\rho$ in the Broadwell model and Euler with heat transfer.}

The time step is chosen according to 
$$
\Delta t^n = \mathrm{CFL}\, \frac{\Delta x}{\lambda_{ {\rm max}}^n},   $$
where $\lambda_{{\rm max}}^n$ represents the maximum, over cells, of the spectral radius of the Jacobian matrix $\partial F/\partial u$.

\subsection{Xin Jin model} \label{X_J_num}
The numerical experiments setup in this section tests the convergence rates of the CAT2 schemes presented before for smooth initial condition and the numerical solutions for non-smooth initial condition related to the Xin Jin model \eqref{Xin_Jin_model}.

\paragraph{Smooth case} The experiment is carried out  with well-prepared smooth initial conditions. 
\begin{figure}[!ht]
    \centering	
    \includegraphics[width = \textwidth]{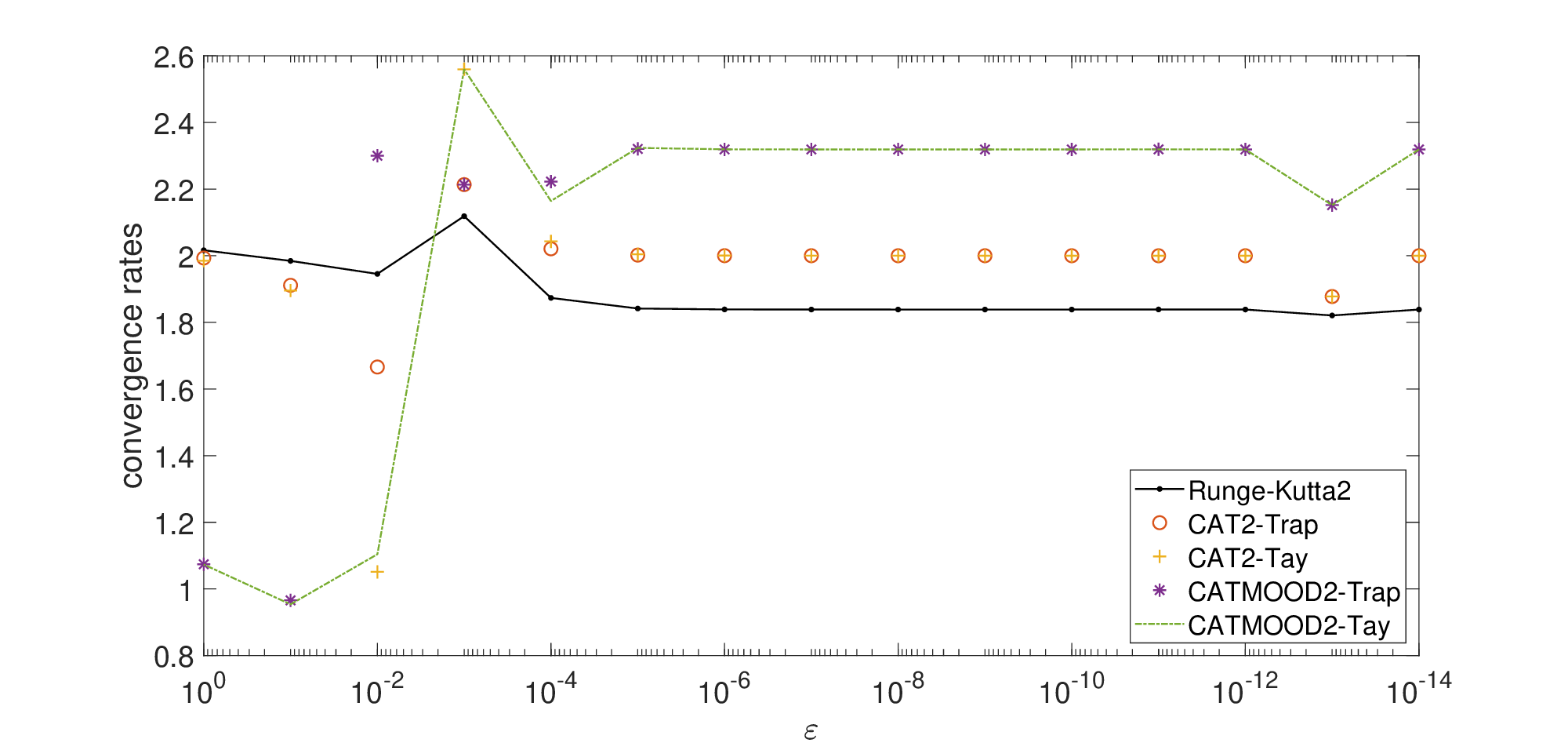}
    \caption{Xin Jin model: smooth case. Convergence rates for the numerical solution $u$ obtained with Runge-Kutta2, CAT2-Trap, CAT2-Tay, CATMOOD2-Trap and CATMOOD2-Tay at time $t = 1$ on the interval $[0,1]$ with CFL$=0.9$. $\epsilon$ takes values on $\{10^{0},\ldots,10^{-14}\}.$ The tolerances $\epsilon_1$ and $\epsilon_2$ are, respectively, set $10^{-3}$ and $10^{-2}.$  }
    \label{Test_X_J_smooth_1}
\end{figure}
\begin{figure}[!ht]
    \centering	
    \includegraphics[width = \textwidth]{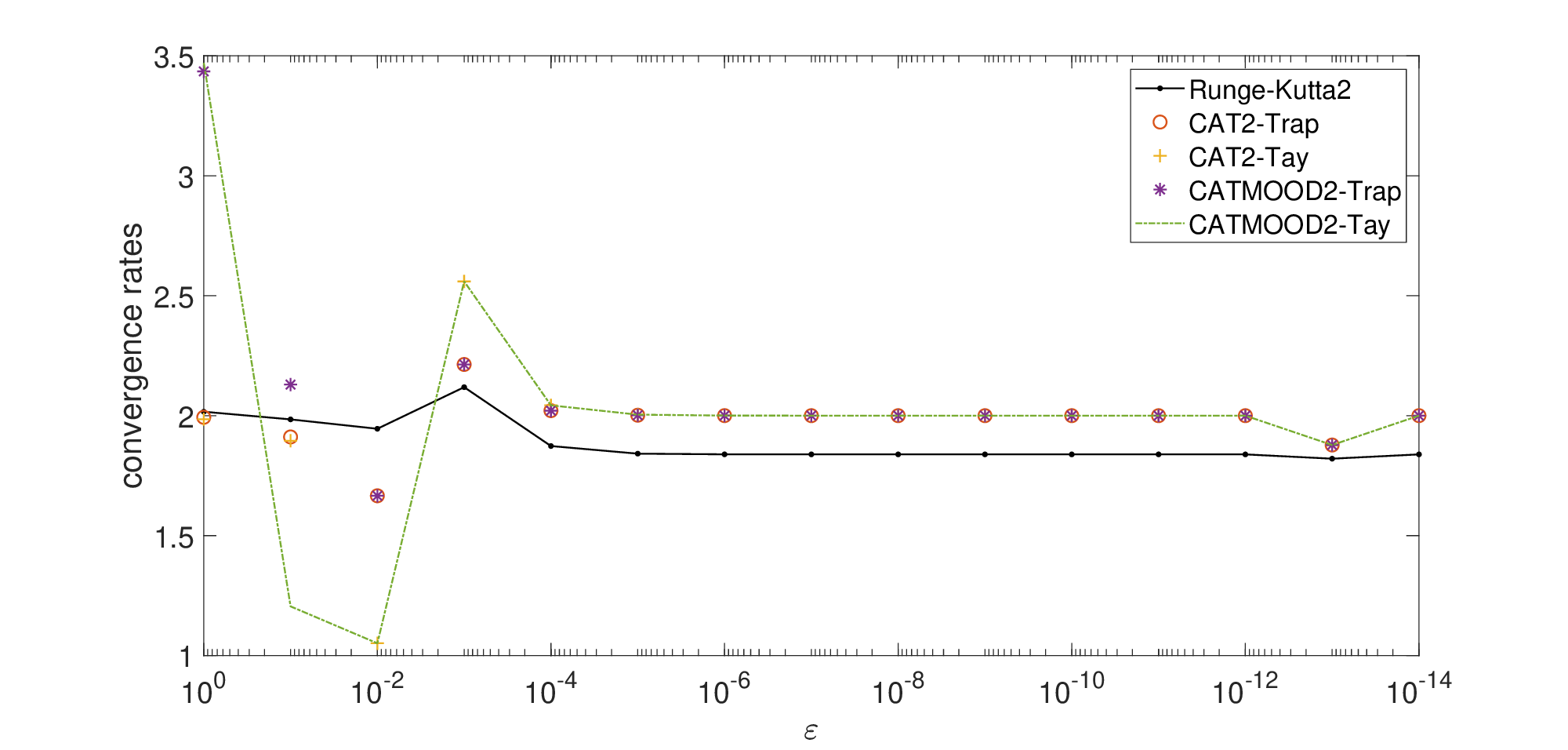}
    \caption{Xin Jin model: smooth case. Convergence rates for the numerical solution $u$ obtained with Runge-Kutta2, CAT2-Trap, CAT2-Tay, CATMOOD2-Trap and CATMOOD2-Tay at time $t = 1$ on the interval $[0,1]$ with CFL$=0.9$. $\epsilon$ takes values on $\{10^{0},\ldots,10^{-14}\}.$ The tolerances $\epsilon_1$ and $\epsilon_2$ are, respectively, set $10^{-2}$ and $10^{-1}.$  }
    \label{Test_X_J_smooth_2}
\end{figure}
We set: $[\alpha,\beta] = [0, 1]$ the computational domain, $t_{\text{fin}} = 1$ the final time, CFL$=0.9$, $N=200$ grid points, periodic boundary conditions, $\epsilon$ is taken within the range $\{10^0, 10^{-14}\}$ with decrements of $1/10$ and the tolerances $\epsilon_1$ and $\epsilon_2$ are chosen  $10^{-3}$ and $10^{-2}$ respectively. 
\begin{table}[!ht]
\resizebox{\columnwidth}{!}{
\begin{tabular}{|c||c|c|c|c|c|c|c|}
\hline \multicolumn{8}{|c|}{\textbf{CPU time in seconds}} \\
\hline   
\hline   
$N,\epsilon$ & \textbf{CM2-Trap$^*$} & \textbf{CM2-Tay$^*$} &\textbf{CM2-Trap} & \textbf{CM2-Tay} & \textbf{C2-Trap} & \textbf{C2-Tay} & \textbf{RK2}\\ \hline
128, $10^{0} $ & 0.0172  & 0.0175 & 0.0694 & 0.0667 & 0.0096 & 0.01 & 0.0192 \\ \hline
4096, $10^{0} $ & 0.9411  & 0.9609 & 0.9483 & 0.9689 & 0.4617 & 0.4753 & 1.5490 \\ \hline
128, $10^{-14} $  & 0.0141  & 0.0161 & 0.0172 & 0.0168 & 0.0081 & 0.0087 & 0.223 \\ \hline
4096, $10^{-14}$ & 0.9646  & 0.9683 & 0.9710 & 0.9683 & 0.4655 & 0.4670 & 1.4858 \\ \hline
\end{tabular}}
\caption{Xin Jin model: smooth case. CPU time expressed in seconds for the variable $u$ obtained through the CATMOOD2-Trap, CATMOOD2-Tay, CAT2-Trap, CAT2-Tay, and {second-order semi-implicit RK schemes} on different mesh and different $\epsilon$. For the CATMOOD2-Trap and CATMOOD2-Tay schemes, the tolerances, $\epsilon_1$ and $\epsilon_2$, are set to $10^{-3}$ and $10^{-2}$ respectively. Meanwhile, for CATMOOD2-Trap$^*$ and CATMOOD2-Tay$^*$ schemes, the tolerances are adjusted to $10^{-2}$ and $10^{-1}$ respectively.}
\label{X_J_table_1}
\end{table}


Here, the well-prepared initial condition $(u_0,v_0) = (u(x,0),v(x,0))$ is given by 
\begin{equation}
    (u_0,v_0) = (1 + \sin(2\pi x), a + a\sin(2\pi x)) 
\end{equation}
where $v_0 = g(u_0) = au_0$ and $a = 0.7.$

Figure~\ref{Test_X_J_smooth_1} displays a the convergence rate for different schemes as we vary the $\epsilon$ parameter. The errors were obtained using the $L^1-$norm, i.e., ${\rm err}= || u_{\rm num} - u_{\rm exact}||_1.$

{Specifically, we consider schemes: Runge-Kutta2, as well as CAT2-Trap \eqref{Met_trap_rule} and CAT2-Tay \eqref{Met_L-stalbe}, respectively. Additionally, we explore CATMOOD2-Trap and CATMOOD2-Tay, which integrate the CAT2 methods \eqref{Met_trap_rule} and \eqref{Met_L-stalbe}  with the a-posteriori MOOD limiter.} 

{In Figure~\ref{Test_X_J_smooth_1}, an interesting trend emerges: adopting a coarse mesh, as $\epsilon$ approaches zero, all the schemes exhibit a second-order convergence behavior. However in the same mesh, when $\epsilon\approx 1$, a subtle shift occurs.} In this case, the Runge-Kutta2 method, along with CAT2-Trap and CAT2-Tay, still maintains its second-order convergence pattern. Anyway, when combining the CAT2-Trap or CAT2-Tay schemes with the MOOD technique, a noticeable drop in the order of convergence is observed.

{This behavior is closely tied to the fundamental characteristics of the system in question. As $\epsilon$ approaches zero, the conservativeness of the system becomes increasingly significant compared to the impact of the source term. This observation aligns with the fundamental principles of CAT schemes. On the other hand, when $\epsilon$ equals 1, the slight oscillations introduced trigger the MOOD technique, resulting in the observed reduction in convergence order.}

To address this challenge, two potential strategies emerge. Firstly, grid refinement emerges as an effective way to mitigate the order loss. Secondly, adjusting the tolerances $\epsilon_1$ and $\epsilon_2$ of the detector (NAD) shows promise, as depicted in Figure~\ref{Test_X_J_smooth_2}.

{Semi-implicit-type CAT2 schemes, whether they incorporate limiters or not, are not uniformly accurate on $\epsilon$ values ($10^{-1} \le \epsilon \le 10^{-2}$). This particular phenomenon aligns with earlier experimental findings, \cite{Boscarino-Filbet, PareschiRusso, BoscarinoRusso}.}

Finally, a crucial aspect deserving attention is highlighted in Table~\ref{X_J_table_1}. The table illustrates the computational cost in seconds, CPU time, associated with the various methods presented earlier, considering different MOOD tolerance settings, varied $\epsilon$ values, and distinct mesh sizes. Notable differences can be observed, influenced not only by the grid features but significantly affected by the value of $\varepsilon$ and hence, the structure of the underlying system.

Particularly, in cases where $\epsilon\gg0$ and a coarser grid is used, the selected tolerance values significantly impact the increment in CPU time. This phenomenon, however, does not hold true for instances where a finer grid is used. Conversely, as $\epsilon\rightarrow0$, regardless of the grid precision, the effect of tolerance values becomes negligible. This is because the conservative part has more influence than the source term, thereby leading to a diminished presence of computational artifacts arising from source term reconstruction.

\begin{remark}
   {Lastly, as mentioned in Sec.~\ref{sec:ode_like_stability}, for this numerical example we observe that if we use smooth well-prepared initial conditions there is no significant difference between the schemes \eqref{Met_trap_rule} and \eqref{Met_L-stalbe}.}
\end{remark}

\begin{remark}
 {   When using a fine grid, the semi-implicit RK method is more computationally intensive than other methods. However, in comparison to CAT2 schemes using the MOOD technique, the computational effort is comparable, although consistently higher than traditional CAT schemes. This difference is mainly because CAT schemes are local and one-step approaches, avoiding the additional intermediate steps to achieve higher-resolution results.}
\end{remark}

\paragraph{Non-smooth case}  The experiment is carried out with well-prepared non-smooth initial conditions:
\begin{align} 
\label{X_J:IC_disc_1}
	&\begin{cases}u(x,0) = \begin{cases}
	    2 \mbox { if } 0.25<x<0.5 \\
        1 \mbox { otherwise}     
	\end{cases} \\
    v(x,0) = au(x,0),
    \end{cases}
\end{align}
with $a=0.7.$
\begin{figure}[!ht]
     \centering
     \begin{subfigure}[b]{0.49\textwidth}
         \centering
         \includegraphics[width=\textwidth]{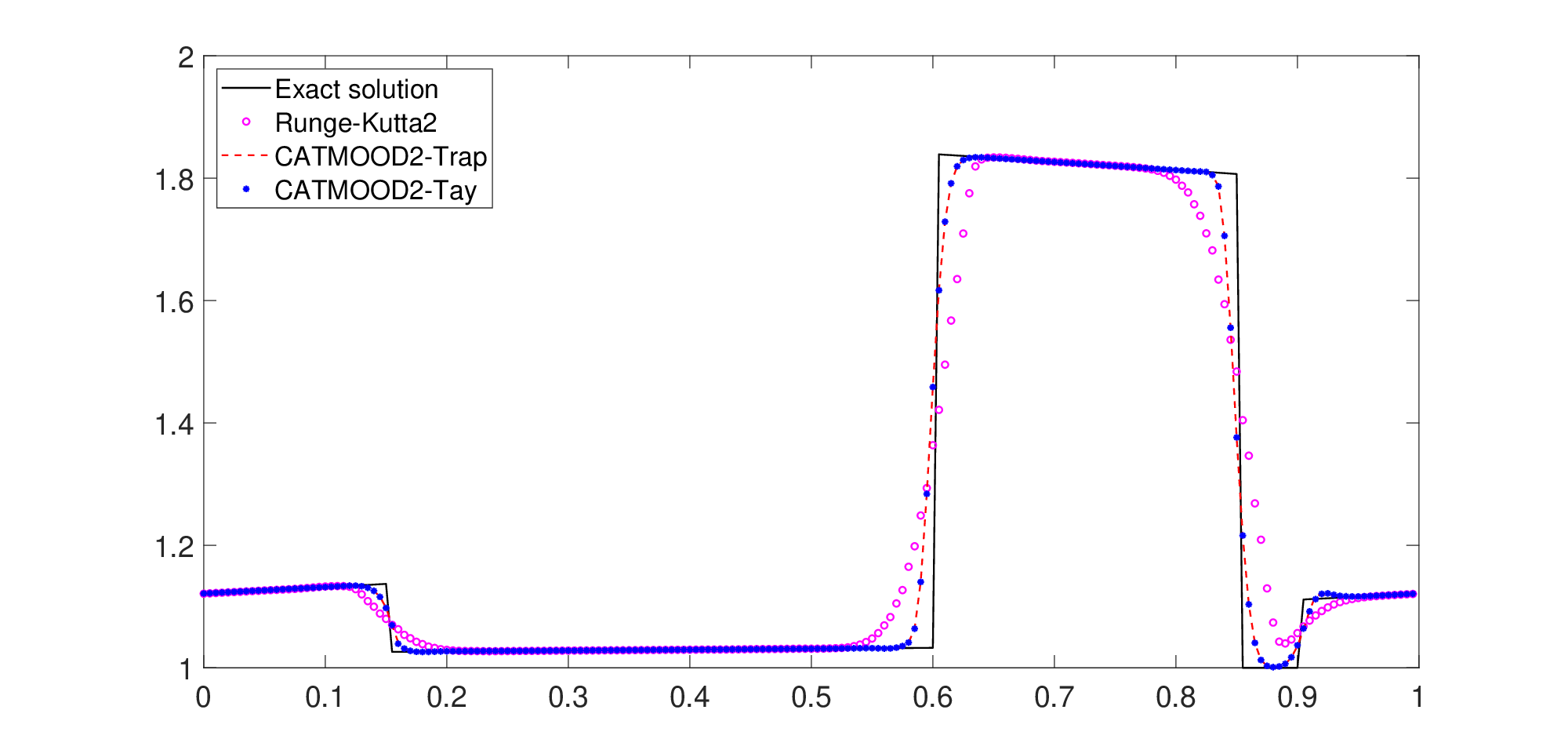}
         \caption{$u.$}
         \label{X_J:disc_u_1}
     \end{subfigure}
     \hfill
     \begin{subfigure}[b]{0.49\textwidth}
         \centering
         \includegraphics[width=\textwidth]{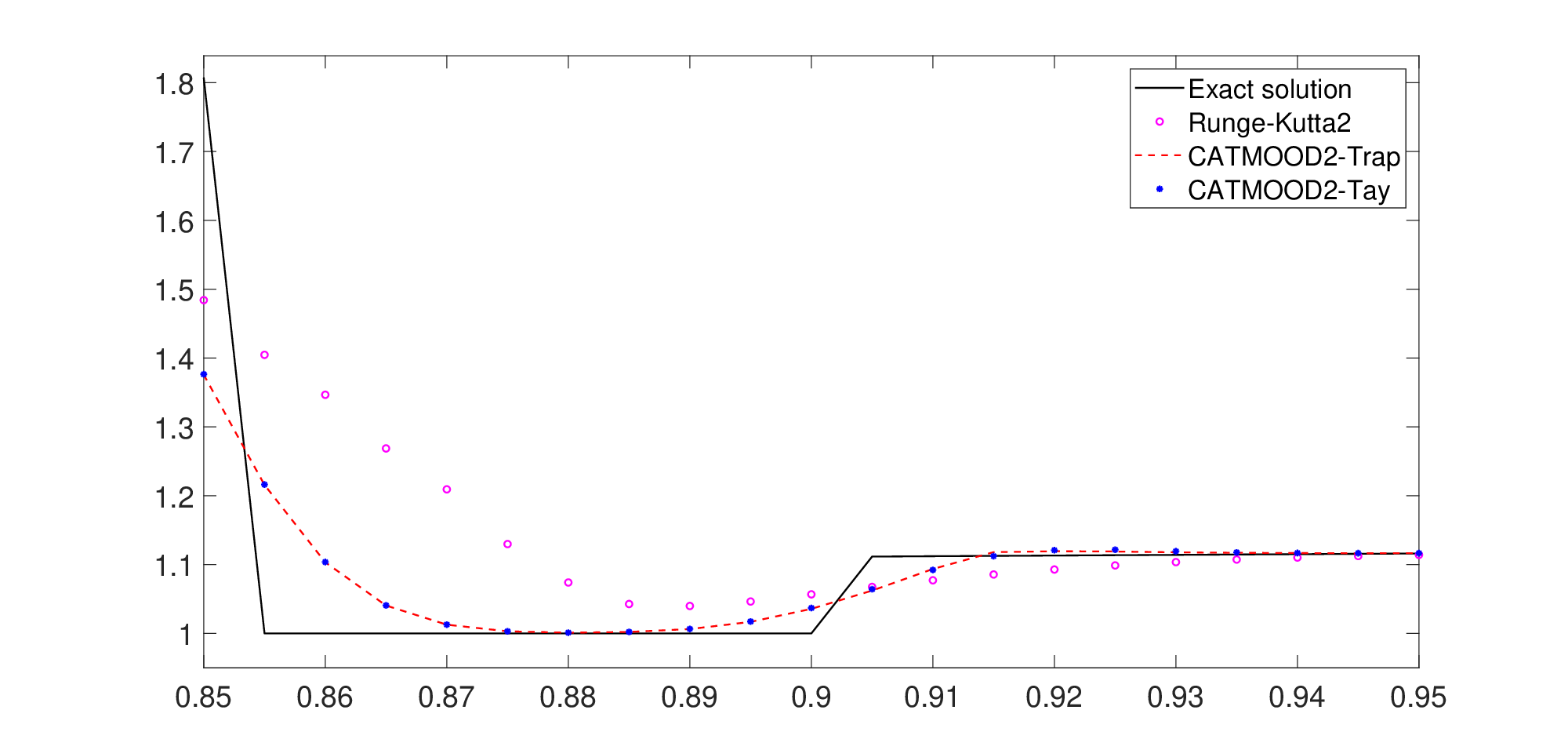}
         \caption{Zoom of $u$ on interval $[0.85,0.95].$}
         \label{X_J:disc_u_1_zoom}
     \end{subfigure} \\ \vspace{0.2cm}
     \begin{subfigure}[b]{0.49\textwidth}
         \centering
         \includegraphics[width=\textwidth]{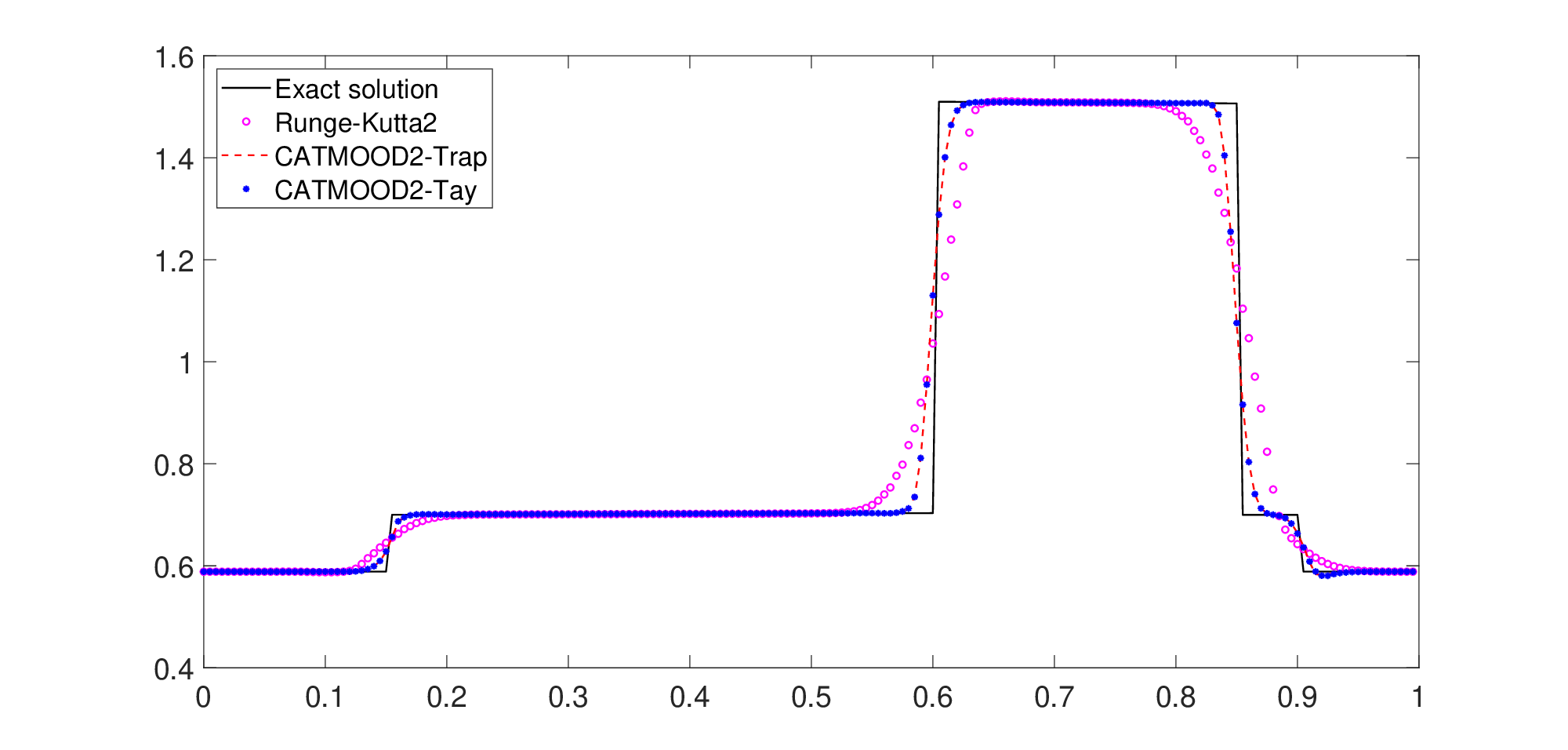}
         \caption{$v.$}
         \label{X_J:disc_v_1}
     \end{subfigure}
     \hfill
     \begin{subfigure}[b]{0.49\textwidth}
         \centering
         \includegraphics[width=\textwidth]{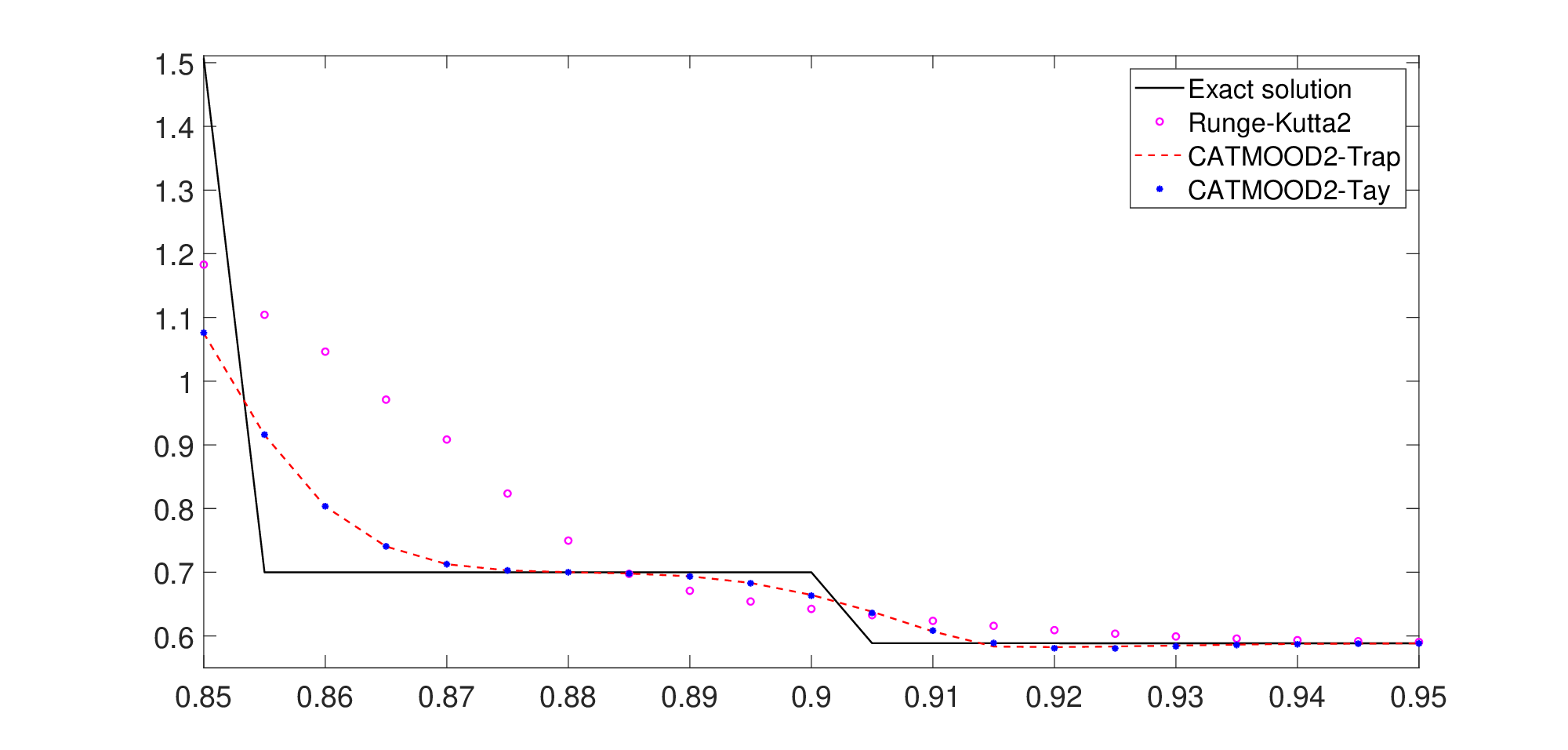}
         \caption{Zoom of $v$ on interval $[0.85,0.95].$}
         \label{X_J:disc_v_1_zoom}
     \end{subfigure}
     \caption{Xin Jin model: non-smooth case. Numerical solution $u$ (top-left) and $v$ (bottom-left) and zoom of $u$ (top-right) and $v$ (bottom-right) obtained with Runge-Kutta2, CATMOOD2-Trap and CATMOOD2-Tay at time $t = 0.35$ on the interval $[0,1]$ with CFL$=0.9$ and $\epsilon$ takes value $1.$ The tolerances $\epsilon_1$ and $\epsilon_2$ in the MOOD technique are, respectively, set $10^{-4}$ and $10^{-3}.$ The reference solution is the exact one.}
     \label{X_J:disc_1}
\end{figure}

 The computational domain $[\alpha,\beta]$ spans from $0$ to $1$, the final time $t_{\text{fin}}$ varies within the set $\{0.35,3\}$, a CFL number of $0.9$ is employed, $N=200$ points, the boundary conditions are periodic, $\epsilon$ is selected from the set $\{1, 10^{-8}\}$; and the tolerance values $\epsilon_1$ and $\epsilon_2$ are chosen $10^{-4}$ and $10^{-3}$ respectively.
\begin{figure}[!ht]
     \centering
     \begin{subfigure}[b]{0.49\textwidth}
         \centering
         \includegraphics[width=\textwidth]{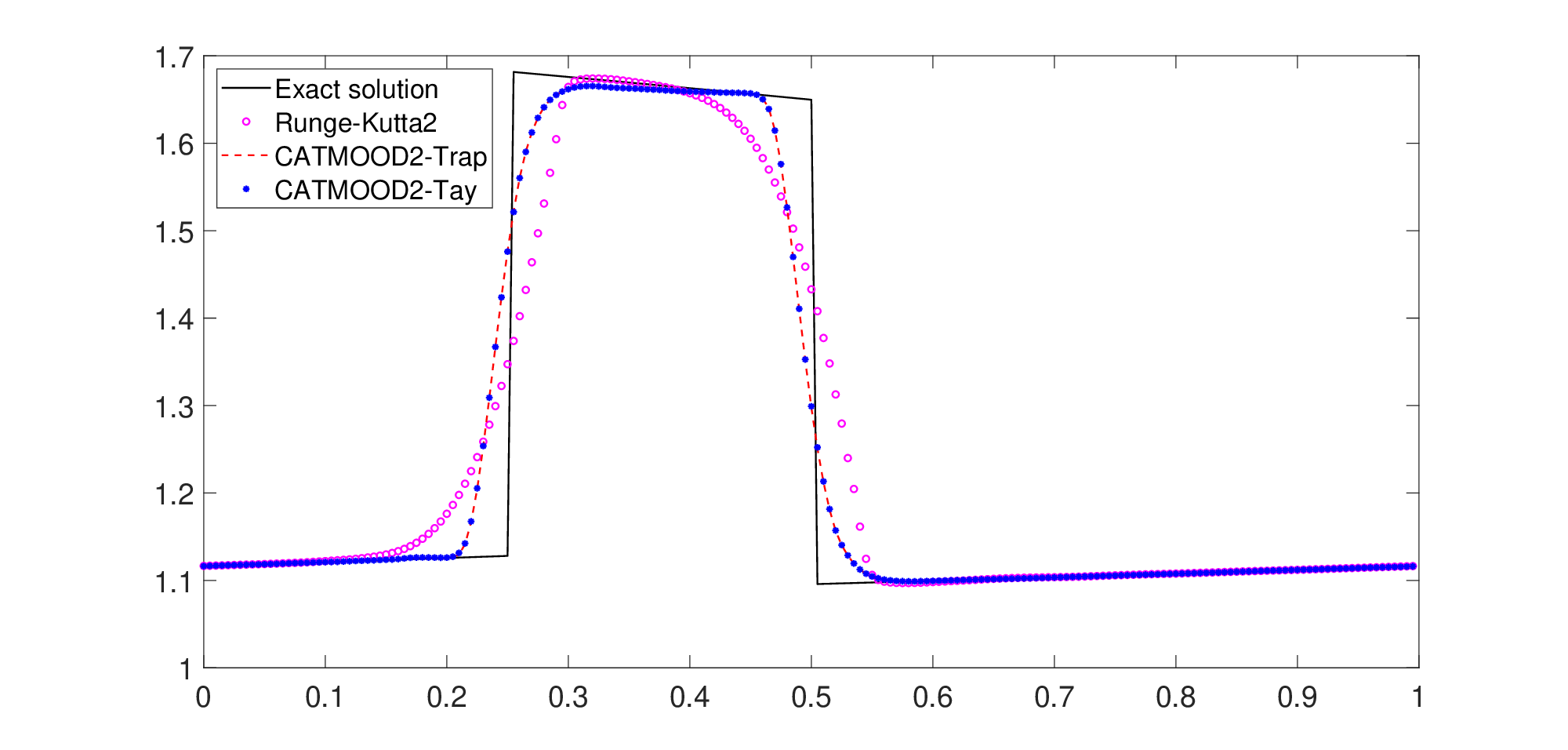}
         \caption{$u$ and $\epsilon = 1$}
         \label{X_J:disc_u_2}
     \end{subfigure}
     \hfill
     \begin{subfigure}[b]{0.49\textwidth}
         \centering
         \includegraphics[width=\textwidth]{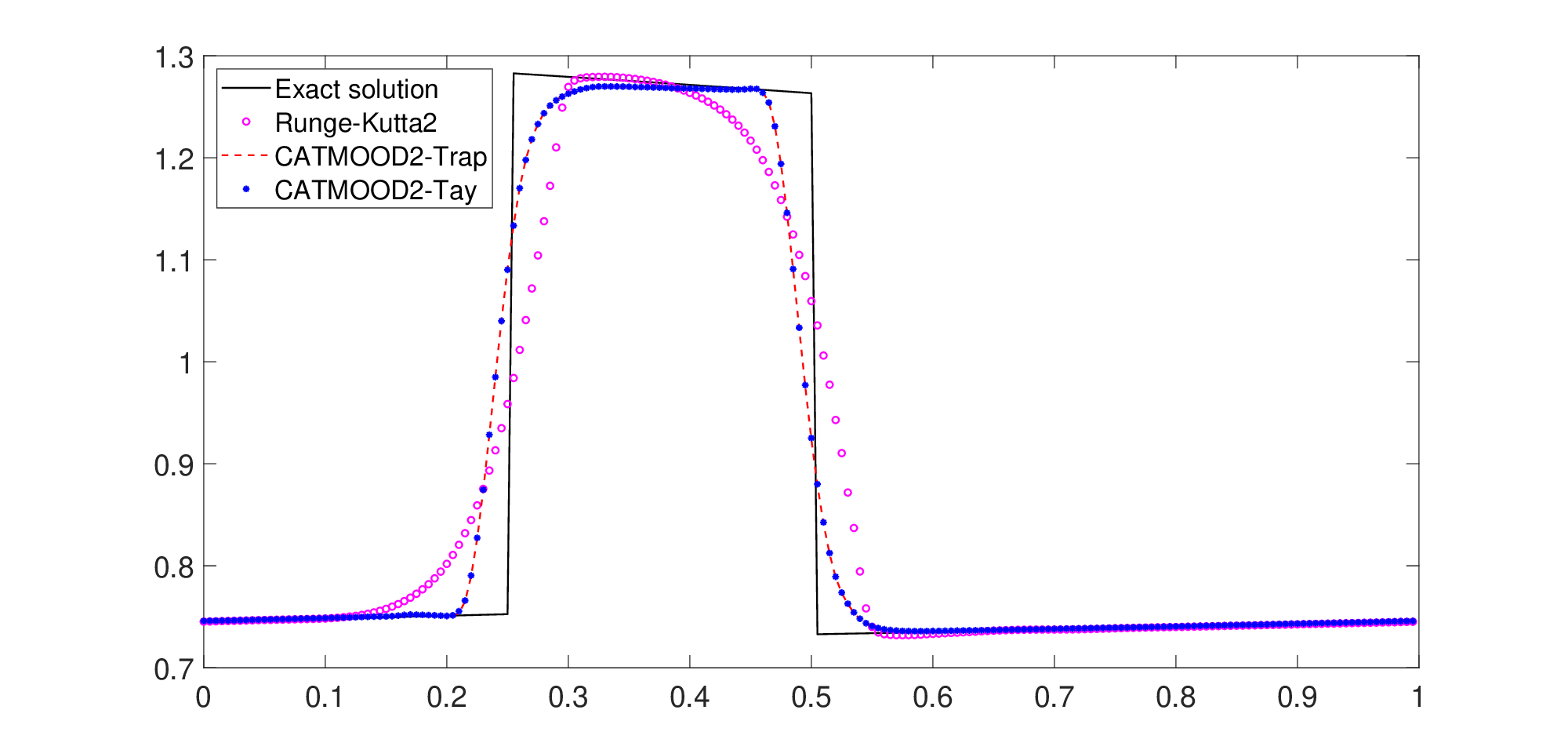}
         \caption{$v$ and $\epsilon = 1$}
         \label{X_J:disc_v_2}
     \end{subfigure} \\ \vspace{0.2cm}
     \begin{subfigure}[b]{0.49\textwidth}
         \centering
         \includegraphics[width=\textwidth]{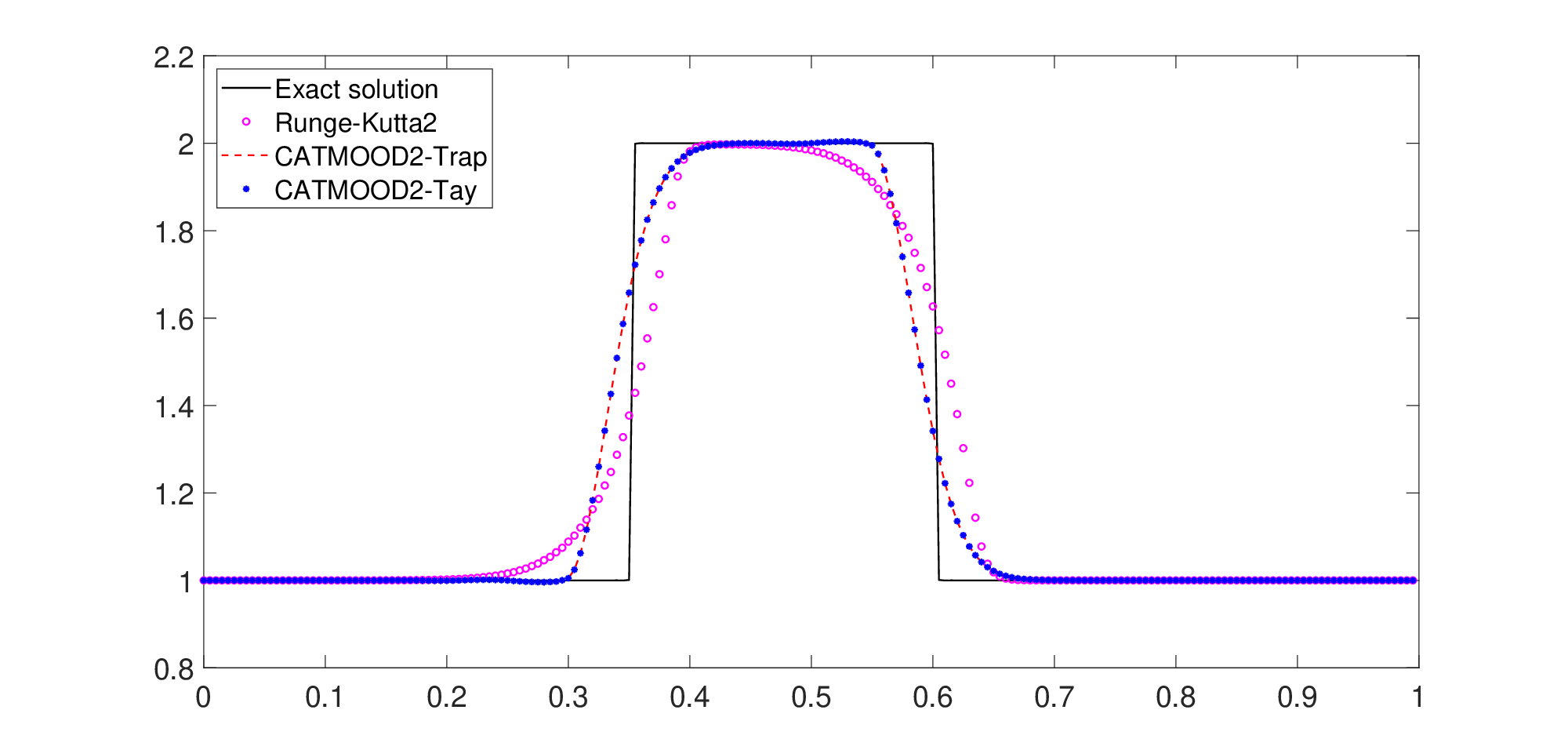}
         \caption{$u$ and $\epsilon = 10^{-8}$}
         \label{X_J:disc_u3}
     \end{subfigure}
     \hfill
     \begin{subfigure}[b]{0.49\textwidth}
         \centering
         \includegraphics[width=\textwidth]{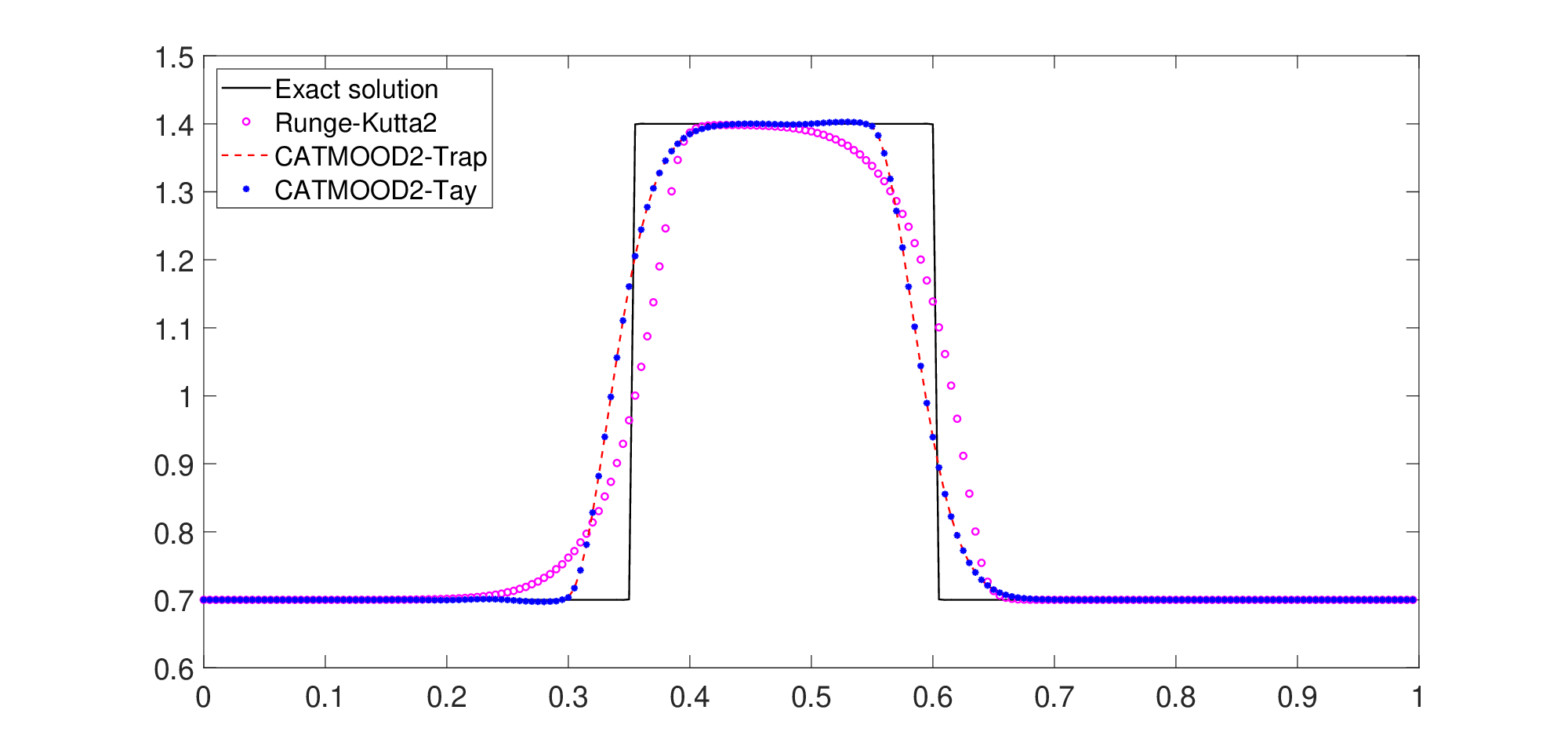}
         \caption{$v$ and $\epsilon = 10^{-8}$}
         \label{X_J:disc_v_3}
     \end{subfigure}
     \caption{Xin Jin model: non-smooth case. Numerical solution $u$ (left) and $v$ (right) obtained with Runge-Kutta2, CATMOOD2-Trap and CATMOOD2-Tay at time $t = 3$ on the interval $[0,1]$ with CFL$=0.9$ and $\epsilon$ takes value $1$ (top) and $10^{-8}$ (bottom). The tolerances $\epsilon_1$ and $\epsilon_2$ in the MOOD technique are, respectively, set $10^{-4}$ and $10^{-3}.$ The reference solution is the exact one.}
     \label{X_J:disc_2}
\end{figure}

In Figure~\ref{X_J:disc_1}, the numerical solutions for the variables $u$ and $v$ are displayed. These solutions are obtained using the second order semi-implicit Runge-Kutta2, CATMOOD2-Trap, and CATMOOD2-Tay schemes and they are compared with the exact one. The simulations conclude at the final time $t_{\text{fin}} = 0.35$. A notable observation is that the solutions generated by the CATMOOD methods exhibit a qualitatively higher level of accuracy when compared to the solutions produced by the Runge-Kutta2 method in presence of shocks. This noteworthy trend is further highlighted in Figures~\ref{X_J:disc_u_1_zoom}-\ref{X_J:disc_v_1_zoom}.

Figure~\ref{X_J:disc_2} illustrates the numerical solution for the variables $u$ and $v$ obtained using the Runge-Kutta2, CATMOOD2-Trap, and CATMOOD2-Tay schemes at the final time $t_{\text{fin}} = 0.35$, considering $\epsilon = 1$ and $\epsilon = 10^{-8}$. We observe that the semi-implicit CATMOOD schemes yield numerically more accurate solutions in comparison to the Runge-Kutta2 solutions.

\begin{remark}
Whether the initial condition is smooth or not, the semi-implicit CAT2 schemes coupled with the MOOD technique outperforms over the second order semi-implicit RK scheme. In the first case, there's a substantial computational advantage, whereas in the second case, a significant improvement in accuracy is obtained.
\end{remark}

\subsection{Broadwell model}
\paragraph{Smooth case} Let us consider the Broadwell model \eqref{Broadwell_model}. 
\begin{figure}[!ht]
    \centering	
    \includegraphics[width = \textwidth]{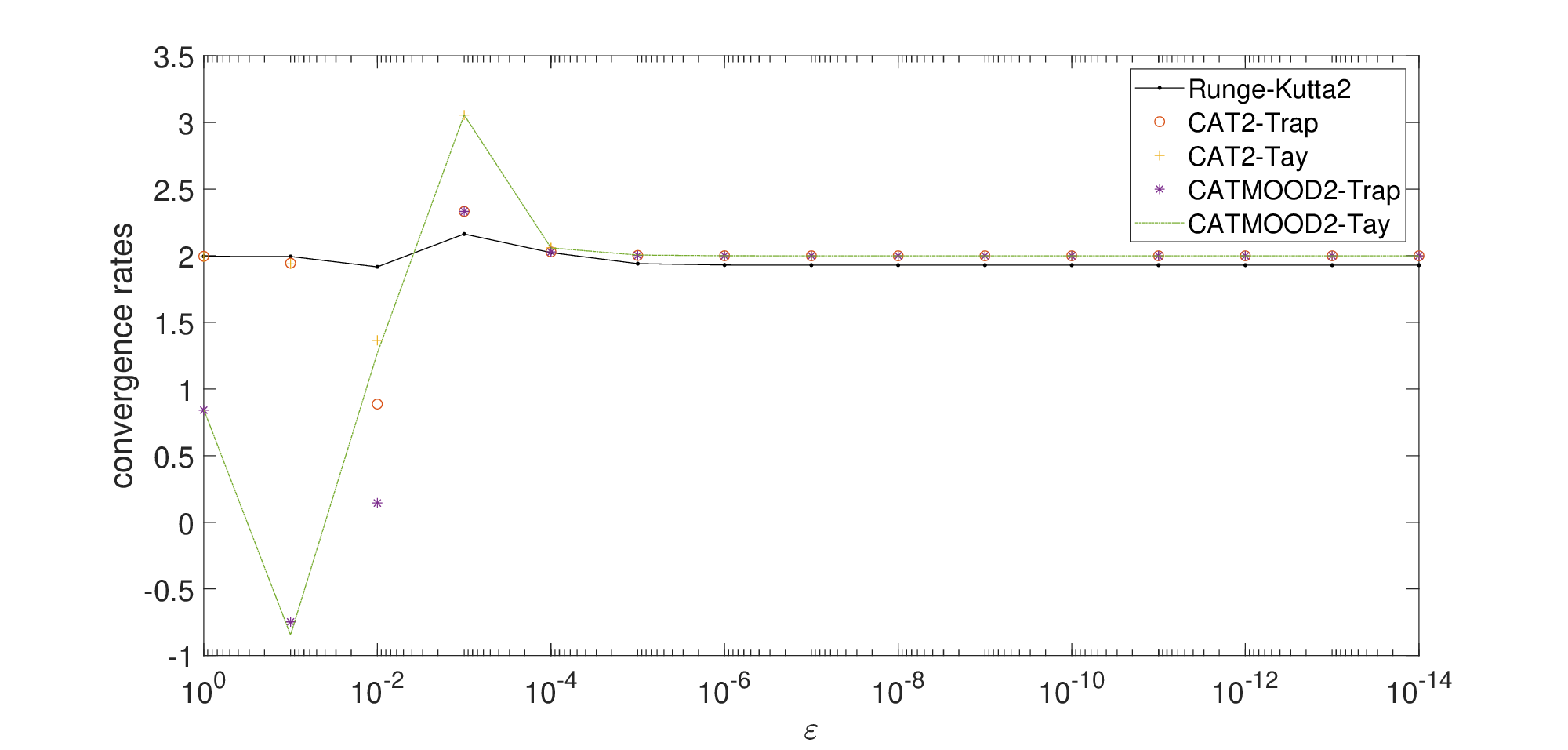}
    \caption{Broadwell model: smooth case. Convergence rates for the numerical solution $\rho$ obtained with Runge-Kutta2, CAT2-Trap, CAT2-Tay, CATMOOD2-Trap and CATMOOD2-Tay at time $t = 1$ on the interval $[0,1]$ with CFL$=0.9$. $\epsilon$ takes values on $\{10^{0},\ldots,10^{-14}\}.$ The tolerances $\epsilon_1$ and $\epsilon_2$ are, respectively, set $10^{-3}$ and $10^{-2}.$  }
    \label{Test_Broad_smooth_1}
\end{figure}
\begin{table}[!ht]
    \centering
\begin{tabular}{|c||c|c|c|c|c|}
    \hline \multicolumn{6}{|c|}{\textbf{CPU time in seconds}} \\
\hline   
\hline   
$N,\epsilon$ & \textbf{CM2-Trap} & \textbf{CM2-Tay} & \textbf{C2-Trap} & \textbf{C2-Tay} & \textbf{RK2}\\ \hline
256, $10^{(0)} $ & 0.157  & 0.163 & 0.067 & 0.074 & 0.152 \\ \hline
4096, $10^{(0)} $ & 21.03  & 24.82 & 11.74 & 14.41 & 28.32 \\ \hline
256, $10^{(-14)} $  & 0.119  & 0.131 & 0.073 & 0.079 & 0.142 \\ \hline
4096, $10^{(-14)} $ & 22.15  & 22.93 & 12.69 & 12.94 & 29.51 \\ \hline
\end{tabular} 
\caption{Broadwell model: Smooth case. The computational time, measured in seconds, is presented for the variable $\rho$ utilizing various schemes: CATMOOD2-Trap, CATMOOD2-Tay, CAT2-Trap, CAT2-Tay, and second-order Runge-Kutta IMEX schemes. These computations are performed on diverse meshes and with varying values of $\epsilon$. In the case of the CATMOOD2-Trap and CATMOOD2-Tay schemes, specific tolerances $\epsilon_1$ and $\epsilon_2$ are set at $10^{-3}$ and $10^{-2}$ respectively.}
\label{Broad_table_1}
\end{table}
We consider well-prepared smooth initial conditions:
\begin{equation}
     (\rho,v,z) = \left(1+0.3\sin(2\pi x),0.5+0.1\sin(2\pi x),0.5\rho(1+v^2)\right),
\end{equation}
with the domain set as $[0,1]$. The final time is $t_{\text{fin}} = 1$, a CFL number of 0.9, a grid nodes $N = 200$, and periodic boundary conditions. The parameter $\epsilon$ ranges from 1 to $10^{-14}$ in increments of $1/10$, while the tolerances $\epsilon_1$ and $\epsilon_2$ for the MOOD detector to control $\rho$, are respectively set to $10^{-3}$ and $10^{-2}$. We observe that the CATMOOD2-Trap (CM2-Trap), CATMOOD2-Tay (CM2-Tay), CAT2-Trap (C2-Trap), and CAT2-Tay (C2-Tay) methods achieve second-order accuracy both in space and time. 

 These methods demonstrate a {small} computational efficiency on fine grids, with respect to the performance of the Runge-Kutta2 method. In Table~\ref{Broad_table_1} we show the CPU time in seconds for different grid sizes $N$ and $\epsilon$.

Furthermore, as shown in Figure~\ref{Test_Broad_smooth_1}, improving the threshold values $\epsilon_1$ and $\epsilon_2$ results in better computational solutions, especially for situations with smooth initial conditions.

\paragraph{Non-smooth case} The experiment is performed over the interval $[-1, 1]$ using initial conditions that are neither well-prepared nor smooth \cite{PareschiRusso}. 
As first test, we consider the Riemann problem
\begin{equation} 
\label{Broad:IC_disc_RP1}
	(\rho, m, z)
	= \left\{
	\begin{array}{ll}
	\displaystyle (2,1,1) & \mbox { if } x \le 0.2, \\[2mm]
	\displaystyle (1,0.13962,1) & \mbox { if } x > 0.2
	\end{array}\right.
\end{equation}
as initial condition.
\begin{figure}[!ht]
     \centering
     \begin{subfigure}[b]{0.48\textwidth}
         \centering
         \includegraphics[width=\textwidth]{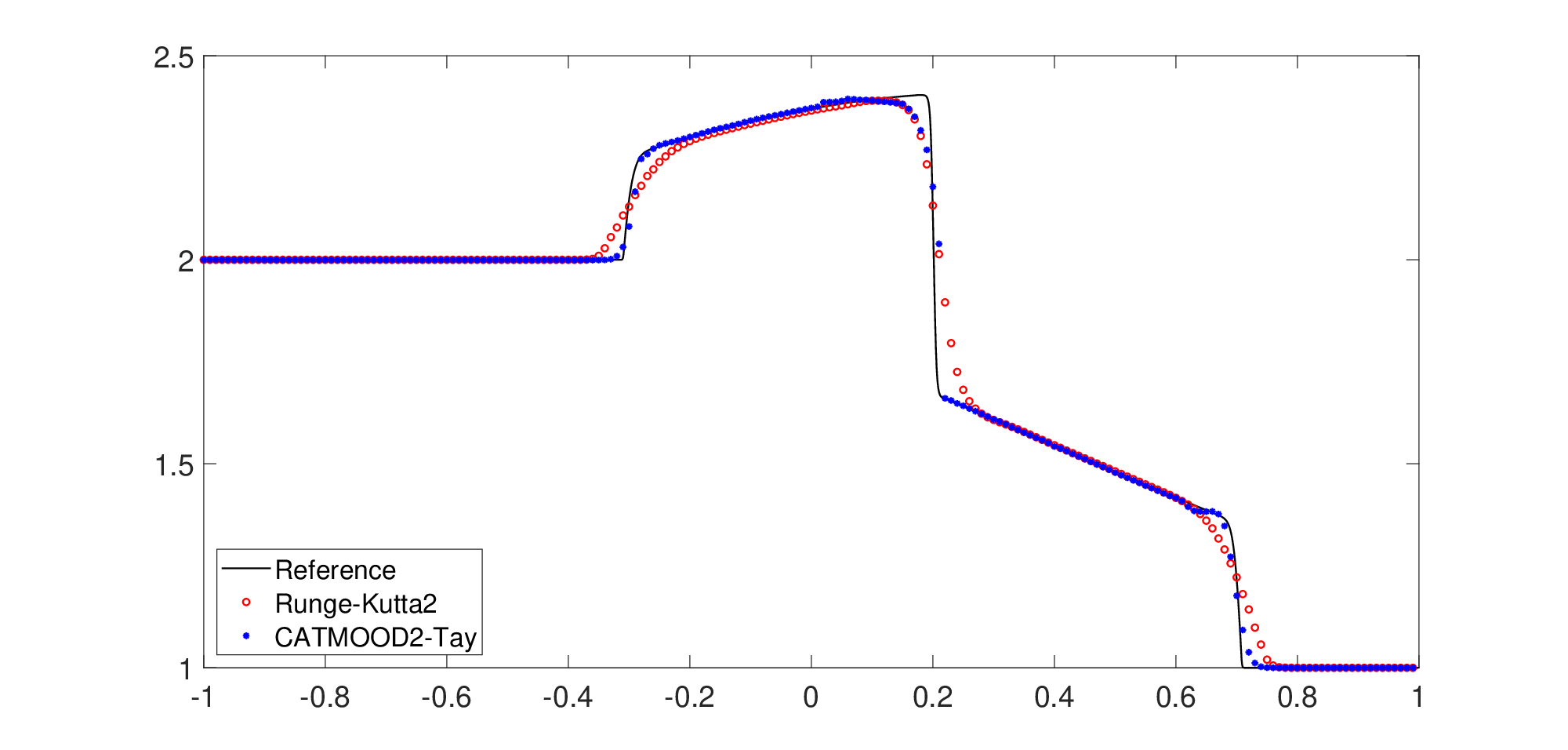}
         \caption{$\rho$}
         \label{Broad:RP1_rho_1}
     \end{subfigure}
     \hfill
     \begin{subfigure}[b]{0.48\textwidth}
         \centering
         \includegraphics[width=\textwidth]{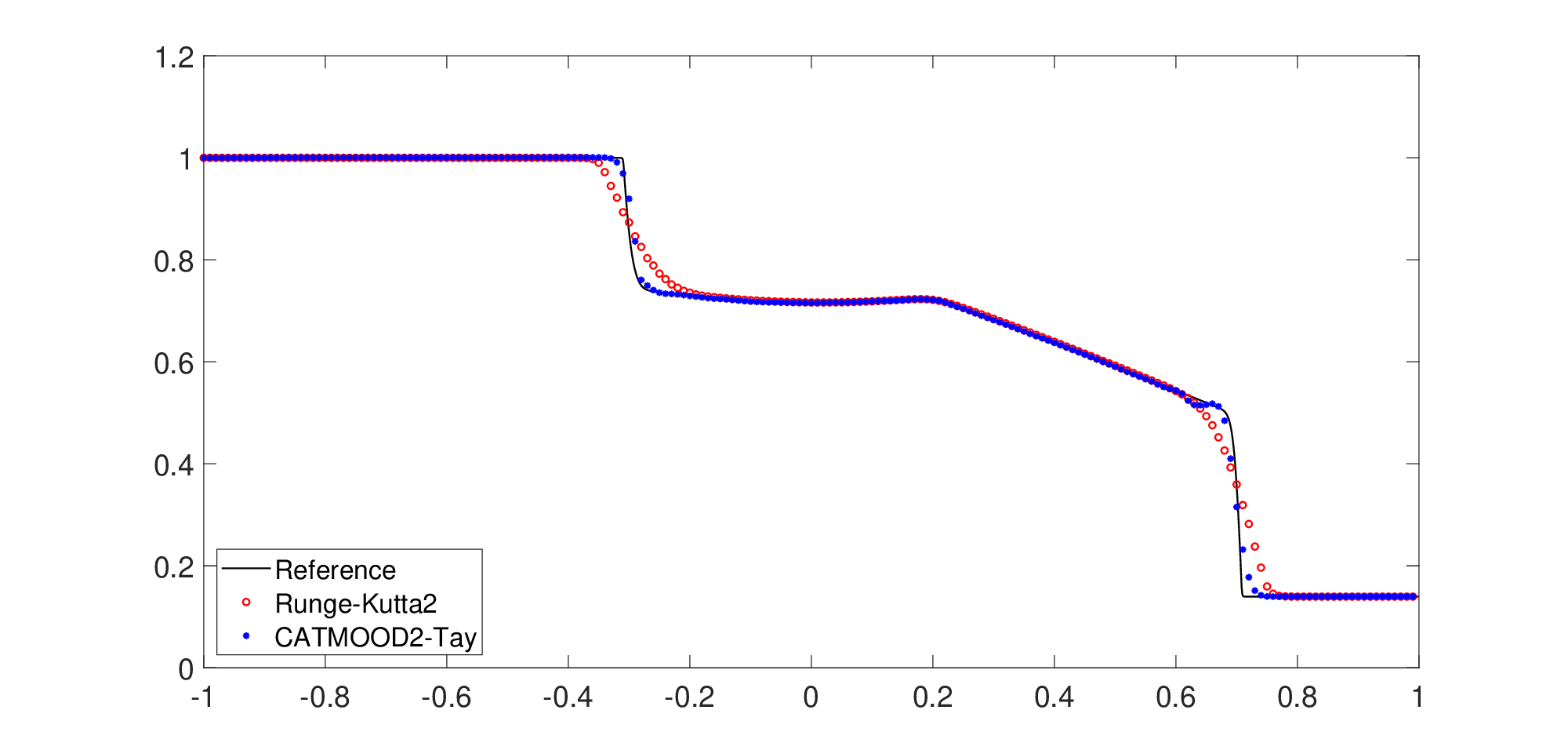}
         \caption{$m$}
         \label{Broad:RP1_m_1}
     \end{subfigure}
     \\
     \begin{subfigure}[b]{0.48\textwidth}
         \centering
         \includegraphics[width=\textwidth]{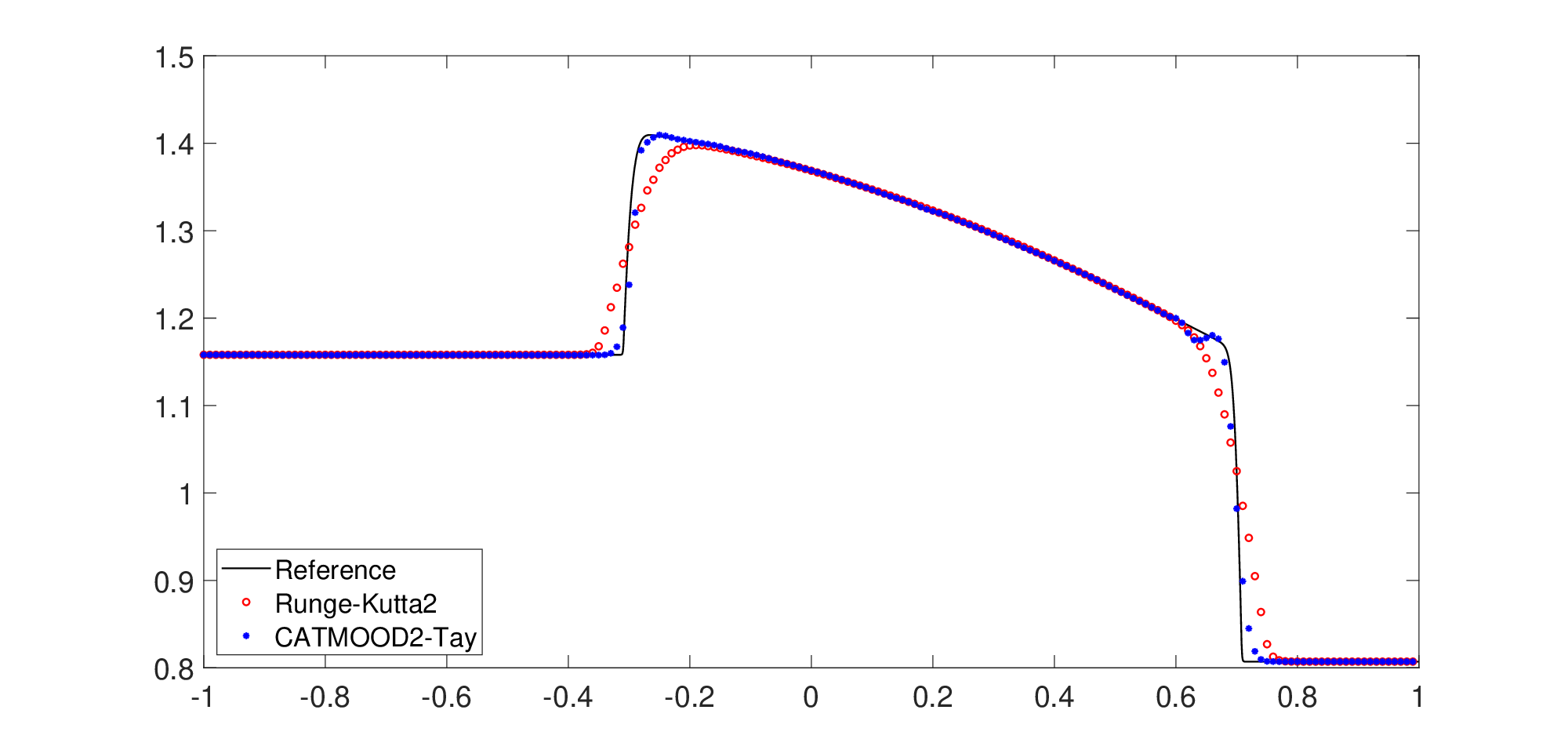}
         \caption{$z$}
         \label{Broad:RP1_z_1}
     \end{subfigure}
     \caption{Broadwell: non-smooth case. Numerical solutions $\rho$ (left), $m$ (central) and $z$ (right) obtained with Runge-Kutta2 and CATMOOD2-Tay at time $t = 0.5$ on the interval $[-1,1]$ with CFL$=0.9$ and $\epsilon = 1$. The tolerances $\epsilon_1$ and $\epsilon_2$ in the MOOD technique are, respectively, set $10^{-4}$ and $10^{-3}.$ The reference solutions have been obtained with Runge-Kutta2 on $2000$ uniform mesh.}
     \label{Broad:RP1_1}
\end{figure}
\begin{figure}[!ht]
     \centering
     \begin{subfigure}[b]{0.48\textwidth}
         \centering
         \includegraphics[width=\textwidth]{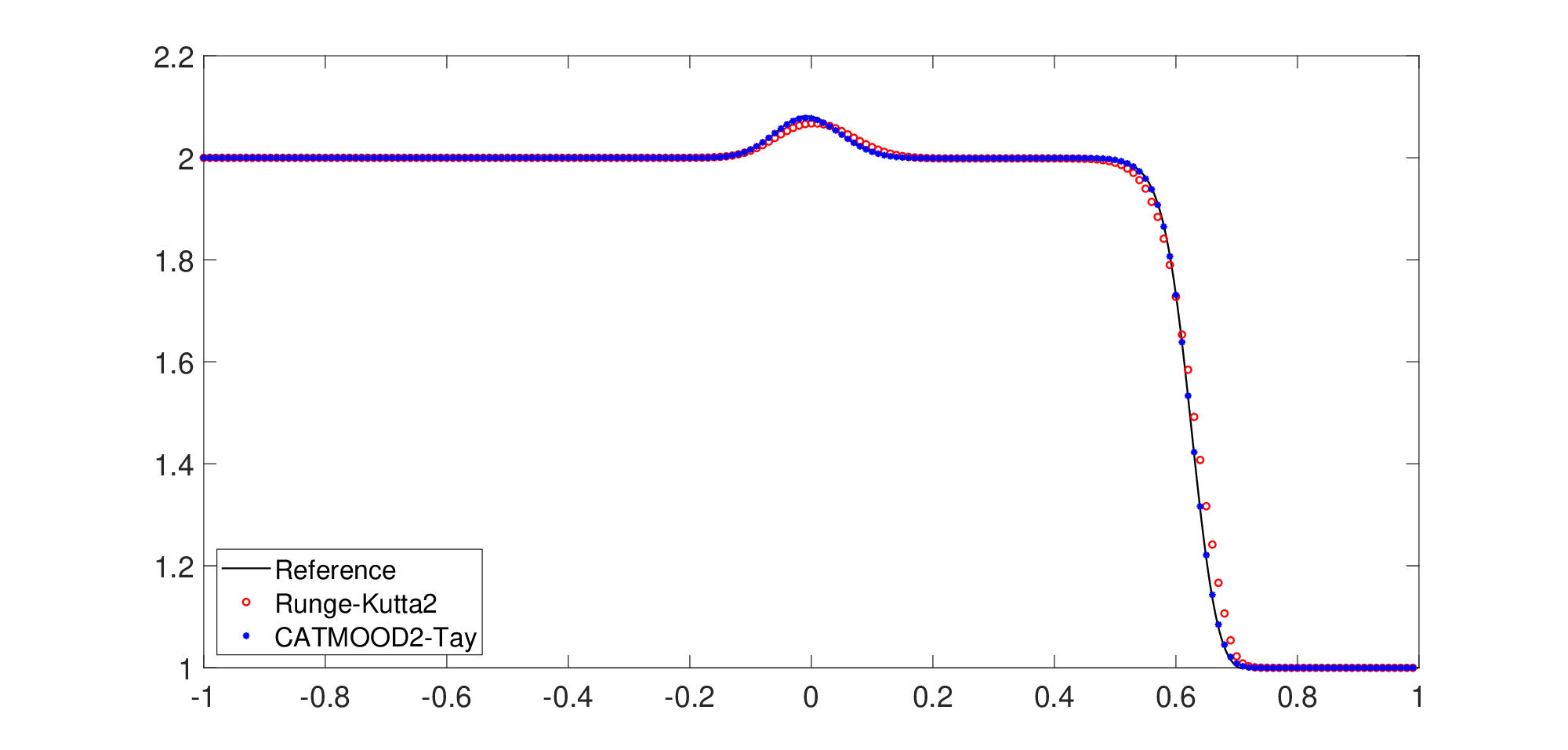}
         \caption{$\rho$}
         \label{Broad:RP1_rho_2}
     \end{subfigure}
     \hfill
     \begin{subfigure}[b]{0.48\textwidth}
         \centering
         \includegraphics[width=\textwidth]{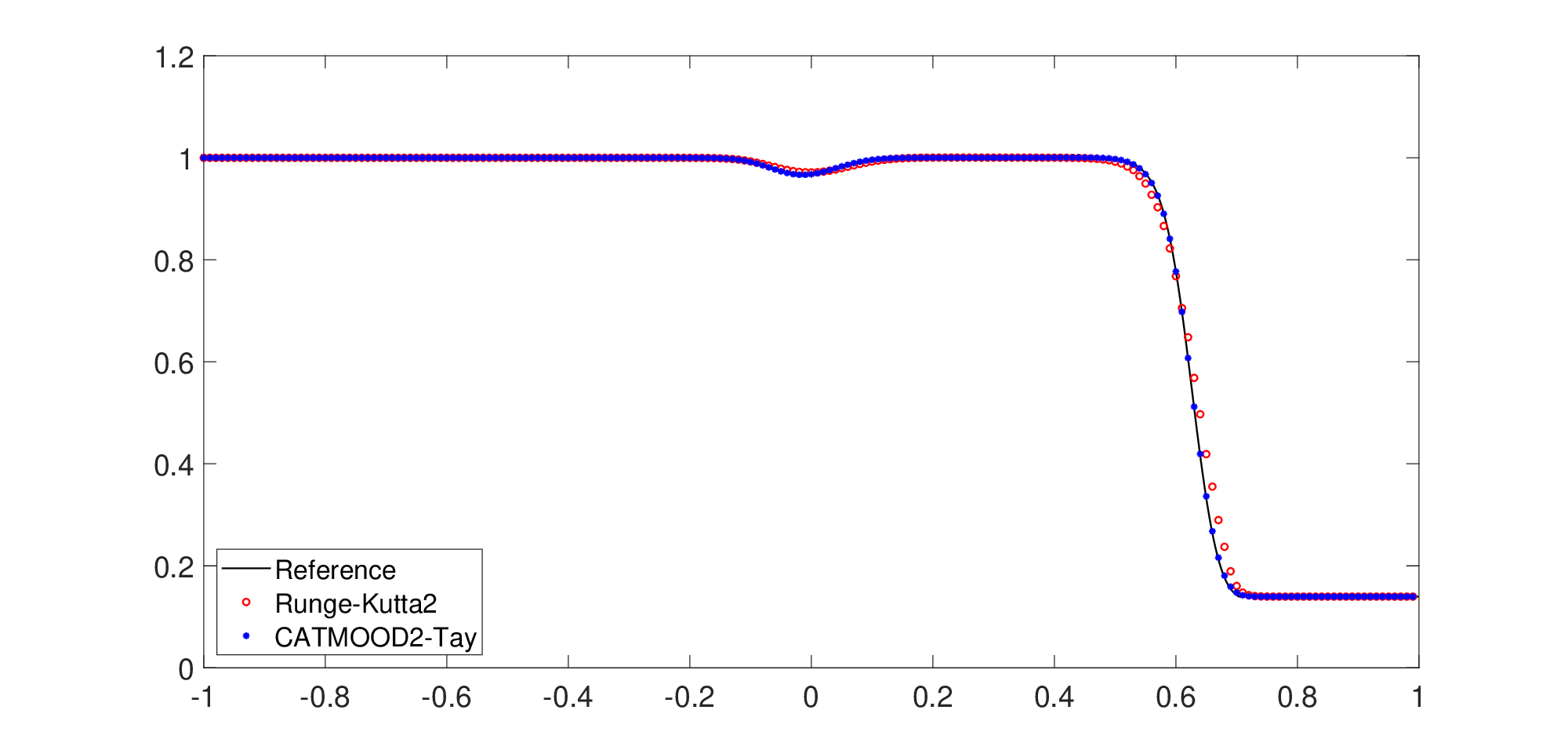}
         \caption{$m$}
         \label{Broad:RP1_m_2}
     \end{subfigure}
     \\
     \begin{subfigure}[b]{0.48\textwidth}
         \centering
         \includegraphics[width=\textwidth]{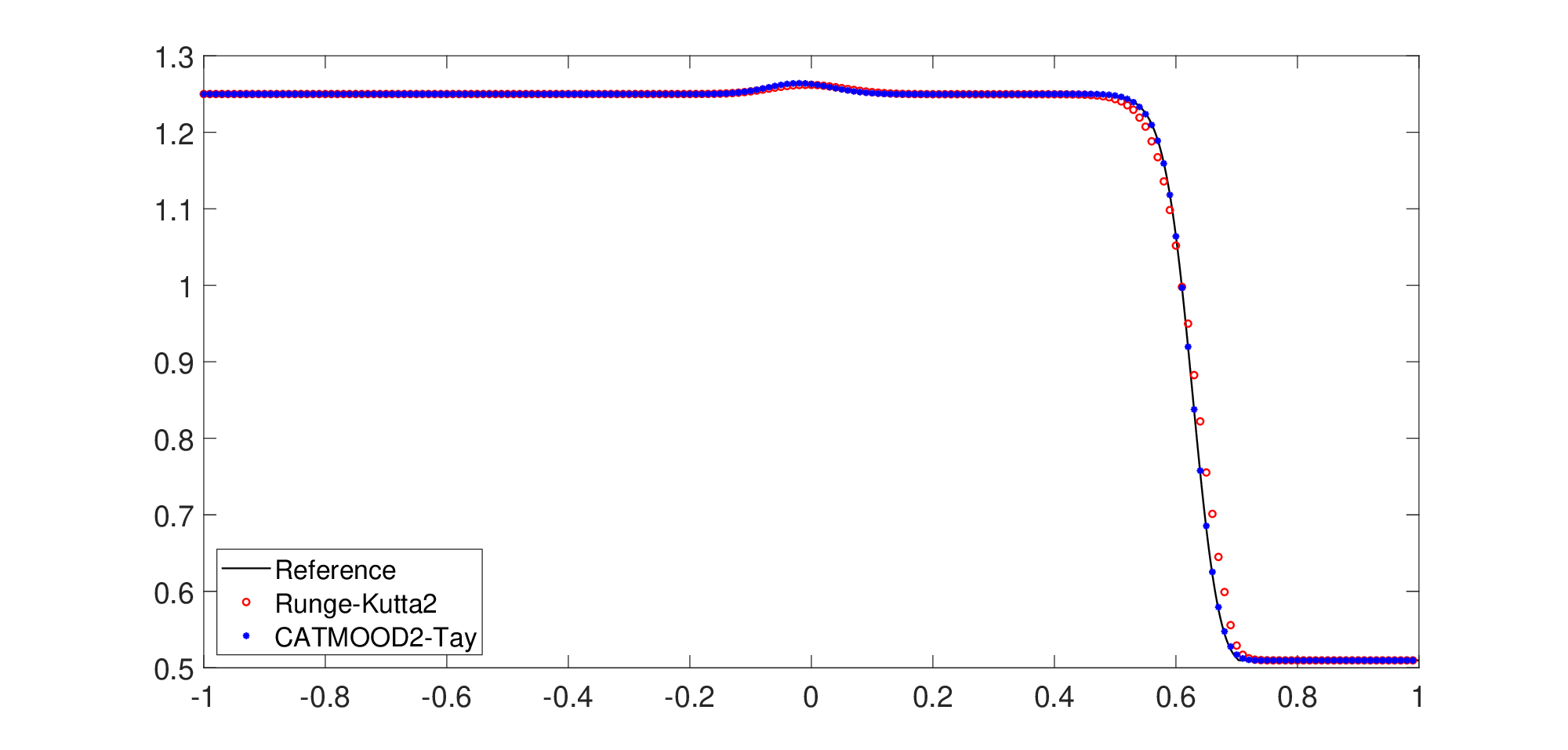}
         \caption{$z$}
         \label{Broad:RP1_z_2}
     \end{subfigure}
     \caption{Broadwell: non-smooth case. Numerical solutions $\rho$ (left), $m$ (central) and $z$ (right) obtained with Runge-Kutta2 and CATMOOD2-Tay at time $t = 0.5$ on the interval $[-1,1]$ with CFL$=0.9$ and $\epsilon = 0.02$. The tolerances $\epsilon_1$ and $\epsilon_2$ in the MOOD technique are, respectively, set $10^{-4}$ and $10^{-3}.$ The reference solutions have been obtained with Runge-Kutta2 on $2000$ uniform mesh.}
     \label{Broad:RP1_2}
\end{figure}
\begin{figure}[!ht]
     \centering
     \begin{subfigure}[b]{0.48\textwidth}
         \centering
         \includegraphics[width=\textwidth]{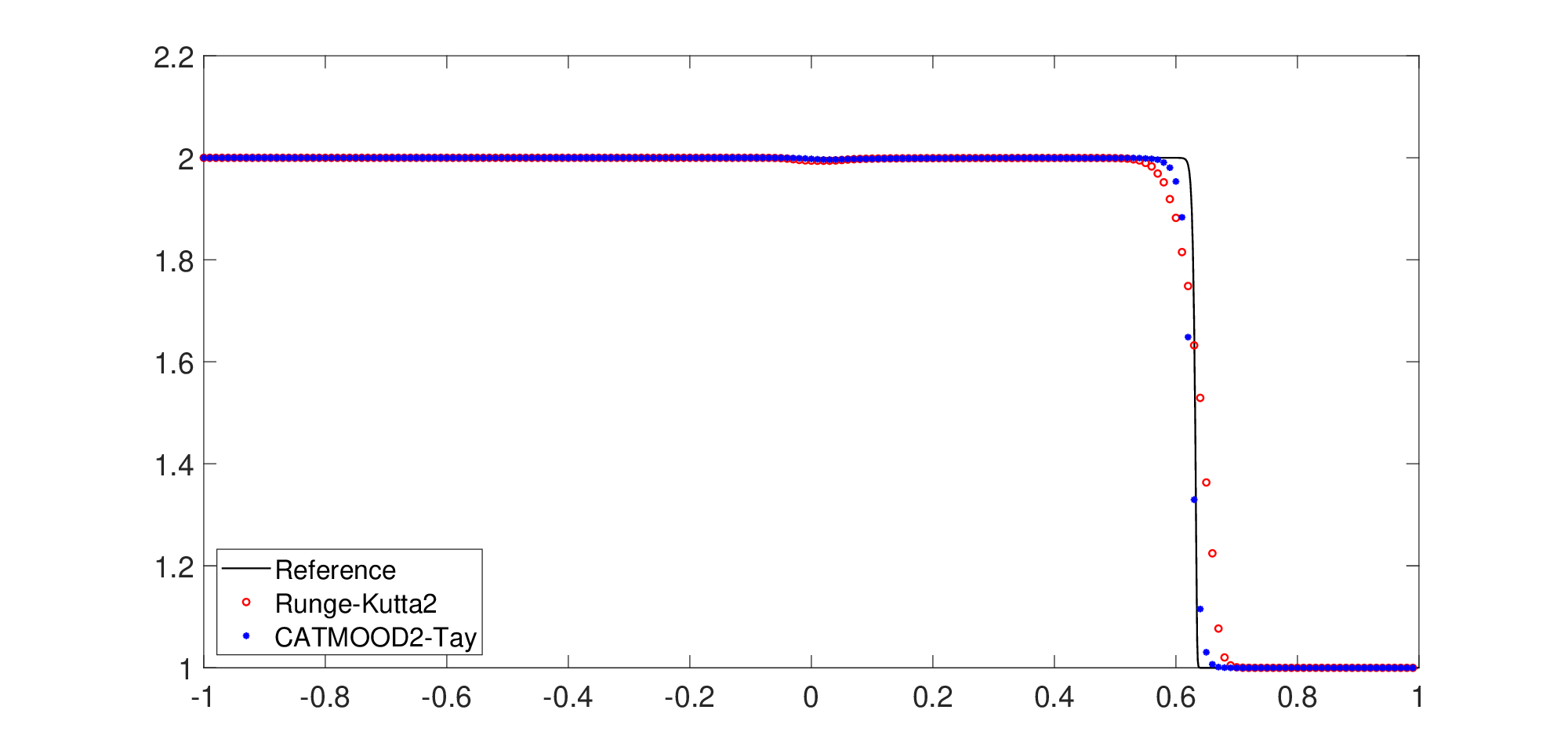}
         \caption{$\rho$}
         \label{Broad:RP1_rho_3}
     \end{subfigure}
     \hfill
     \begin{subfigure}[b]{0.48\textwidth}
         \centering
         \includegraphics[width=\textwidth]{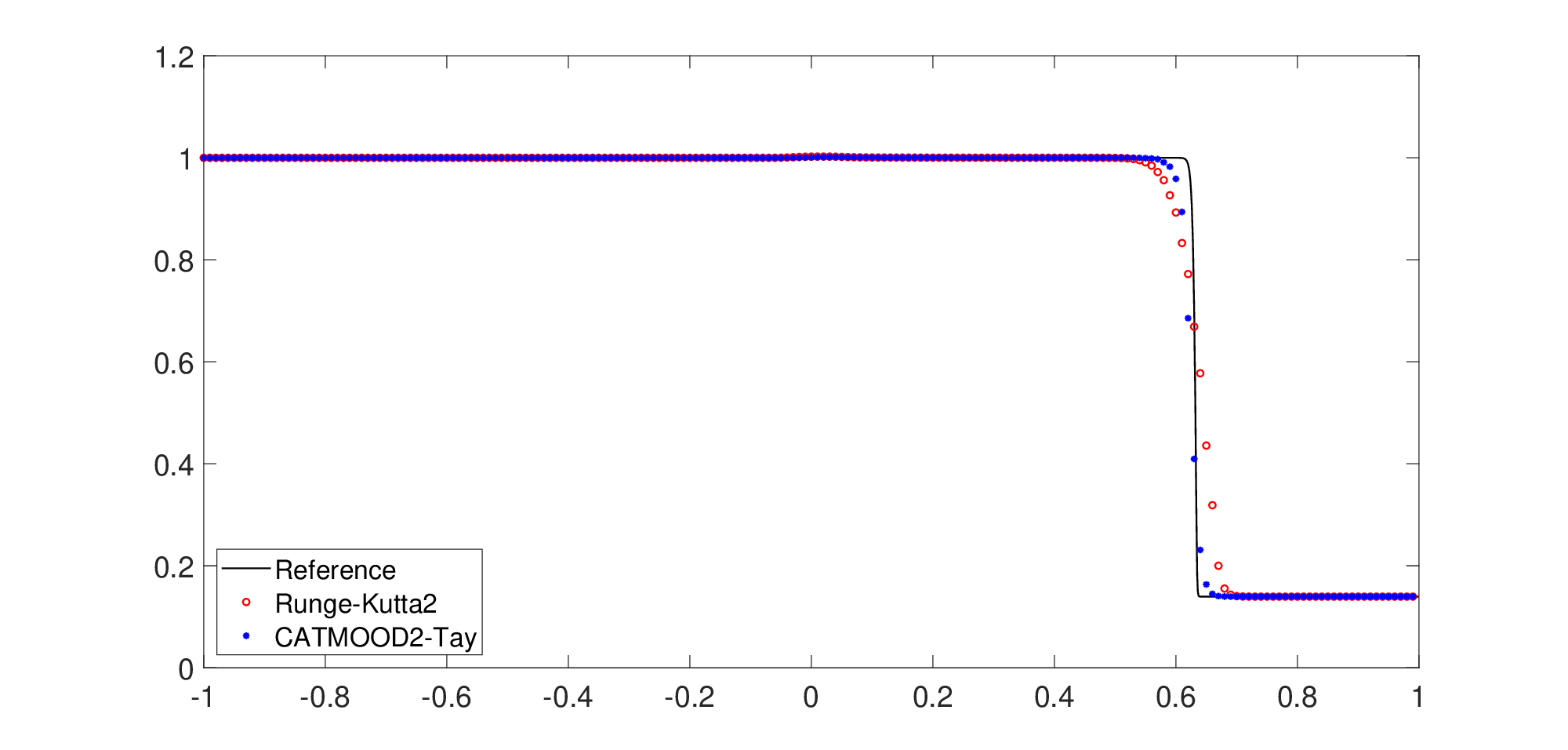}
         \caption{$m$}
         \label{Broad:RP1_m_3}
     \end{subfigure}
     \\
     \begin{subfigure}[b]{0.48\textwidth}
         \centering
         \includegraphics[width=\textwidth]{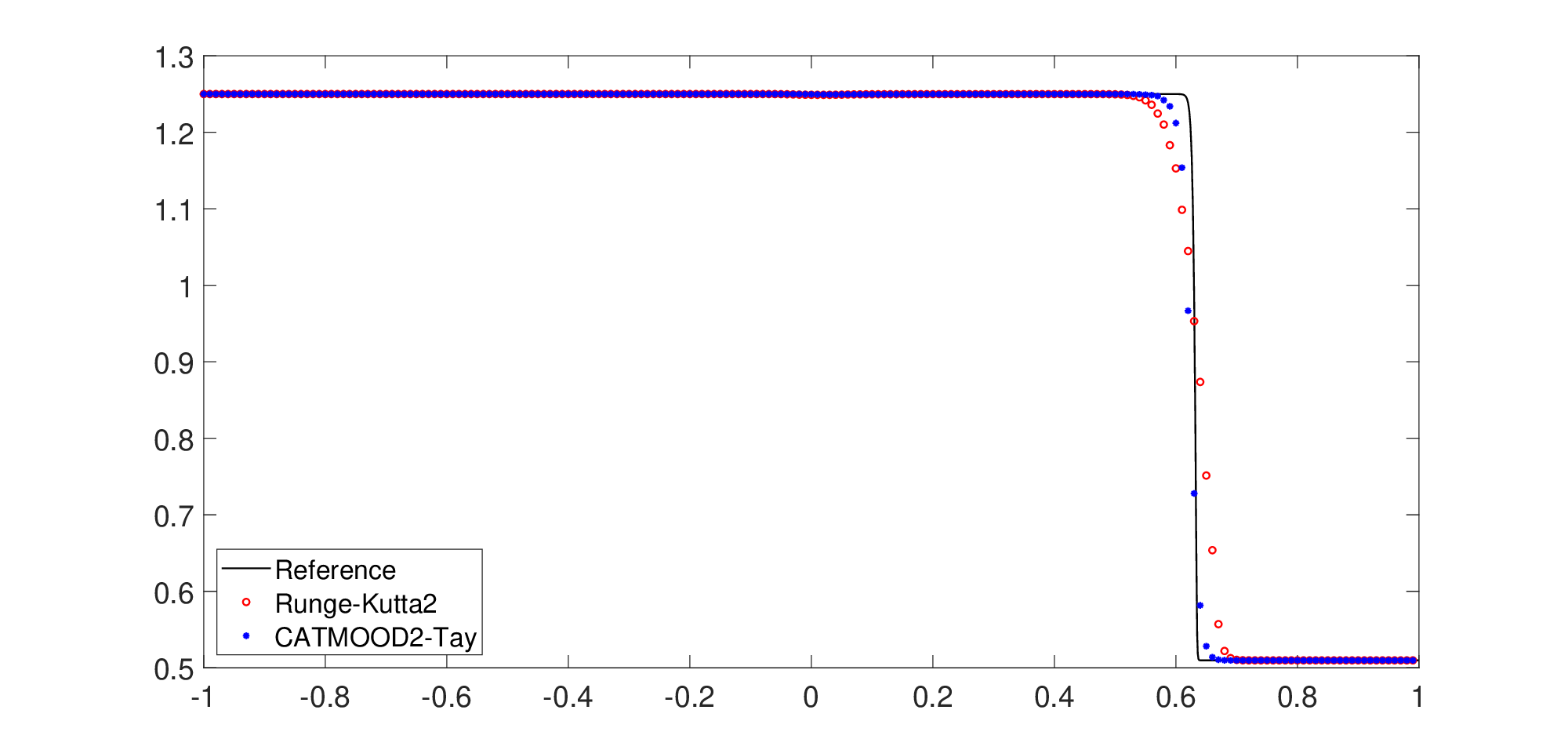}
         \caption{$z$}
         \label{Broad:RP1_z_3}
     \end{subfigure}
     \caption{Broadwell: non-smooth case. Numerical solutions $\rho$ (left), $m$ (central) and $z$ (right) obtained with Runge-Kutta2 and CATMOOD2-Tay at time $t = 0.5$ on the interval $[-1,1]$ with CFL$=0.9$ and $\epsilon = 10^{-8}$. The tolerances $\epsilon_1$ and $\epsilon_2$ in the MOOD technique are, respectively, set $10^{-4}$ and $10^{-3}.$ The reference solutions have been obtained with Runge-Kutta2 on $2000$ uniform mesh.}
     \label{Broad:RP1_3}
\end{figure}
\begin{figure}[!ht]
     \centering
     \begin{subfigure}[b]{0.49\textwidth}
         \centering
         \includegraphics[width=\textwidth]{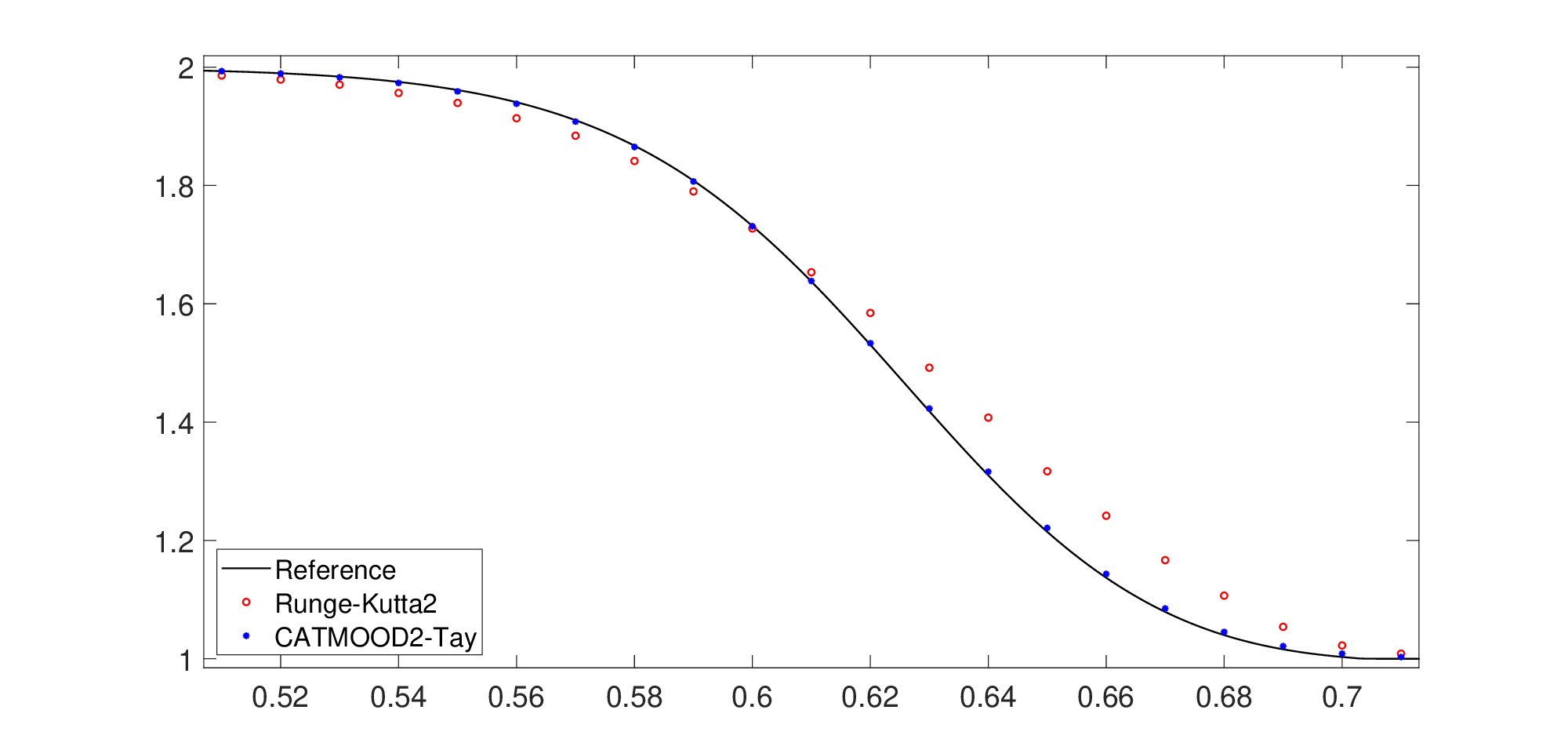}
         \caption{$\rho$ when $\epsilon = 0.02.$}
         \label{Broad:RP1_rho_2_zoom}
     \end{subfigure}
     \hfill
     \begin{subfigure}[b]{0.49\textwidth}
         \centering
         \includegraphics[width=\textwidth]{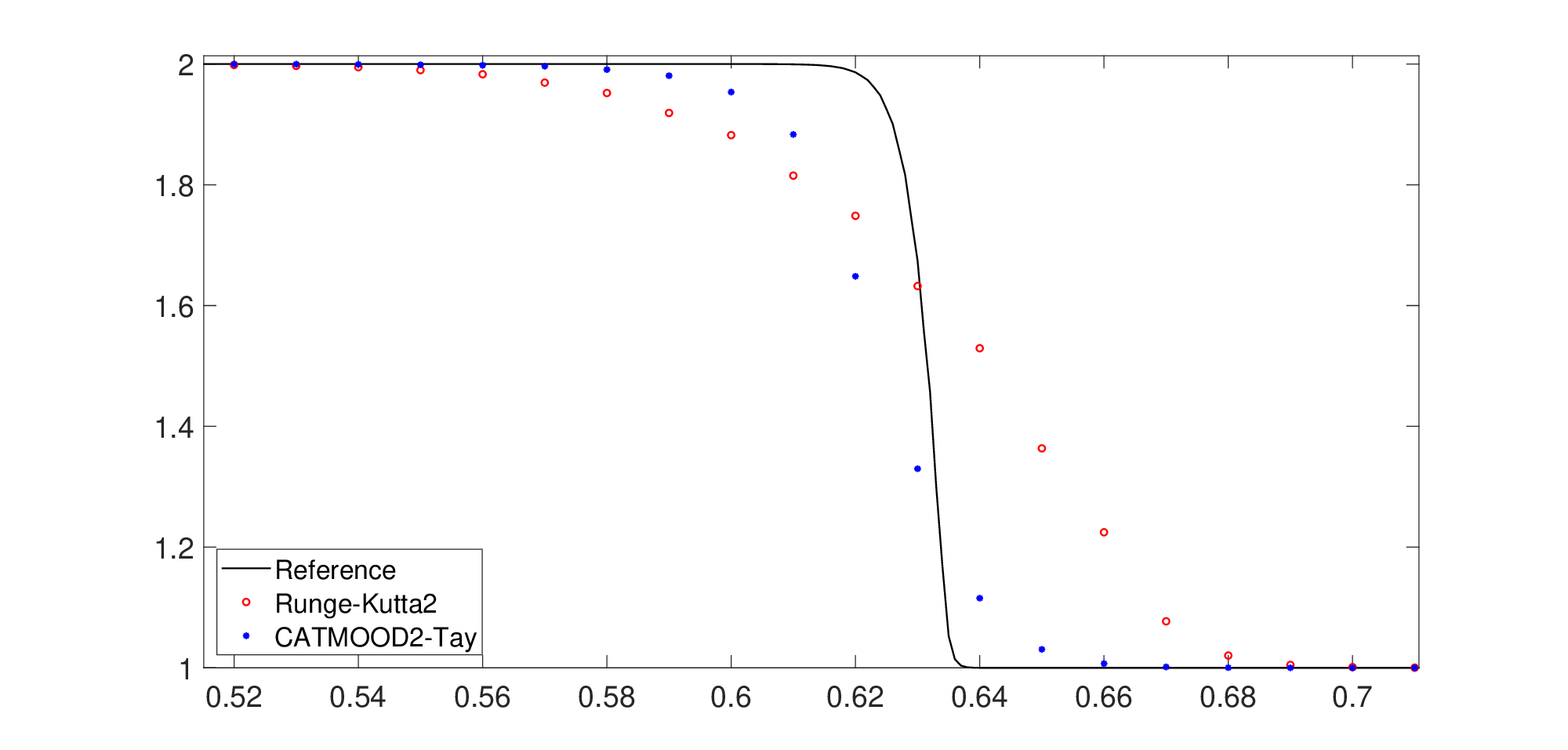}
         \caption{$\rho$ when $\epsilon = 10^{-8}.$}
         \label{Broad:RP1_rho_3_zoom}
     \end{subfigure}
     \caption{Broadwell: non-smooth case. Zoom of the numerical solutions for $\rho$ close the shock obtained with Runge-Kutta2 and CATMOOD2-Tay at time $t = 0.5$ on the interval $[-1,1]$ with CFL$=0.9$ and $\epsilon$ values $0.02$ (left) and $10^{-8}$ (right). The tolerances $\epsilon_1$ and $\epsilon_2$ in the MOOD technique are, respectively, set $10^{-4}$ and $10^{-3}.$ The reference solutions have been obtained with Runge-Kutta2 on $2000$ uniform mesh.}
     \label{Broad:RP1_2_3_zoom}
\end{figure}
The computational domain $[\alpha,\beta]$ spans from $-1$ to $1$, the final time $t_{\text{fin}}=0.5$, CFL$=0.9$, $N=200$ points and Neumann zero conditions are imposed at boundaries. The $\epsilon$ values are drawn from the set $\{1,0.02, 10^{-8}\}$, and finally, the tolerance values $\epsilon_1$ and $\epsilon_2$ are set to $10^{-4}$ and $10^{-3}$ respectively. The reference solutions are obtained using the second-order semi-implicit RK method on a grid consisting of 2000 points.

\begin{figure}[!ht]
     \centering
     \begin{subfigure}[b]{0.48\textwidth}
         \centering
         \includegraphics[width=\textwidth]{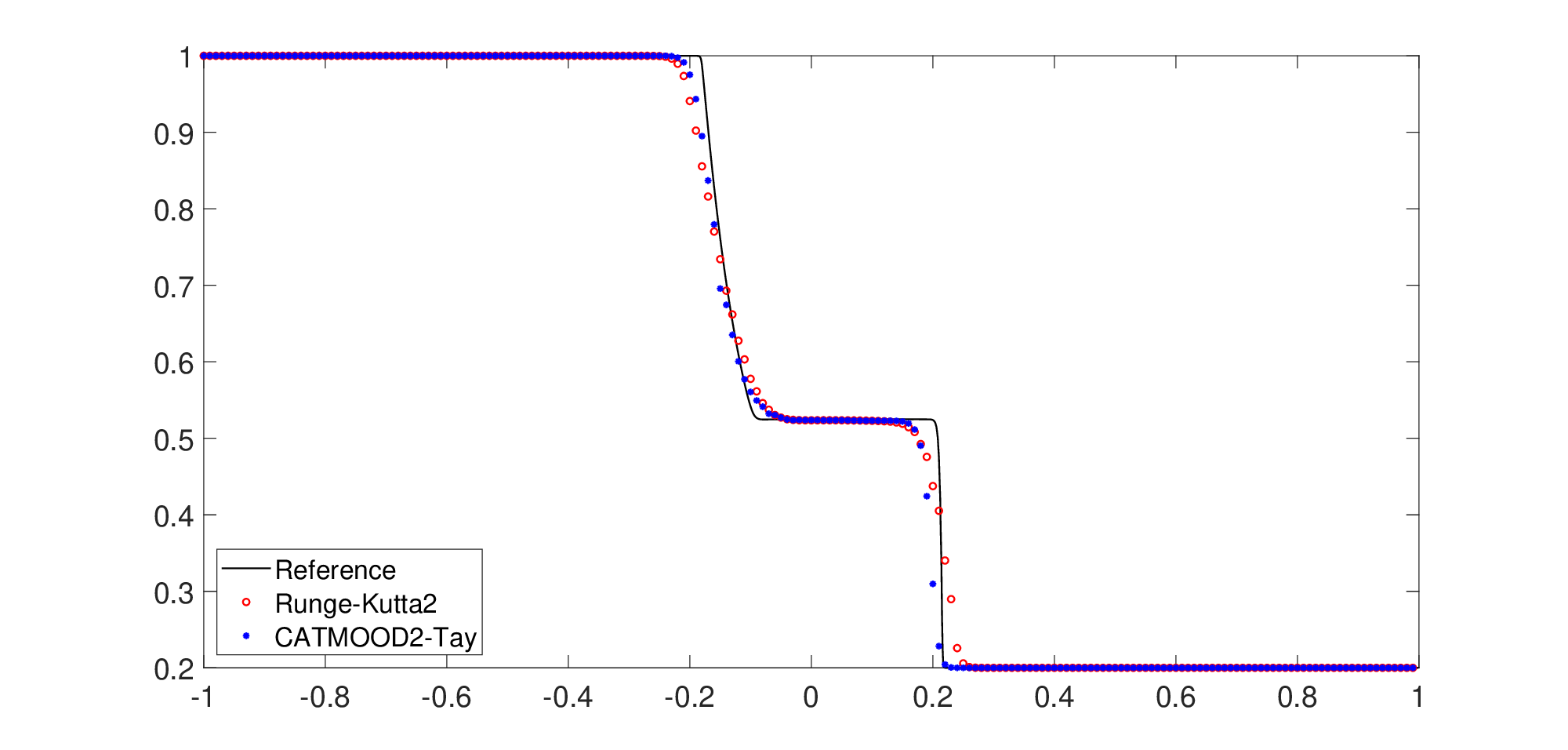}
         \caption{$\rho$}
         \label{Broad:RP2_rho}
     \end{subfigure}
     \hfill
     \begin{subfigure}[b]{0.48\textwidth}
         \centering
         \includegraphics[width=\textwidth]{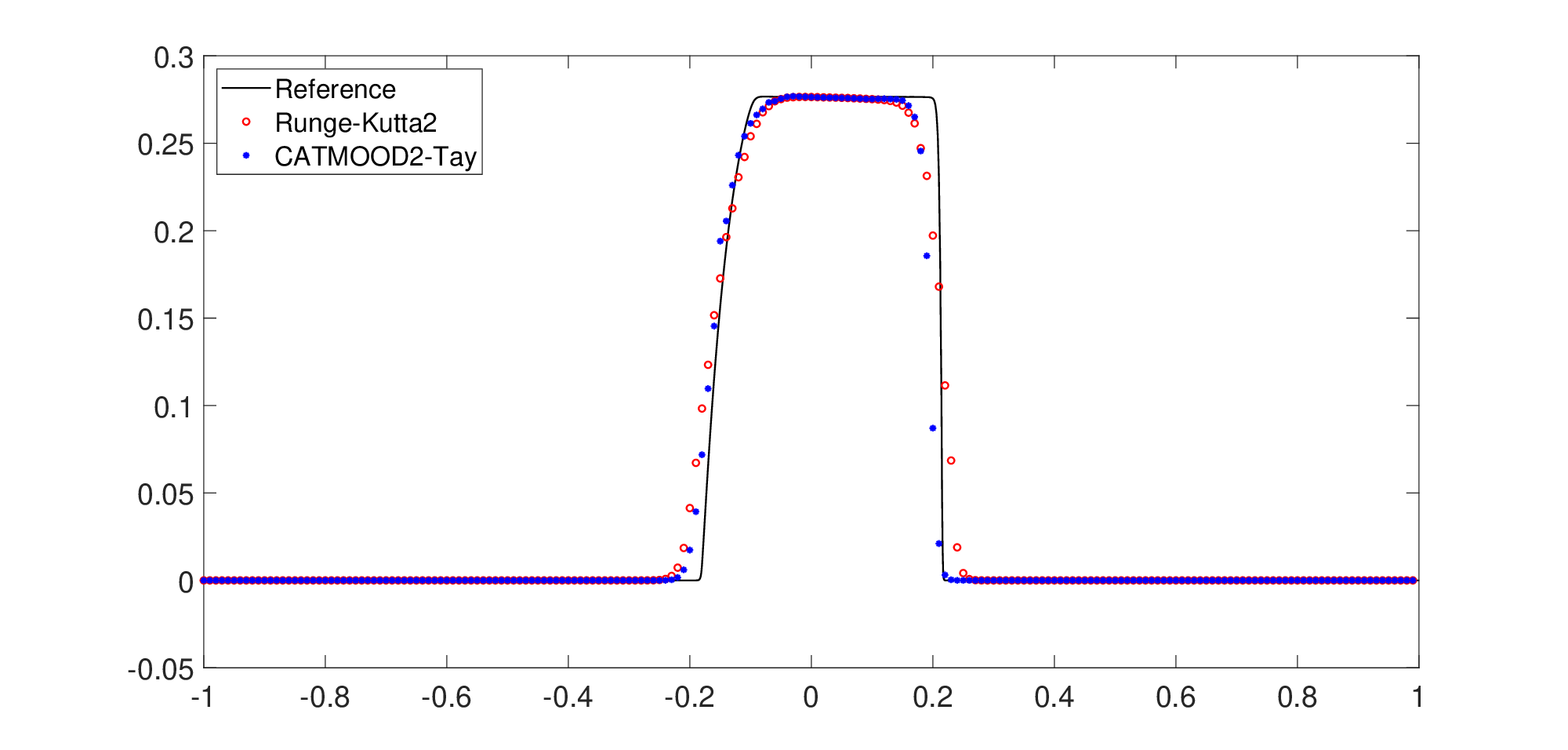}
         \caption{$m$}
         \label{Broad:RP2_m}
     \end{subfigure}
     \\
     \begin{subfigure}[b]{0.48\textwidth}
         \centering
         \includegraphics[width=\textwidth]{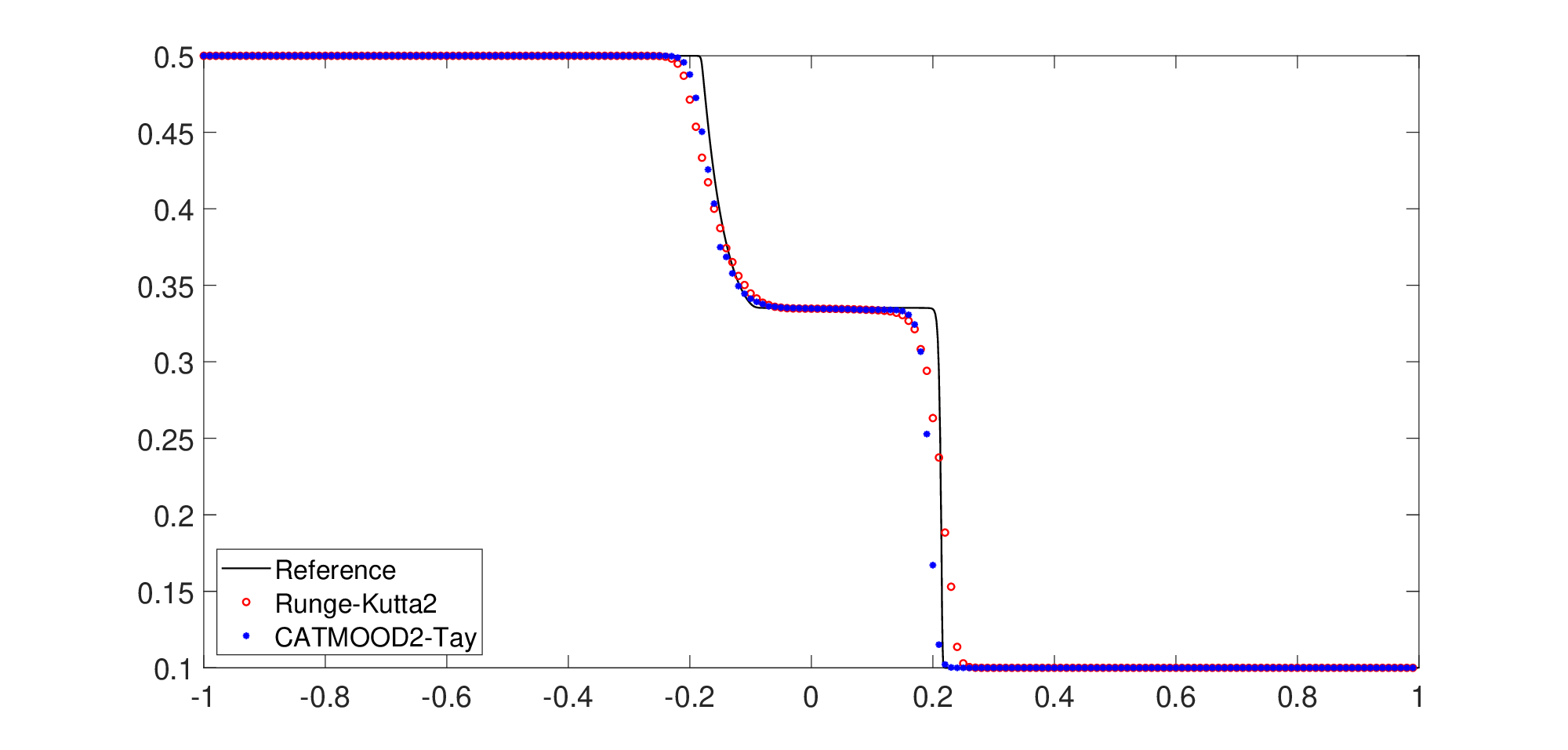}
         \caption{$z$}
         \label{Broad:RP2_z}
     \end{subfigure}
     \caption{Broadwell: non-smooth case. Numerical solutions $\rho$ (left), $m$ (central) and $z$ (right) obtained with Runge-Kutta2 and CATMOOD2-Tay at time $t = 0.25$ on the interval $[-1,1]$ with CFL$=0.9$ and $\epsilon = 10^{-8}$. The tolerances $\epsilon_1$ and $\epsilon_2$ in the MOOD technique are, respectively, set $10^{-4}$ and $10^{-3}.$ The reference solutions have been obtained with Runge-Kutta2 on $2000$ uniform mesh.}
     \label{Broad:RP2_1}
\end{figure}

Here we used two different methods, Runge-Kutta2 and CATMOOD2-Tay, and try them out with different $\epsilon$ values.  In Figures~\ref{Broad:RP1_1}-\ref{Broad:RP1_3}, we display the variables $\rho$, $m$, and $z$. Figure~\ref{Broad:RP1_1} shows a big difference between the two methods close the shocks where CATMOOD2-Tay scheme is much more accurate compared to Runge-Kutta2.

On the contrary, in Figures~\ref{Broad:RP1_2}-\ref{Broad:RP1_3}, the difference between the two methods isn't as strong. However, in Figure~\ref{Broad:RP1_2_3_zoom} from a zoom of the numerical solutions, (for the variable $\rho$), we can see that CATMOOD2-Tay performs better solution compared to Runge-Kutta2.

\begin{remark}
    Given that the initial condition \eqref{Broad:IC_disc_RP1} is not well-prepared, and taking into account the stability analysis discussed in Section~\ref{sec:ode_like_stability}, it becomes evident that the semi-implicit-type CAT2-Trap scheme does not perform well in the limit case (it's not AP). Consequently, this scheme is not suitable for application in this particular case when $\epsilon \to 0$. 

\end{remark}

For the second test, we examine the Riemann problem described by the following initial condition:
\begin{equation} 
\label{Broad:IC_disc_RP2}
	(\rho, m, z)
	= \left\{
	\begin{array}{ll}
	\displaystyle (1,0,1) & \mbox { for } x \le 0, \\[2mm]
	\displaystyle (0.2,0,1) & \mbox { for } x > 0
	\end{array}\right.
\end{equation}

 The computational domain $[\alpha,\beta]$ extends from -1 to 1, with a final simulation time $t_{\text{fin}}=0.25$. A CFL number of 0.9 is applied, and the grid consists of $N=200$ points. Neumann zero conditions are enforced at the boundaries. The value of $\epsilon$ is set to $10^{-8}$, and the tolerance values $\epsilon_1$ and $\epsilon_2$ are established at $10^{-4}$ and $10^{-3}$ respectively. It's worth noting that reference solutions are obtained using the second-order semi-implicit RK method with a grid of 2000 points.

Figure~\ref{Broad:RP2_1} provides a representation of the numerical solution for the variables $\rho$ (left), $m$ (center), and $z$ (right). These solutions are generated using the Runge-Kutta2 and CATMOOD2-Tay schemes.

\subsection{Euler with heat transfer model}
This experiment uses initial conditions that are unprepared and discontinuous, as described in
\cite{jin1995runge}.
\begin{figure}[!ht]
     \centering
     \begin{subfigure}[b]{0.48\textwidth}
         \centering
         \includegraphics[width=\textwidth]{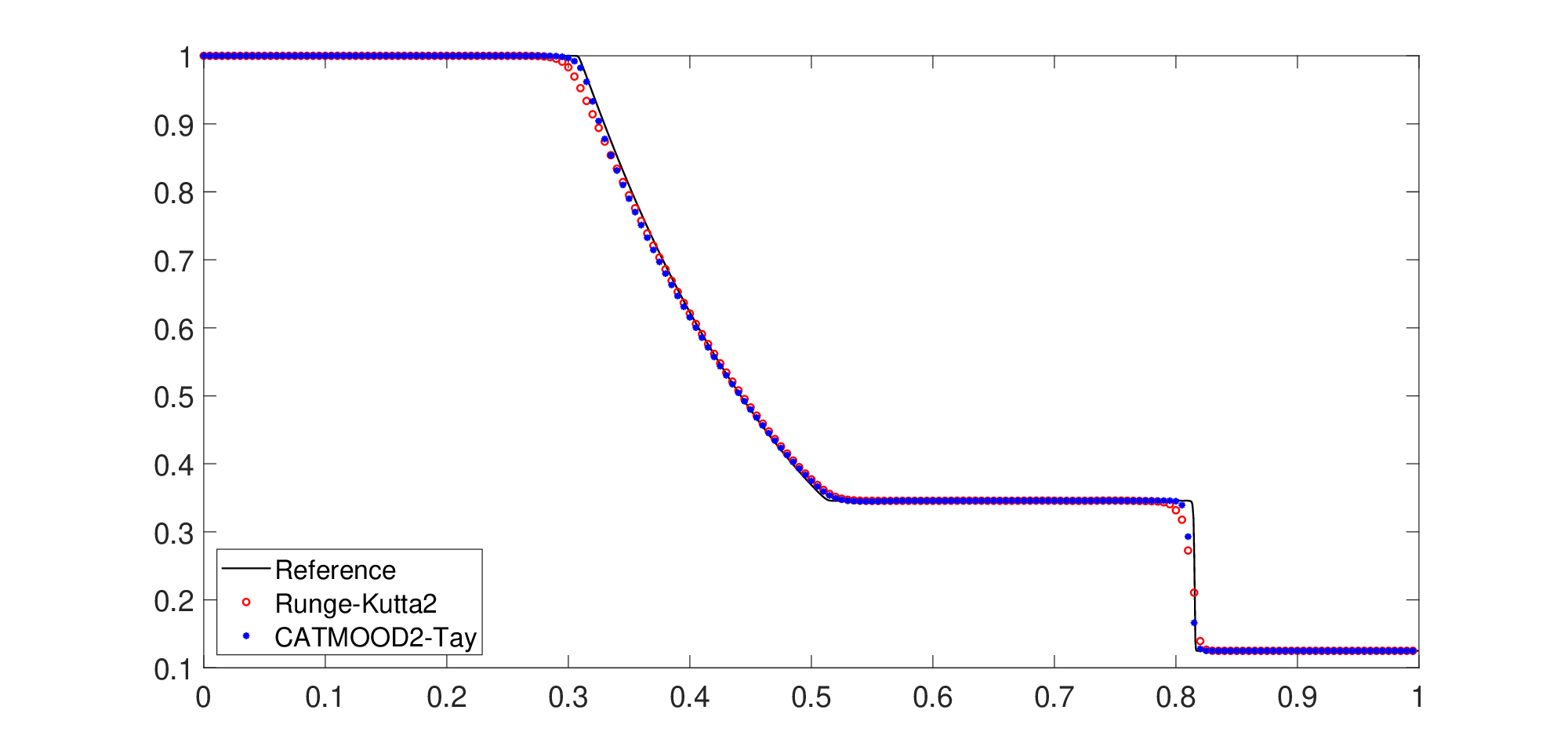}
         \caption{$\rho$}
         \label{EHeat:RP4_rho_1}
     \end{subfigure}
     \hfill
     \begin{subfigure}[b]{0.48\textwidth}
         \centering
         \includegraphics[width=\textwidth]{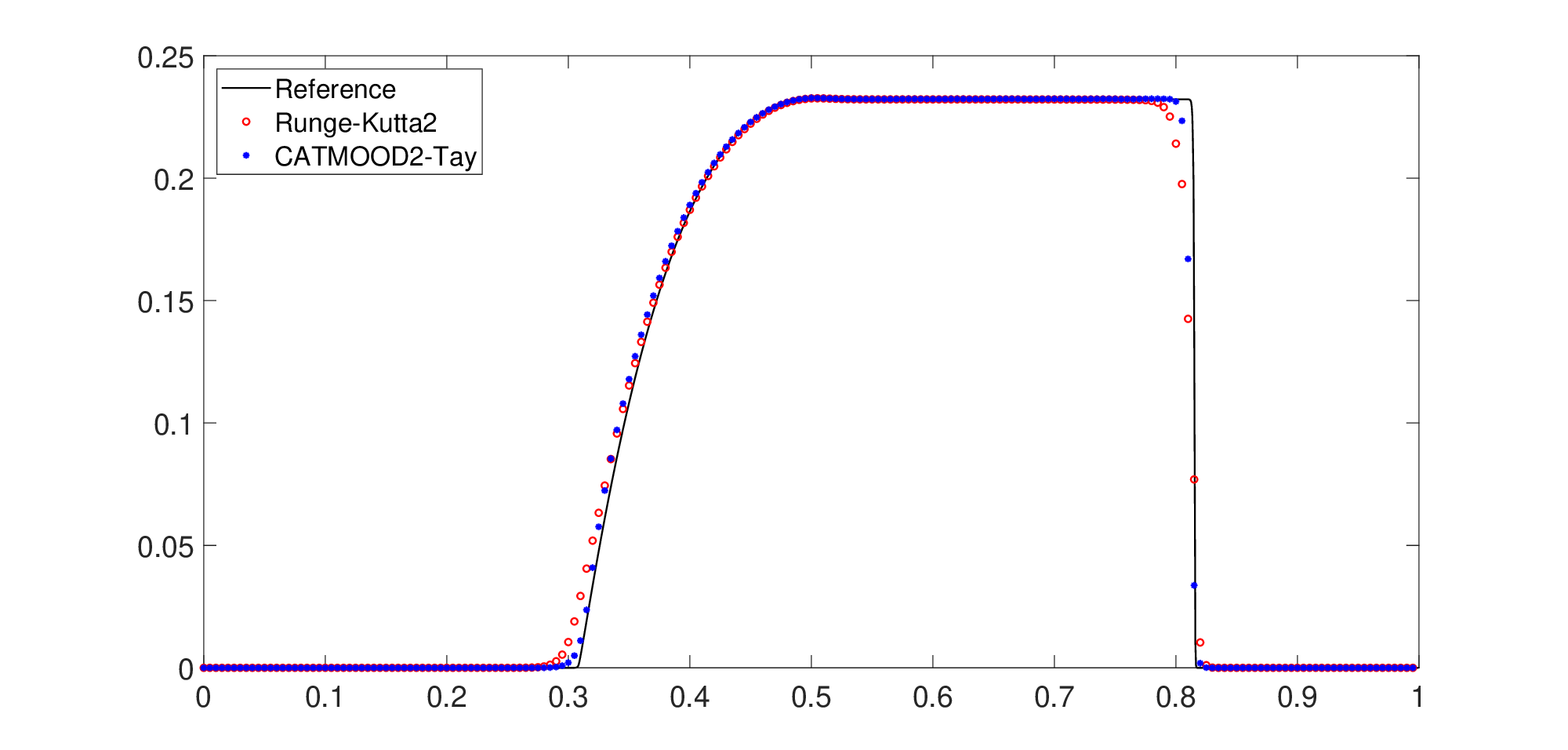}
         \caption{$m$}
         \label{EHeat:RP4_m_1}
     \end{subfigure}
     \\
     \begin{subfigure}[b]{0.48\textwidth}
         \centering
         \includegraphics[width=\textwidth]{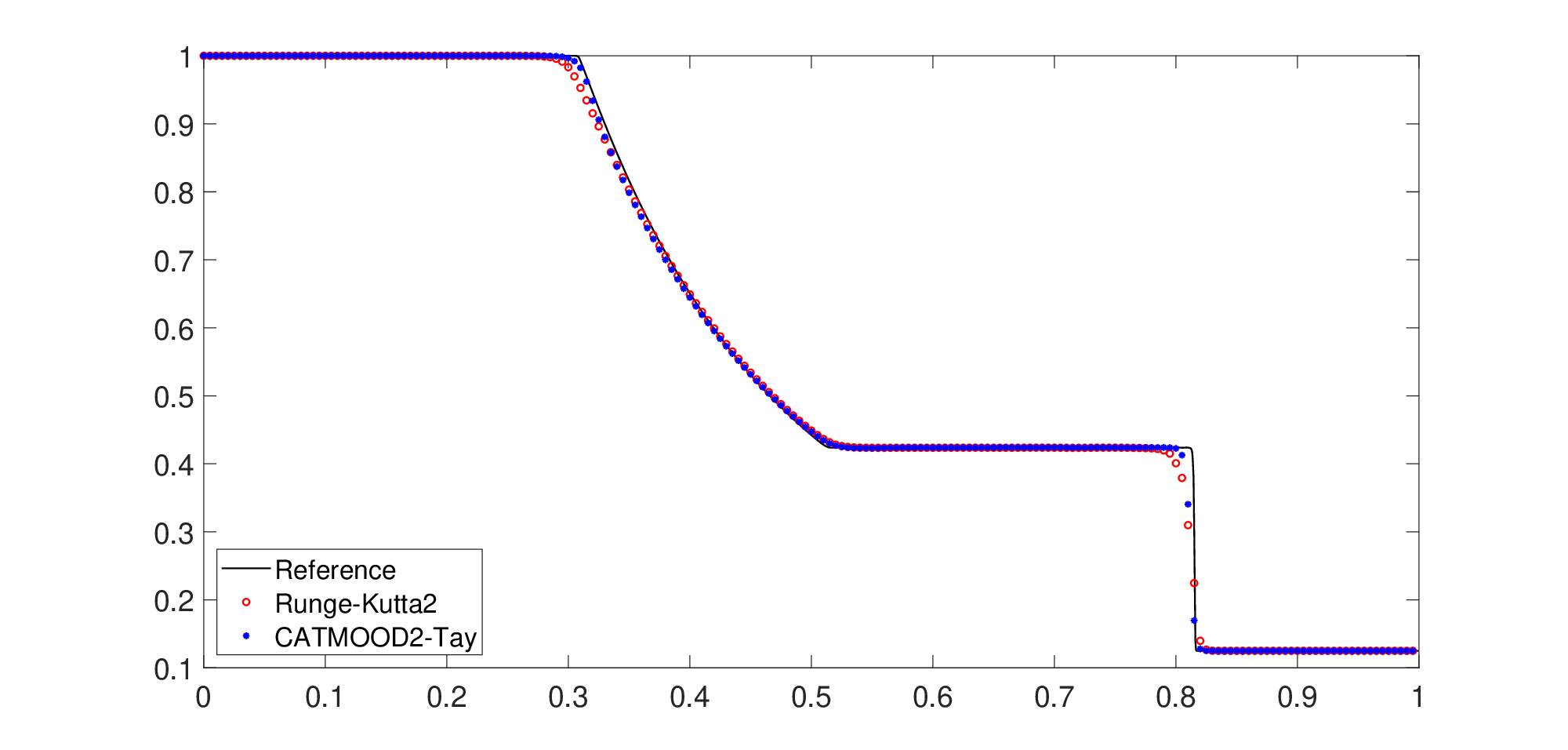}
         \caption{$z$}
         \label{EHeat:RP4_z_1}
     \end{subfigure}
     \caption{Euler with heat transfer: non-smooth scenario. The numerical solutions for $\rho$ (left), $m$ (center), and $z$ (right) are computed using Runge-Kutta2 and CATMOOD2-Tay at time $t = 0.3$ within the interval $[0,1]$. A CFL number of $0.7$ is employed for CATMOOD2-Tay, while Runge-Kutta2 uses a CFL of $0.4$. The value of $\epsilon$ is set to $10^{-8}$. In the MOOD technique, the tolerances $\epsilon_1$ and $\epsilon_2$ are established at $10^{-4}$ and $10^{-3}$ respectively. It's important to note that reference solutions are obtained using Runge-Kutta2 with a uniform mesh of 2000 points.}
     \label{EHeat:RP4_1}
\end{figure}
\begin{figure}[!ht]
     \centering
     \begin{subfigure}[b]{0.48\textwidth}
         \centering
         \includegraphics[width=\textwidth]{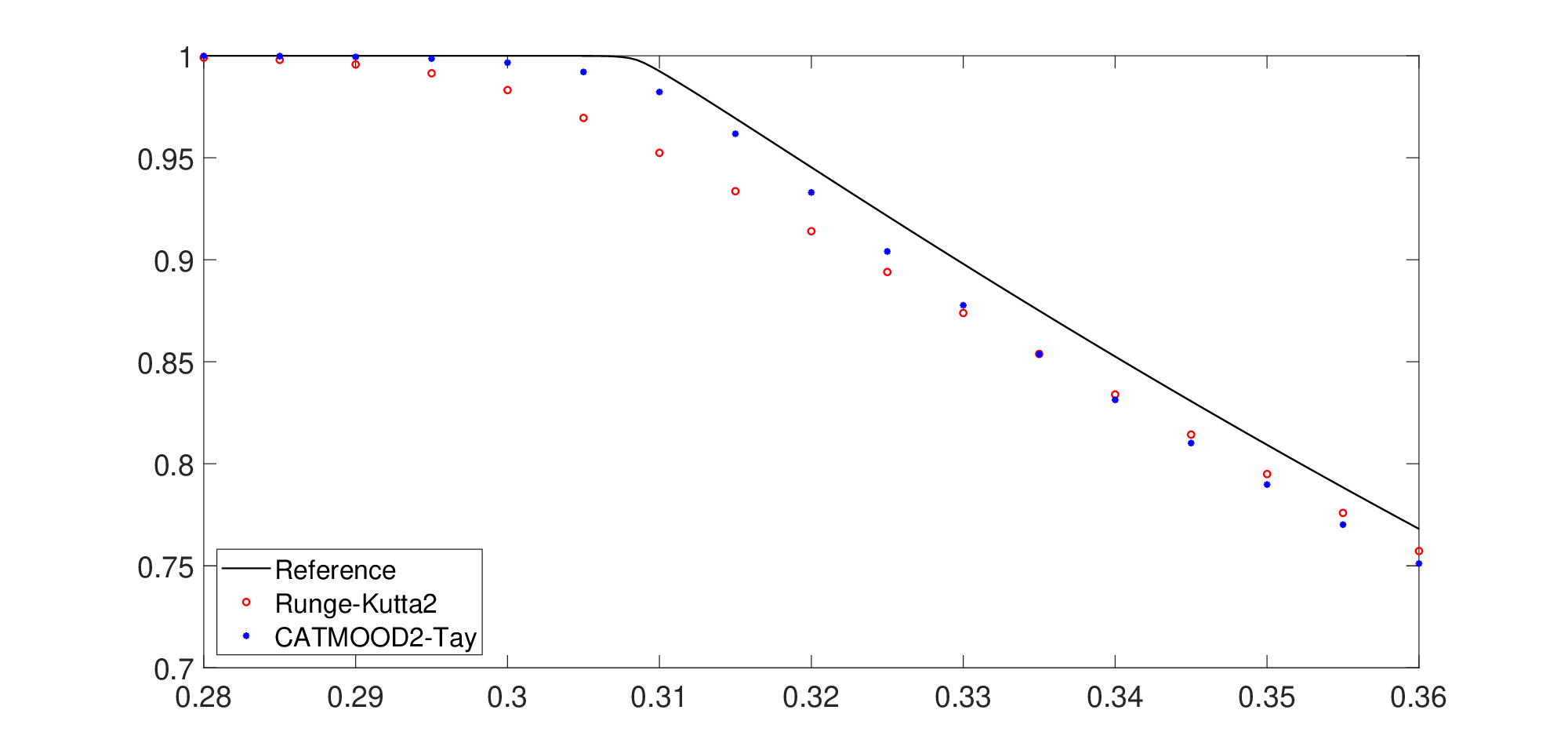}
         \caption{$\rho$}
         \label{EHeat:RP4_rho_1_zoom}
     \end{subfigure}
     \hfill
     \begin{subfigure}[b]{0.48\textwidth}
         \centering
         \includegraphics[width=\textwidth]{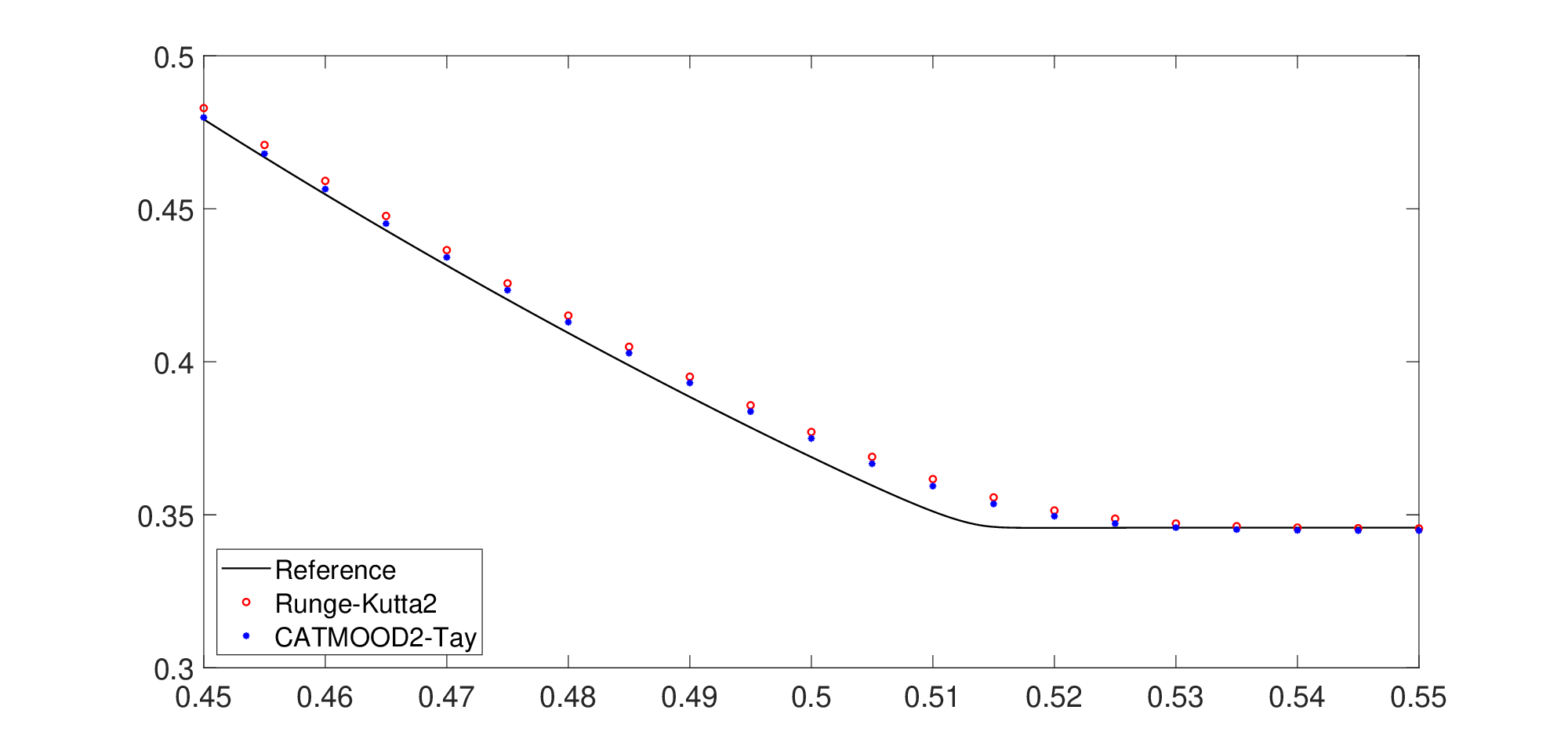}
         \caption{$m$}
         \label{EHeat:RP4_rho_2_zoom}
     \end{subfigure}
     \\
     \begin{subfigure}[b]{0.48\textwidth}
         \centering
         \includegraphics[width=\textwidth]{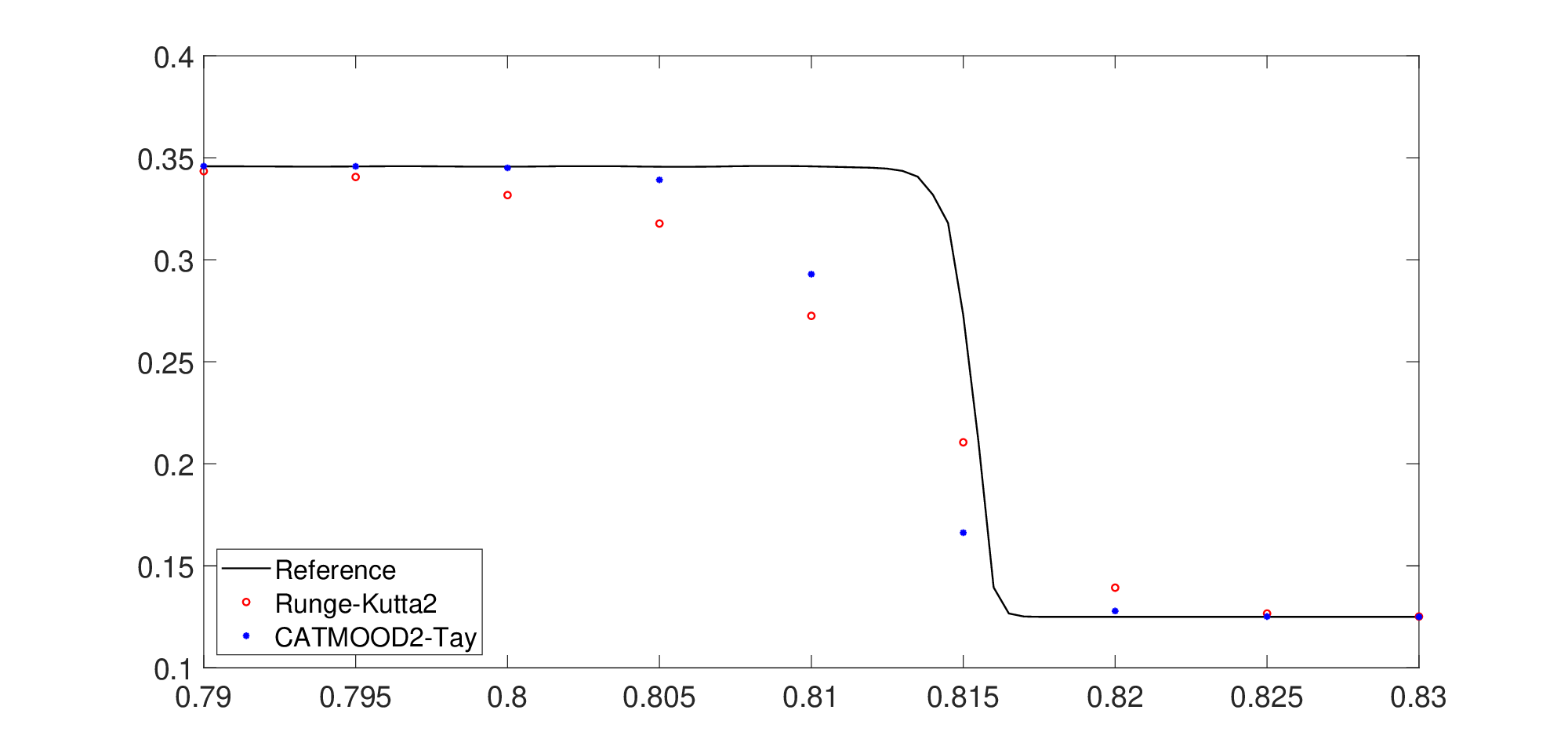}
         \caption{$z$}
         \label{EHeat:RP4_rho_3_zoom}
     \end{subfigure}
     \caption{Euler with heat transfer: non-smooth scenario. A close examination of the numerical solutions for $\rho$ at three critical points is presented. These solutions are computed at time $t = 0.3$ within the interval $[0,1]$, employing both Runge-Kutta2 and CATMOOD2-Tay. The CFL number for CATMOOD2-Tay is set at $0.7$, while Runge-Kutta2 utilizes a CFL of $0.4$. The value of $\epsilon$ is specified as $10^{-8}$. Within the MOOD technique, the tolerances $\epsilon_1$ and $\epsilon_2$ are defined at $10^{-4}$ and $10^{-3}$ respectively. The reference solutions are generated using Runge-Kutta2 with a uniform mesh consisting of 2000 points.}
     \label{EHeat:RP4_rho_zoom}
\end{figure}
We consider the computational domain $[\alpha,\beta]$ extends from $0$ to $1$, with a final simulation time $t_{\text{fin}}=0.3$ and Neumann zero boundary conditions. For the CATMOOD2-Tay scheme, a CFL number of $0.7$ is applied, while for the second-order semi-implicit RK scheme a CFL of $0.4$ is used. The grid is composed of $N=200$ points. The value of $\epsilon$ is set at $10^{-8}$, and the tolerance levels $\epsilon_1$ and $\epsilon_2$ are established at $10^{-4}$ and $10^{-3}$ respectively. Notably, the reference solution is computed using the second-order semi-implicit RK method on a grid containing 2000 points with a CFL of $0.4$.

Figure~\ref{EHeat:RP4_1} illustrates the numerical solutions for density ($\rho$), momentum ($m$), and energy per unit mass ($\rho E$) obtained through the Runge-Kutta2 and CATMOOD2-Tay schemes.  Initially, these solutions look quite similar. However, if we zoom in on key areas of density ($\rho$) in Figure~\ref{EHeat:RP4_rho_zoom}, it becomes clear that the CATMOOD2-Tay solution is slightly more accurate than the one produced by Runge-Kutta2.

Anyway, it's important to mention that there's a difference in computational cost between these two schemes. Specifically, the CPU time in seconds for Runge-Kutta2 is $0.0405$, whereas for CATMOOD2-Tay takes to $0.0505$. This variation primarily arises from the utilization of the MOOD technique and the flexibility of adjusting certain parameters. By modifying $\epsilon_1$ and $\epsilon_2$ to $10^{-3}$ and $10^{-2}$, it's possible to achieve a solution that's a bit less accurate but takes a similar amount of time to compute as Runge-Kutta2.

\section{Conclusions and perspectives}
In this article, we introduced a semi-implicit-type technique designed for hyperbolic equation models featuring stiff source terms. We extended the CAT2 scheme, originally designed for conservation laws, to balance laws with stiff sources, using two different approaches. We combined this approach with the local order-adaptive MOOD (Multi-dimensional Optimal Order Detection) technique.

To evaluate the effectiveness of these schemes, we selected some models from the existing literature. These models allowed us to assess various properties, including the accuracy of the schemes under smooth conditions, their ability to preserve the asymptotic preserving property of the numerical solutions, and their essential non-oscillatory behavior when confronted with shocks or large gradients. We conducted comparative analyses, naming the schemes CATMOOD2-Trap and CATMOOD2-Tay based on the source term type reconstruction. 

Our results demonstrated that the semi-implicit CATMOOD2 schemes outperform in terms of both accuracy and computational efficiency, when dealing with well-prepared and smooth initial conditions compared to the Runge-Kutta2 scheme. In cases involving unprepared and non-smooth initial conditions, similar behaviors were observed between CATMOOD2-Tay and Runge-Kutta2.

However, as anticipated from its derivation in Sec.~\ref{sec:ode_like_stability} and in the numerical tests section, CATMOOD2-Trap is an unstable scheme as $\epsilon\rightarrow0$ (the stiff parameter), i.e., the scheme is not asymptotic preserving. 

Anyway,  in scenarios featuring well-prepared and non-smooth initial conditions, the solutions derived from CATMOOD-Trap or CATMOOD-Tay exhibited greater or similar accuracy near shocks compared to Runge-Kutta2, consistently across varying $\epsilon$ values.

{It is worth noting that the performance of the semi-implicit-type CATMOOD schemes is depends significantly on the MOOD technique. By fine-tuning the parameters of this technique can leads to solutions that show a slight reduction in accuracy for significant computational cost savings. Therefore, it is crucial to consider both the experimental setup and its associated parameters.}

In conclusion, the current versions of the semi-implicit-type CATMOOD schemes has certain limitations, including second-order accuracy, applicability solely to balance laws featuring a stiff source term, and the absence of well-balanced properties. These aspects, coupled with the extension to 2D scenarios, represent important areas for further exploration in our future research in this field.

\section*{Acknowledgements} 
This research has received funding from the European Union’s NextGenerationUE – Project: Centro Nazionale HPC, Big Data e Quantum Computing, “Spoke 1” (No. CUP E63C22001000006). E. Macca was partially supported by GNCS No. CUP E53C22001930001 Research Project “Metodi numerici per problemi differenziali multiscala: schemi di alto ordine, ottimizzazione, controllo”. E. Macca and S. Boscarino would like to thank the Italian Ministry of Instruction, University and Research (MIUR) to support this research with funds coming from PRIN Project 2022  (2022KA3JBA, entitled “Advanced numerical methods for time dependent parametric partial differential euqations and applications”). Sebastiano Boscarino has been supported for this work from Italian Ministerial grant PRIN 2022 PNRR “FIN4GEO: Forward and Inverse Numerical Modeling of hydrothermal systems in volcanic regions with application to geothermal energy exploitation.”, No. P2022BNB97. E. Macca and S. Boscarino are members of the INdAM Research group GNCS.

\bibliographystyle{plain}
\bibliography{biblio}

\begin{thebibliography}{10}

\bibitem{MaccaExner}
S.~Avgerinos, M.J. Castro, E.~Macca, and G.~Russo.
\newblock A semi-implicit finite volume method for the {E}xner model of
  sediment transport.
\newblock {\em Journal of Computational Physics, Major revision}, 2023.

\bibitem{ImexBosca}
S.~Boscarino.
\newblock Error analysis of {IMEX} {R}unge–{K}utta methods derived from
  {D}ifferential-{A}lgebraic systems.
\newblock {\em SIAM J. Numer. Anal.}, 45(4):1600--1621, 2006.

\bibitem{Boscarino-Filbet}
S.~Boscarino, F.~Filbet, and G.~Russo.
\newblock High order semi-implicit schemes for time dependent partial
  differential equations.
\newblock {\em Journal of Scientific Computing}, 68(8):975--1001, 2016.

\bibitem{BoscarinoRusso}
S.~Boscarino, L.~Pareschi, and G.~Russo.
\newblock A unified imex runge-kutta approach for hyperbolic systems with
  multiscale relaxation.
\newblock {\em SIAM J. Numer. Anal.}, 55(4):2017, 2085-2109.

\bibitem{broadwell1964shock}
J.E. Broadwell.
\newblock Shock structure in a simple discrete velocity gas.
\newblock {\em The Physics of Fluids}, 7(8):1243--1247, 1964.

\bibitem{MCPR2022}
H.~Carrillo, E.~Macca, C.~Parés, and G.~Russo.
\newblock Well-{B}alanced {A}daptive {C}ompact {A}pproximate {T}aylor methods
  for systems of balance laws.
\newblock {\em Journal of Computational Physics}, 478, 2023.

\bibitem{CPZMR2020}
H.~Carrillo, E.~Macca, C.~Parés, G.~Russo, and D.~Zorío.
\newblock An order-adaptive {C}ompact {A}pproximate {T}aylor method for systems
  of conservation law.
\newblock {\em Journal of Computational Physics}, 438:31, 2021.

\bibitem{Carrillo-Pares}
H.~Carrillo and C.~Parés.
\newblock Compact approximate taylor methods for systems of conservation laws.
\newblock {\em J. Sci. Comput.}, 80:1832--1866, 2019.

\bibitem{Ciarlet}
P.G. Ciarlet.
\newblock Discrete maximum principle for finite-difference operators.
\newblock {\em Aeq. Math.}, 4:338--352, 1970.

\bibitem{CDL1}
S.~Clain, S.~Diot, and R.~Loub{\`e}re.
\newblock A high-order finite volume method for systems of conservation laws --
  multi-dimensional optimal order detection ({MOOD}).
\newblock {\em Journal of Computational Scientific}, 230(10):4028 -- 4050,
  2011.

\bibitem{CDL0_FVCA}
S.~Clain, S.~Diot, and R.~Loub\`ere.
\newblock Multi-dimensional optimal order detection (mood) — a very
  high-order finite volume scheme for conservation laws on unstructured meshes.
\newblock In Fort F{\"u}rst Halama Herbin~Hubert (Eds.), editor, {\em FVCA 6,
  International Symposium, Prague, June 6-10}, volume~4 of {\em Series:
  Springer Proceedings in Mathematics}, 2011.
\newblock 1st Edition. XVII, 1065 p. 106 illus. in color.

\bibitem{CFL}
R.~Courant, K.~Friedrichs, and H.~Lewy.
\newblock Über die partiellen differenzengleichungen der mathematischen
  physik.
\newblock {\em Mathematische Annalen}, 100(1):32--74, 1928.

\bibitem{CDL2}
S.~Diot, S.~Clain, and R.~Loub{\`e}re.
\newblock Improved detection criteria for the multi-dimensional optimal order
  detection ({MOOD}) on unstructured meshes with very high-order polynomials.
\newblock {\em Computers and Fluids}, 64:43 -- 63, 2012.

\bibitem{CDL3}
S.~Diot, R.~Loub{\`e}re, and S.~Clain.
\newblock The {MOOD} method in the three-dimensional case: Very-high-order
  finite volume method for hyperbolic systems.
\newblock {\em International Journal of Numerical Methods in Fluids},
  73:362--392, 2013.

\bibitem{jin1995runge}
Shi Jin.
\newblock {R}unge-{K}utta methods for hyperbolic conservation laws with stiff
  relaxation terms.
\newblock {\em Journal of Computational Physics}, 122(1):51--67, 1995.

\bibitem{jin2010asymptotic}
Shi Jin.
\newblock Asymptotic preserving (ap) schemes for multiscale kinetic and
  hyperbolic equations: a review.
\newblock {\em Lecture notes for summer school on methods and models of kinetic
  theory (M\&MKT), Porto Ercole (Grosseto, Italy)}, pages 177--216, 2010.

\bibitem{LaxWendroff}
P.~Lax and B.~Wendroff.
\newblock Systems of conservation laws.
\newblock {\em Communications Pure and Applied Mathematics}, 13(2):217--237,
  1960.

\bibitem{LeVeque2007}
R.J. LeVeque.
\newblock {\em Finite difference methods for ordinary and partial differential
  equations: steady-state and time-dependent problems (Classics in Applied
  Mathematics)}.
\newblock Society for Industrial and Applied Mathematics, Philadelpia, PA.
  USA., 1 edition, 2007.

\bibitem{CATMOOD_1D}
R.~Loub{\`e}re, E.~Macca, C.~Par{\'e}s, and G.~Russo.
\newblock {CAT-MOOD} methods for conservation laws in one space dimension.
\newblock In C.~Parés, M.J.Castro, M.L. Muñoz, and T.~Morales de~Luna,
  editors, {\em Theory, Numerics and Applications of Hyperbolic Problems},
  SEMA-SIMAI Springer Series, 2023.
\newblock Proceedings of HYP2022.

\bibitem{TesiPhD}
E.~Macca.
\newblock {\em Shock-{C}apturing methods: {W}ell-{B}alanced {A}pproximate
  {T}aylor and {S}emi-{I}mplicit schemes}.
\newblock PhD thesis, Università degli Studi di Palermo, Palermo, 2022.

\bibitem{CATMOOD_2}
E.~Macca, R.~Loubere, C.~Pares, and G.~Russo.
\newblock An almost fail-safe a-posteriori limited high-order {CAT} scheme.
\newblock {\em Journal of Computational Physics, Submitted}, 2023.

\bibitem{MaccaRussoBumi}
E.~Macca and G.~Russo.
\newblock Boundary effects on wave trains in the {E}xner model of sedimental
  transport.
\newblock {\em Bollettino {U}nione {M}atematica {I}taliana}, 2023.

\bibitem{MacCormack}
R.W. MacCormack.
\newblock The effect of viscosity in hypervelocity impact cratering.
\newblock {\em AIAA Pape, Cincinnati, Ohio}, pages 69--354, 1969.

\bibitem{PareschiRusso}
L.~Pareschi and G.~Russo.
\newblock Implicit–{E}xplicit {R}unge–{K}utta {S}chemes and {A}pplications
  to {H}yperbolic {S}ystems with {R}elaxation.
\newblock {\em Journal of Scientific Computing}, 25:129--155, 2005.

\bibitem{Qiu-Shu}
J.~Qiu and C.~W. Shu.
\newblock Finite difference weno schemes with lax-wendroff-type time
  discretizations.
\newblock {\em SIAM J. Sci. Comput.}, 24(6):2185--2198, 2003.

\bibitem{Richtmeyer}
R.~D. Richtmyer and K.~W. Morton.
\newblock Difference methods for initial-value problems.
\newblock {\em Interscience Tracts in Pure and Appl. Math. Interscience, New
  York}, 1, 1967.

\bibitem{Schwartzkopff}
T.~Schwartzkopff, C.~Munz, and E.~Toro.
\newblock Ader: a high-order approach for linear hyperbolic systems in 2d.
\newblock {\em J. Sci. Comput.}, 17:231--240, 2002.

\bibitem{Shu1997}
C.~W. Shu.
\newblock Essentially non-oscillatory and weighted essentially non--oscillatory
  schemes for hyperbolic conservation laws.
\newblock Technical report, Institute for Computer Applications in Science and
  Engineering (ICASE), 1997.

\bibitem{Titarev}
V.~Titarev and E.~Toro.
\newblock Ader: arbitrary high order godunov approach.
\newblock {\em J. Sci. Comput.}, 17:609--618, 2002.

\bibitem{Toro2009}
E.F. Toro.
\newblock {\em {{Riemann}} Solvers and Numerical Methods for Fluid Dynamics}.
\newblock Springer, third edition, 2009.

\bibitem{wanner1996solving}
G.~Wanner and E.~Hairer.
\newblock {\em Solving ordinary differential equations II: Stiff and
  Differential-Algebraic Problems}.
\newblock Springer Berlin Heidelberg, 1996.

\bibitem{Zorio}
D.~Zorío, A.~Baeza, and P.~Mulet.
\newblock An approximate lax-wendroff-type procedure for high order accurate
  scheme for hyperbolic conservation laws.
\newblock {\em J. Sci. Comput.}, 71(1):246--273, 2017.

\end{thebibliography}

\end{document}